\documentclass[11pt]{article}
\usepackage{amsmath,amssymb,color}

\definecolor{mycolorred}{rgb}{1, 0, 0}

\def\<{\langle}
\def\>{\rangle}
\def\e{{\mathbf{e}}}

\newtheorem{theorem}{Theorem}[section]

\newtheorem{corollary}[theorem]{Corollary}

\newtheorem{lemma}[theorem]{Lemma}

\newtheorem{proposition}[theorem]{Proposition}
\newtheorem{remark}[theorem]{Remark}

\numberwithin{equation}{section}

\begin{document}

\title{Convergence and regularity of probability laws\\
by using an interpolation method}
\author{ \textsc{Vlad Bally}\thanks{%
Universit\'e Paris-Est, LAMA (UMR CNRS, UPEMLV, UPEC), INRIA, F-77454
Marne-la-Vall\'ee, France. Email: \texttt{bally@univ-mlv.fr} }\smallskip \\
%EndAName
\textsc{Lucia Caramellino}\thanks{%
Dipartimento di Matematica, Universit\`a di Roma - Tor Vergata, Via della
Ricerca Scientifica 1, I-00133 Roma, Italy. Email: \texttt{%
caramell@mat.uniroma2.it}}\smallskip\\
}
\date{}
\maketitle

\begin{abstract}
In \cite{bib:[FP]} Fournier and Printems establish a methodology which
allows to prove the absolute continuity of the law of the solution of some
stochastic equations with H\"{o}lder continuous coefficients. This is of
course out of reach by using already classical probabilistic methods based
on Malliavin calculus. In \cite{DR} Debussche and Romito employ some Besov
space technics in order to substantially improve the result of Fournier and
Printems. In our paper we show that this kind of problem naturally fits in
the framework of interpolation spaces: we prove an interpolation inequality
(see Proposition \ref{8bis}) which allows to state (and even to slightly
improve) the above absolute continuity result. Moreover it turns out that
the above interpolation inequality has applications in a completely
different framework: we use it in order to estimate the error in total
variance distance in some convergence theorems.
\end{abstract}

%\tableofcontents

\parindent 0pt

\medskip


\medskip


\section{Introduction}

In this paper we prove an interpolation type inequality which leads to three
main applications. First we give a criteria for the regularity of the law $%
\mu $ of a random variable. This was the first aim of the integration by
parts formulas constructed in the Malliavin calculus (in the Gaussian
framework, and of many other variants of this calculus, in a more general
case). But our starting point was the paper of N. Fournier and J. Printems
\cite{bib:[FP]} who noticed that some regularity of the law may be obtained
even if no integration by parts formula holds for $\mu $ itself: they just
use a sequence $\mu _{n}\rightarrow \mu $ and assume that an integration by
parts formula of type $\int f^{\prime }d\mu _{n}=\int fh_{n}d\mu _{n}$ holds
for each $\mu _{n}.$ If $\sup_{n}\int \left\vert h_{n}\right\vert d\mu
_{n}<\infty $ we are close to Malliavin calculus. But the interesting point
is that one may obtain some regularity for $\mu $ even if $\sup_{n}\int
\left\vert h_{n}\right\vert d\mu _{n}=\infty $ - so we are out of the domain
of application of Malliavin calculus. The key point is that one establishes
an equilibrium between the speed of convergence of $\mu _{n}\rightarrow \mu $
and the blow up $\int \left\vert h_{n}\right\vert d\mu _{n}\uparrow \infty .$
The approach of Fournier and Printems is based on Fourier transforms and
more recently Debussche and Romito \cite{DR} obtained a much more powerful
version of this type of criteria based on Besov space technics. This
methodology has been used in several recent papers (see \cite{bib:[BCl1]},
\cite{bib:[BCl2]}, \cite{bib:[BF]}, \cite{bib:[De]}, \cite{DF} and \cite{F1}%
) in order to obtain the absolute continuity of the law of the solution of
some stochastic equations with weak regularity assumptions on the
coefficients: as a typical example, one proves that, under uniform
ellipticity conditions, diffusion processes with H\"{o}lder continuous
coefficients have absolute continuous law at any time $t>0$. In the present
paper we use a different approach, based on an interpolation argument and on
Orlicz spaces, which allows one to go further and to treat, for example,
diffusion processes with log-H\"{o}lder coefficients.

The second application concerns the regularity of the density with respect
to a parameter. We illustrate this direction by giving sufficient conditions
in order that $(x,y)\rightarrow p_{t}(x,y)$ is smooth with respect to $(x,y)$
where $p_{t}(x,y)$ is the density of the law of $X_{t}(x)$ which is a
piecewise deterministic Markov process starting from $x.$

The third application concerns estimates of the speed of convergence $\mu
_{n}\rightarrow \mu $ in total variation distance, and under some stronger
assumptions, the speed of convergence of the derivatives of the densities of
$\mu _{n}$ to the corresponding derivative of the density of $\mu .$ Such
results appear in a natural way as soon as the suited interpolation
framework is settled.

Let us give our main results. We work with the following weighted Sobolev
norms on $C^{\infty }({\mathbb{R}}^{d};{\mathbb{R}})$:
\begin{equation*}
\left\Vert f\right\Vert _{k,m,p}=\sum_{0\leq \left\vert \alpha \right\vert
\leq k}\Big(\int (1+\left\vert x\right\vert )^{m}\left\vert \partial
_{\alpha }f(x)\right\vert ^{p}dx\Big)^{1/p},\qquad p>1,
\end{equation*}%
where $\alpha $ is a multi index, $|\alpha |$ denotes its length and $%
\partial _{\alpha }$ is the corresponding derivative. In the case $m=0$ we
have the standard Sobolev norm that we denote by $\left\Vert f\right\Vert
_{k,p}.$ We will also consider the weaker norm

\begin{equation*}
\left\Vert f\right\Vert _{k,m,1+}=\sum_{0\leq \left\vert \alpha \right\vert
\leq k}\int (1+\left\vert x\right\vert )^{m}\left\vert \partial _{\alpha
}f(x)\right\vert (1+\ln ^{+}\left\vert x\right\vert +\ln ^{+}\left\vert
f(x)\right\vert )dx,
\end{equation*}%
with $\ln ^{+}(x)=\max \{0,\ln \left\vert x\right\vert \}$. Moreover, for
two measures $\mu $ and $\nu $ we consider the distances
\begin{equation*}
d_{k}(\mu ,\nu )=\sup \Big\{\Big\vert\int fd\mu -\int fd\nu \Big\vert%
:\sum_{0\leq |\alpha |\leq k}\Vert \partial _{\alpha }f\Vert _{\infty }\leq 1%
\Big\}.
\end{equation*}%
For $k=0$ this is the total variation distance and for $k=1$ this is the
Fortet Mourier distance.

Our key estimate is the following. Let $m,q,k\in {\mathbb{N}}$ and $p>1$ be
given and let $p_{\ast }$ be the conjugate of $p.$\ We consider a function $%
f\in C^{q+2m}({\mathbb{R}}^{d})$ and a sequence of functions $f_{n}\in
C^{q+2m}({\mathbb{R}}^{d}),n\in {\mathbb{N}} $ and we denote $\mu
(dx)=f(x)dx $ and $\mu _{n}(dx)=f_{n}(x)dx.$ We prove that there exists a
universal constant $C$ such that%
\begin{equation}
\left\Vert f\right\Vert _{q,p}\leq C\Big(\sum_{n=0}^{\infty
}2^{n(q+k+d/p_{\ast })}d_{k}(\mu ,\mu _{n})+\sum_{n=0}^{\infty }\frac{1}{%
2^{2mn}}\left\Vert f_{n}\right\Vert _{q+2m,2m,p}\Big)  \label{ii1}
\end{equation}%
and
\begin{equation}
\left\Vert f\right\Vert _{q,1+}\leq C\Big(\sum_{n=0}^{\infty
}n2^{n(q+k)}d_{k}(\mu ,\mu _{n})+\sum_{n=0}^{\infty }\frac{1}{2^{2mn}}%
\left\Vert f_{n}\right\Vert _{q+2m,2m,,1+}\Big)  \label{ii1'}
\end{equation}

This is Proposition \ref{8bis} and the proof is based on a development in
Hermite series and on a powerful estimate for mixtures of Hermite kernels
inspired from \cite{bib:[PY]}. This inequality fits in the general theory of
interpolation spaces (we thank to D. Elworthy for a useful remark in this
sense). Many interpolation results between Sobolev spaces of positive and
negative indexes are known but they are not relevant from a probabilistic
point of view: convergence in distribution is characterized by the Fortet
Mourier distance and this amounts to convergence in the dual of $W^{1,\infty
}.$ So we are not concerned with Sobolev spaces associated to $L^{p}$ norms
but to $L^{\infty }$ norms. This is a limit case which is more delicate and
we have not found in the literature classical interpolation results which
may be used in our framework.

Once we have (\ref{ii1}) and(\ref{ii1'}) we obtain the following regularity
criteria. Let $\mu $ be a finite non negative measure. Suppose that there
exists a sequence of functions $f_{n}\in C^{q+2m}({\mathbb{R}}^{d}),n\in {%
\mathbb{N}}$ such that
\begin{equation}
d_{k}(\mu ,\mu _{n})\times \left\Vert f_{n}\right\Vert
_{1+q+2m,2m,p}^{\alpha }\leq C,\qquad \alpha >\frac{q+k+d/p_{\ast }}{2m}.
\label{ii2}
\end{equation}%
with $\mu _{n}(dx)=f_{n}(x)dx.$ Then $\mu (dx)=f(x)dx$ and $f\in W^{q,p}$
(the standard Sobolev space).

In terms of $\left\Vert f\right\Vert _{q,m,,1+}$ the statement is the
following: suppose that there exists $m\in {\mathbb{N}}$ such that
\begin{equation}
d_{1}(\mu ,\mu _{n})\times \left\Vert f_{n}\right\Vert
_{2m,2m,1+}^{1/2m}\leq \frac{C}{(\ln n)^{2+1/2m}}.  \label{ii2'}
\end{equation}%
Then $\mu $ is absolutely continuous with respect to the Lebesgue measure.

The statement of the corresponding results are Theorem \ref{ThLp} and
Theorem \ref{ThLog} respectively. These are two significant particular cases
of a more general result stated in terms of Orlicz norms in Theorem \ref{2C}%
. The proof is, roughly speaking, as follows: let $\gamma _{\varepsilon }$
be the Gaussian density of variance $\varepsilon >0$ and let $\mu
^{\varepsilon }=\mu \ast \gamma _{\varepsilon }$ and $\mu _{n}^{\varepsilon
}=\mu _{n}\ast \gamma _{\varepsilon }.$ Then $\mu ^{\varepsilon
}(dx)=f^{\varepsilon }(x)dx$ and $\mu _{n}^{\varepsilon
}(x)=f_{n}^{\varepsilon }(x)dx.$ Using (\ref{ii1}) for $f^{\varepsilon }$
and $f_{n}^{\varepsilon },n\in {\mathbb{N}}$ one proves that $%
\sup_{\varepsilon }\left\Vert f^{\varepsilon }\right\Vert _{q,p}<\infty .$
And then one employs a relatively compactness argument in $W^{q,p}$ in order
to produce the density $f$ of $\mu .$

We give now the convergence result (see Theorem \ref{Conv}). Suppose that (%
\ref{ii2}) holds for some $\alpha >\frac{q+k+d/p_{\ast }}{m}.$ Then $\mu
(dx)=f(x)dx$ and, for every $n\in {\mathbb{N}},$
\begin{equation}
\left\Vert f-f_{n}\right\Vert _{W^{q,p}}\leq Cd_{k}^{\theta }(\mu ,\mu
_{n})\quad \mbox{with}\quad \theta =\frac{1}{\alpha }\wedge (1-\frac{%
q+k+d/p_{\ast }}{\alpha m}).  \label{ii3}
\end{equation}%
Roughly speaking this inequality is obtained by using (\ref{ii1}) with $\mu $
replaced by $\mu -\mu _{n}.$

In the statements of (\ref{ii2}) we do not use $d_{k}(\mu ,\mu _{n})$ and $%
\left\Vert f_{n}\right\Vert _{1+q+2m,2m,p}$ directly, but some function $%
\lambda $ which have some nice properties and such that $\lambda (1/n)\geq
\left\Vert f_{n}\right\Vert _{1+q+2m,2m,p}$. But this is a technical point
which we leave out in this introduction.

The paper is organized as follows. In Section \ref{sect-results} we
introduce the Orlicz spaces, we give the general result and the criteria
concerning the absolute continuity and the regularity of the density. We
also give in Section \ref{sect:convLp} the convergence criteria mentioned
above. In Section \ref{sect-RV} we translate the results in terms of
integration by parts formulae. In Section \ref{sect-pathdep} (respectively
Section \ref{sect-heat}) we prove absolute continuity for the law of the
solution to a SDE (respectively to a SPDE) with log-H\"{o}lder continuous
coefficients. Moreover, in Section \ref{jumps} we discuss an example
concerning piecewise deterministic Markov processes: we assume that the
coefficients are smooth and we prove existence of the density of the law of
the solution together with regularity with respect to the initial condition.
We also consider an approximation scheme and we use (\ref{ii3}) in order to
estimate the error. Finally, we add some appendices containing technical
results: Appendix \ref{sect-Hermite} is devoted to the proof of the main
estimate (\ref{ii1}) based on a development in Hermite series; in Appendix %
\ref{sect-interp} we discuss the relation with interpolation spaces; in
Appendix \ref{app-superkernels} we give some auxiliary estimates concerning
super-kernels.

\section{Criterion for the regularity of a probability law}

\label{sect-results}

\subsection{Notations}

We work on ${\mathbb{R}}^{d}$ and we denote by $\mathcal{M}$ the set of the
finite signed measures on ${\mathbb{R}}^{d}$ with the Borel $\sigma $
algebra. Moreover $\mathcal{M}_{a}\mathcal{\subset M}$ is the set of the
measures which are absolutely continuous with respect to the Lebesgue
measure. For $\mu \in \mathcal{M}_{a}$ we denote by $p_{\mu }$ the density
of $\mu $ with respect to the Lebesgue measure. And for a measure $\mu \in
\mathcal{M}$ we denote by $L_{\mu }^{p}$\ the space of the measurable
functions $f:{\mathbb{R}}^{d}\rightarrow {\mathbb{R}}$ such that $\int
\left\vert f\right\vert ^{p}d\left\vert \mu \right\vert <\infty .$ For $f\in
L_{\mu }^{1}$ we denote $f\mu $ the measure $(f\mu )(A)=\int_{A}fd\mu .$ For
a bounded function $\phi :{\mathbb{R}}^{d}\rightarrow {\mathbb{R}}$ we
denote $\mu \ast \phi $ the measure defined by $\int fd\mu \ast \phi =\int
f\ast \phi d\mu =\int \int \phi (x-y)f(y)dyd\mu (x).$ Then $\mu \ast \phi
\in \mathcal{M}_{a}$ and $p_{\mu \ast \phi }(x)=\int \phi (x-y)d\mu (y).$

We denote by $\alpha =(\alpha _{1},...,\alpha _{d})\in {\mathbb{N}}^{d}$ a
multi index and we put $\left\vert \alpha \right\vert =\sum_{i=1}^{d}\alpha
_{i}.$ Here ${\mathbb{N}}=\{0,1,2,...\}$ is the set of non negative integers
and we put ${\mathbb{N}}_{\ast }={\mathbb{N}}\setminus \{0\}.$ For a multi
index $\alpha $ with $\left\vert \alpha \right\vert =k$ we denote $\partial
_{\alpha }$ the corresponding derivative that is $\partial _{x_{1}}^{\alpha
_{1}}...\partial _{x_{d}}^{\alpha _{d}}$ with the convention that $\partial
_{x_{i}}^{\alpha _{i}}f=f$ if $\alpha _{i}=0.$ In particular if $\alpha $ is
the null multi index then $\partial _{\alpha }f=f.$

We denote by $\left\Vert f\right\Vert _{p}=(\int \left\vert f(x)\right\vert
^{p}dx)^{1/p},p\geq 1$ and $\left\Vert f\right\Vert _{\infty }=\sup_{x\in {%
\mathbb{R}}^{d}}\left\vert f(x)\right\vert .$ Then $L^{p}=\{f:\left\Vert
f\right\Vert _{p}<\infty \}$ are the standard $L^{p}$ spaces with respect to
the Lebesgue measure.

\subsection{Orlicz spaces}

\label{sect-orlicz} In the following we will work in Orlicz spaces, so we
briefly recall the notation and the results we will use, for which we refer
to \cite{bib:[K.M-W]}.

A function $\mathbf{e}:{\mathbb{R}}\rightarrow {\mathbb{R}}_{+}$ is a Young
function if it is symmetric, strictly convex, non negative and $\mathbf{e}%
(0)=0.$ In the following we will consider Young functions having the two
supplementary properties:%
\begin{equation}
\begin{array}{ll}
i)\quad & \text{there exists }\lambda >0\text{ such that }\mathbf{e}(2s)\leq
\lambda \mathbf{e}(s),\smallskip \\
ii)\quad & s\mapsto \displaystyle\frac{\mathbf{e}(s)}{s}\text{ is non
decreasing.}%
\end{array}
\label{Y}
\end{equation}%
The property $i)$ is known as the $\Delta _{2}$ condition or doubling
condition (see \cite{bib:[K.M-W]}). Through the whole paper we work with
Young functions which satisfy (\ref{Y}). We set $\mathcal{E}$ the space of
these functions:
\begin{equation}
\mathcal{E}=\{\mathbf{e}\,:\,%
\mbox{$\e$ is a Young function satisfying
\eqref{Y}.}\}  \label{Espace}
\end{equation}%
For $\mathbf{e}\in \mathcal{E}$ and $f:{\mathbb{R}}^{d}\rightarrow {\mathbb{R%
}}$, we define the norm
\begin{equation}
\left\Vert f\right\Vert _{\mathbf{e}}=\inf \Big\{c>0:\int \mathbf{e}\Big(%
\frac{1}{c}f(x)\Big)dx\leq 1\Big\}.  \label{O1}
\end{equation}%
This is the so called Luxembourg norm which is equivalent to the Orlicz norm
(see \cite{bib:[K.M-W]} p 227 Th 7.5.4). It is convenient for us to work
with this norm (instead of the Orlicz norm). The space $L^{\mathbf{e}%
}=\{f:\left\Vert f\right\Vert _{\mathbf{e}}<\infty \}$ is the Orlicz space.

\begin{remark}
\label{U}Let $u_{l}(x)=(1+\left\vert x\right\vert )^{-l}.$ As a consequence
of (\ref{Y}) ii), for every $\mathbf{e}\in \mathcal{E}$ and $l>d$ one has $%
u_{l}\in L^{\mathbf{e}}$ and moreover,
\begin{equation*}
\left\Vert u_{l}\right\Vert _{\mathbf{e}}\leq (\mathbf{e}(1)\left\Vert
u_{l}\right\Vert _{1})\vee 1<\infty .
\end{equation*}%
Indeed (\ref{Y}) ii) implies that for $t\leq 1$ one has $\mathbf{e}(t)\leq
\mathbf{e}(1)t.$ For $c\geq (\mathbf{e}(1)\left\Vert u_{l}\right\Vert
_{1})\vee 1$ one has $\frac{1}{c}u_{l}(x)\leq u_{l}(x)\leq 1$ so that
\begin{equation*}
\int \mathbf{e}\Big(\frac{1}{c}u_{l}(x)\Big)dx\leq \frac{\mathbf{e}(1)}{c}%
\int u_{l}(x)dx=\frac{\mathbf{e}(1)}{c}\left\Vert u_{l}\right\Vert _{1}\leq
1.
\end{equation*}
\end{remark}

For $a>0$, we define $\mathbf{e}^{-1}(a)=\sup \{c:\mathbf{e}(c)\leq a\}$
and:
\begin{equation}
\phi _{\mathbf{e}}(r)=\frac{1}{\mathbf{e}^{-1}\big(\frac{1}{r}\big)}\quad %
\mbox{and}\quad \beta _{\mathbf{e}}(R)=\frac{R}{\mathbf{e}^{-1}(R)}=R\phi _{%
\mathbf{e}}\Big(\frac{1}{R}\Big),\quad r,R>0.  \label{O3}
\end{equation}

\begin{remark}
\label{Increasing} The function $\phi _{\mathbf{e}}$ is the
\textquotedblleft fundamental function\textquotedblright\ of $L^{\mathbf{e}}$
equipped with the Luxembourg norm (see \cite{bib:[BS]} Lemma 8.17 pg 276).
In particular $\frac{1}{r}\phi _{\mathbf{e}}(r)$ is decreasing (see \cite%
{bib:[BS]} Corollary 5.2 pg 67). It follows that $\beta _{\mathbf{e}}$ is
increasing. For the sake of completeness we give here the argument. By (\ref%
{Y})$,ii)$, if $a>1$ then $\mathbf{e}(ax)\geq a\mathbf{e}(x)$ so that $%
ax\geq \mathbf{e}^{-1}(a\mathbf{e}(x)).$ Taking $y=\mathbf{e}(x)$ we obtain $%
a\mathbf{e}^{-1}(y)\geq \mathbf{e}^{-1}(ay)$ which gives
\begin{equation*}
\beta _{\mathbf{e}}(ay)=\frac{ay}{\mathbf{e}^{-1}(ay)}\geq \frac{ay}{a%
\mathbf{e}^{-1}(y)}=\beta _{\mathbf{e}}(y).
\end{equation*}
\end{remark}

One defines the conjugate of $\mathbf{e}$ by
\begin{equation*}
\mathbf{e}_{\ast }(s)=\sup \{st-\mathbf{e}(t):t\in {\mathbb{R}}\}.
\end{equation*}
$\mathbf{e}_{\ast }$ is a Young function as well, so the corresponding
Luxembourg norm $\left\Vert f\right\Vert _{\mathbf{e}_{\ast }}$ is given by (%
\ref{O1}) with $\mathbf{e}$ replaced by $\mathbf{e}_{\ast }$. And one has
the following H\"{o}lder inequality:%
\begin{equation}
\left\vert \int fg(x)dx\right\vert \leq 2\left\Vert f\right\Vert _{\mathbf{e}%
}\left\Vert g\right\Vert _{\mathbf{e}_{\ast }}.  \label{O2}
\end{equation}%
(see Theorem 7.2.1 at p 215 in \cite{bib:[K.M-W]}; we stress that the factor
$2$ does not appear in that reference but in the right hand side of the
inequality in the statement of Theorem 7.2.1 in \cite{bib:[K.M-W]} one has
the Orlicz norm of $g$ and by using the equivalence between the Orlicz and
the Luxembourg norm we can replace the Orlicz norm by $2\left\Vert
g\right\Vert _{\mathbf{e}_{\ast }}$).

If $\mathbf{e}$ satisfies the $\Delta _{2}$ condition (that is (\ref{Y}) $i)$%
) then $L^{\mathbf{e}}$ is reflexive (see \cite{bib:[K.M-W]}, Theorem 7.7.1,
p 234). In particular, in this case, any bounded subset of $L^{\mathbf{e}}$
is weakly relatively compact.

For $f\,\in C^{\infty }\,({\mathbb{R}}^{d},{\mathbb{R)}}$, we introduce the
norms%
\begin{equation}
\left\Vert f\right\Vert _{k,\mathbf{e}}=\sum_{0\leq \left\vert \alpha
\right\vert \leq k}\left\Vert \partial _{\alpha }f\right\Vert _{\mathbf{e}%
}\quad \mbox{and}\quad \left\Vert f\right\Vert _{k,\infty }=\sum_{0\leq
\left\vert \alpha \right\vert \leq k}\left\Vert \partial _{\alpha
}f\right\Vert _{\infty }  \label{O4}
\end{equation}%
and we denote
\begin{equation*}
W^{k,\mathbf{e}}=\{f:\left\Vert f\right\Vert _{k,\mathbf{e}}<\infty \}\qquad %
\mbox{and}\qquad W^{k,\infty }=\{f:\left\Vert f\right\Vert _{k,\infty
}<\infty \}.
\end{equation*}%
For a multi index $\gamma $ we denote $x^{\gamma
}=\prod_{i=1}^{d}x_{i}^{\gamma _{i}}$ and for two multi indexes $\alpha
,\gamma $ we denote $f_{\alpha ,\gamma }$ the function
\begin{equation*}
f_{\alpha ,\gamma }(x)=x^{\gamma }\partial _{\alpha }f(x).
\end{equation*}%
Then we consider the norms
\begin{equation}
\left\Vert f\right\Vert _{k,l,\mathbf{e}}=\sum_{0\leq \left\vert \alpha
\right\vert \leq k}\sum_{0\leq \left\vert \gamma \right\vert \leq
l}\left\Vert f_{\alpha ,\gamma }\right\Vert _{\mathbf{e}}\quad \mbox{and}%
\quad W^{k,l,\mathbf{e}}=\{f:\left\Vert f\right\Vert _{k,l,\mathbf{e}%
}<\infty \}.  \label{O4a}
\end{equation}

We stress that in $\Vert \cdot \Vert _{k,l,\mathbf{e}}$ the first index $k$
is related to the order of the derivatives which are involved while the
second index $l$ is connected to the power of the polynomial multiplying the
function and its derivatives up to order $k$.

Let us propose two examples of Young functions, that represent the leading
ones in our approach.

\medskip

\textbf{Example 1.} If we take $\mathbf{e}_{p}(x)=\left\vert x\right\vert
^{p},p>1,$ then $\left\Vert f\right\Vert _{\mathbf{e}_{p}}$ is the usual $%
L^{p}$ norm and the corresponding Orlicz space is the standard $L^{p}$ space
on ${\mathbb{R}}^{d}$. Clearly $\beta _{{\mathbf{e}}_p}(t)=t^{1/p_{\ast }}$
with $p_{\ast }$ the conjugate of $p.$

\medskip

\textbf{Example 2.} %\begin{example}\label{ex-2}
Set $\mathbf{e}_{\log }(t)=(1+\left\vert t\right\vert )\ln (1+\left\vert
t\right\vert ).$

\medskip

Since the norm from ${\mathbf{e}}_{\log }$ is not explicit we replace it by
the following quantities:%
\begin{equation}
\begin{array}{l}
\displaystyle\left\Vert f\right\Vert _{p,1+}=\int (1+\left\vert x\right\vert
)^{p}\left\vert f(x)\right\vert (1+\ln ^{+}|x|+\ln ^{+}\left\vert
f(x)\right\vert )dx\smallskip \\
\displaystyle\left\Vert f\right\Vert _{k,p,1+}=\sum_{0\leq \left\vert \alpha
\right\vert \leq k}\left\Vert \partial _{\alpha }f\right\Vert _{p,1+}%
\end{array}
\label{E22}
\end{equation}%
with $\ln ^{+}(x)=\max \{0,\ln \left\vert x\right\vert \}.$ We stress that $%
\left\Vert f\right\Vert _{p,1+}$ is not a norm.

We will need the following:

\begin{lemma}
For each $k\in {\mathbb{N}}$ and $p\geq 0$ there exists a constant $C$
depending on $k,p$ only such that
\begin{equation}
\left\Vert f\right\Vert _{k,p,{\mathbf{e}}_{\log }}\leq C(1\vee \left\Vert
f\right\Vert _{k,p,1+}).  \label{E23}
\end{equation}%
Moreover%
\begin{equation}
\limsup_{t\rightarrow \infty }\frac{\beta _{{\mathbf{e}}_{\log }}(t)}{\ln t}%
\leq 2.  \label{E1}
\end{equation}
\end{lemma}

\textbf{Proof}. The inequality (\ref{E23}) is an immediate consequence of
the following simpler one:
\begin{equation}
\left\Vert f\right\Vert _{{\mathbf{e}}_{\log }}\leq 2\Big(1\vee \int
\left\vert f(x)\right\vert (1+\ln ^{+}\left\vert f(x)\right\vert )dx\Big).
\label{E24}
\end{equation}%
Let us prove it. We assume that $f\geq 0$ and we take $c\geq 2$ and we write%
\begin{equation*}
\int {\mathbf{e}}_{\log }\Big(\frac{1}{c}f(x)\Big)dx\leq \int_{\{f\leq c\}}{%
\mathbf{e}}_{\log }\Big(\frac{1}{c}f(x)\Big)dx+\int_{\{f>c\}}{\mathbf{e}}%
_{\log }\Big(\frac{1}{c}f(x)\Big)dx=:I+J.
\end{equation*}%
Using the inequality $\ln (1+y)\leq y$ we obtain $I\leq 2\int \ln (1+\frac{1%
}{c}f)\leq \frac{2}{c}\int f.$ And if $f\geq c\geq 2$ then $\frac{f}{c}%
+1\leq \frac{2}{c}f\leq f.$ Then ${\mathbf{e}}_{\log }(\frac{1}{c}f(x))\leq
\frac{2}{c}f\ln f.$ It follows that $J\leq \frac{2}{c}\int_{\{f>c\}}f\ln
^{+}f$ and finally $\int {\mathbf{e}}_{\log }(\frac{1}{c}f))\leq \frac{2}{c}%
\int_{\{f>c\}}(1+f)\ln ^{+}f.$ We conclude that for $c\geq 2\int f(1+\ln
^{+}f)$ we have $\int {\mathbf{e}}_{\log }(\frac{1}{c}f)\leq 1$ which by the
very definition means that $\left\Vert f\right\Vert _{{\mathbf{e}}_{\log
}}\leq 2\int f(1+\ln ^{+}f)$.

Let us prove (\ref{E1}). We denote $e(t)=2t\ln (2t)$ and we notice that for
large $t$ one has ${\mathbf{e}}_{\log }(t)\leq e(t).$ It follows that%
\begin{equation*}
\beta _{{\mathbf{e}}_{\log }}(t)\leq \frac{t}{e^{-1}(t)}.
\end{equation*}%
Using the change of variable $R=e(t)$ we obtain
\begin{equation*}
\lim_{R\rightarrow \infty }\frac{R}{e^{-1}(R)\ln R}=\lim_{t\rightarrow
\infty }\frac{e(t)}{t\ln e(t)}=2.
\end{equation*}%
So for large $R$ we have $\beta _{{\mathbf{e}}_{\log }}(R)\leq
R/e^{-1}(R)\leq 2\ln R.$ $\square $

\begin{remark}
We recall that the $LlogL$ space of Zygmund is the space of the functions $f
$ such that $\int \left\vert f(x)\right\vert \ln ^{+}\left\vert
f(x)\right\vert dx<\infty $ (see \cite{bib:[BS]}). Then $L^{{\mathbf{e}}%
_{\log }}=L^{1}\cap LlogL.$ The inequality (\ref{E24}) already gives one
inclusion. The converse inclusion is a consequence of the following
inequalities. Let $\varepsilon _{\ast }>0$ be such that $t\leq 2\ln (1+t)$
for $0<t\leq \varepsilon _{\ast }$ and let $C_{\ast }=2+1/\ln (1+\varepsilon
_{\ast }).$ Then%
\begin{equation}  \label{E24bis}
\begin{array}{ll}
i) & \qquad \displaystyle \int \left\vert f(x)\right\vert dx \leq C_{\ast
}\left\Vert f\right\Vert _{{\mathbf{e}}_{\log }}\quad \mbox{and}\smallskip
\\
ii) & \qquad \displaystyle \int \left\vert f(x)\right\vert \ln
^{+}\left\vert f(x)\right\vert dx \leq \left\Vert f\right\Vert _{{\mathbf{e}}%
_{\log }}(1+2C_{\ast }\ln ^{+}\left\Vert f\right\Vert _{{\mathbf{e}}_{\log
}}).%
\end{array}%
\end{equation}%
In order to prove $i)$ we denote $g=\left\Vert f\right\Vert _{{\mathbf{e}}%
_{\log }}^{-1}\left\vert f\right\vert $ and we write
\begin{align*}
\int g &=\int_{\{g\leq \varepsilon _{\ast }\}}g+\int_{\{g>\varepsilon _{\ast
}\}}g\leq 2\int_{\{g\leq \varepsilon _{\ast }\}}\ln (1+g)+\frac{1}{\ln
(1+\varepsilon _{\ast })}\int_{\{g>\varepsilon _{\ast }\}}g\ln (1+g) \\
&\leq C_{\ast }\int (1+g)\ln (1+g)=C_{\ast }\int {\mathbf{e}}_{\log
}(g)=C_{\ast }.
\end{align*}%
In order to prove $ii)$ we notice that $\int g\ln ^{+}g\leq \int {\mathbf{e}}%
_{\log }(g)=1$ so that
\begin{equation*}
\int \left\vert f\right\vert \ln ^{+}\frac{\left\vert f\right\vert }{%
\left\Vert f\right\Vert _{{\mathbf{e}}_{\log }}}\leq \left\Vert f\right\Vert
_{{\mathbf{e}}_{\log }}.
\end{equation*}%
Then we write%
\begin{equation*}
\int \left\vert f\right\vert \ln ^{+}\left\vert f\right\vert
=\int_{\{\left\vert f\right\vert \geq 1\vee \left\Vert f\right\Vert _{{%
\mathbf{e}}_{\log }}\}}\left\vert f\right\vert \ln ^{+}\left\vert
f\right\vert +\int_{\{\left\vert f\right\vert <1\vee \left\Vert f\right\Vert
_{{\mathbf{e}}_{\log }}\}}\left\vert f\right\vert \ln ^{+}\left\vert
f\right\vert =:I+J.
\end{equation*}%
If $\left\vert f\right\vert \geq 1\vee \left\Vert f\right\Vert _{{\mathbf{e}}%
_{\log }}$ then $\ln ^{+}| f|=\ln |f|=\ln^+(\frac{|f|}{\|f\|_{{\mathbf{e}}%
_{\log}}})+\ln\|f\|_{{\mathbf{e}}_{\log}}$. So, by using the previous
inequality,
\begin{equation*}
I\leq \left\Vert f\right\Vert _{{\mathbf{e}}_{\log }}+\ln \left\Vert
f\right\Vert _{{\mathbf{e}}_{\log }}\int \left\vert f\right\vert \leq
\left\Vert f\right\Vert _{{\mathbf{e}}_{\log }}(1+C_{\ast }\ln \left\Vert
f\right\Vert _{{\mathbf{e}}_{\log }})
\end{equation*}%
the last inequality being a consequence of $i).$ And%
\begin{equation*}
J\leq \ln ^{+}\left\Vert f\right\Vert _{{\mathbf{e}}_{\log }}\int \left\vert
f\right\vert \leq C_{\ast }\left\Vert f\right\Vert _{{\mathbf{e}}_{\log
}}\ln ^{+}\left\Vert f\right\Vert _{{\mathbf{e}}_{\log }}.
\end{equation*}
\end{remark}

\subsection{Main results}

\label{main-res}

We consider the following distances between two measures $\mu ,\nu \in
\mathcal{M}$: for $k\in {\mathbb{N}}$, we set
\begin{equation}
d_{k}(\mu ,\nu )=\sup \Big\{\Big\vert\int \phi d\mu -\int \phi d\nu \Big\vert%
:\phi \in C^{\infty }({\mathbb{R}}^{d}),\left\Vert \phi \right\Vert
_{k,\infty }\leq 1\Big\}.  \label{O6}
\end{equation}%
Notice that $d_{0}$ is the total variation distance and $d_{1}$ is the
bounded variation distance (also called Fort\'{e}t Mourier distance). We
recall that the Wasserstein distance (which is more popular) is $d_{W}(\mu
,\nu )=\sup \{\left\vert \int \phi d\mu -\int \phi d\nu \right\vert :\phi
\in C^{1}({\mathbb{R}}^{d}),\left\Vert \nabla \phi \right\Vert _{\infty
}\leq 1\}$, so that $d_{1}(\mu ,\nu )\leq d_{W}(\mu ,\nu ).$ It follows that
all the results proved with respect to $d_{1}$ will be a fortiori true for $%
d_{W}.$ The Wasserstein distance is relevant from a probabilistic point of
view because it characterizes the convergence in law of probability
measures. The distances $d_{k}$ with $k\geq 2$ are less often used. We
mention however that people working in approximation theory (for diffusion
process for example - see \cite{bib:[TT]} or \cite{bib:[VN]}) use such
distances in an implicit way: indeed, they study the speed of convergence of
certain schemes but they are able to obtain their estimates for test
functions $f\in C^{k}$ with $k$ sufficiently large - so $d_{k}$ comes on. We
also recall that for $k=1,2,3$, $d_k$ plays an important role in the
so-called Stein's method for normal approximation (see e.g. \cite{bib:NP}).

We fix now a Young function $\mathbf{e}\in \mathcal{E}$ (see \eqref{Espace}%
), and we recall the function $\beta _{\mathbf{e}}$ (see \eqref{O3} and
Remark \ref{Increasing} respectively).

Let $q,k\in {\mathbb{N}}$ and $m\in {\mathbb{N}}_{\ast }.$ For $\mu \in
\mathcal{M}$ and for a sequence $\mu _{n}\in \mathcal{M}_{a},n\in {\mathbb{N}%
}$ we define
\begin{equation}
\pi _{q,k,m,\mathbf{e}}(\mu ,(\mu _{n})_{n})=\sum_{n=0}^{\infty
}2^{n(q+k)}\beta _{\mathbf{e}}(2^{nd})d_{k}(\mu ,\mu
_{n})+\sum_{n=0}^{\infty }\frac{1}{2^{2nm}}\left\Vert p_{\mu
_{n}}\right\Vert _{2m+q,2m,\mathbf{e}}.  \label{O7}
\end{equation}%
Moreover we define
\begin{equation}
\rho _{q,k,m,\mathbf{e}}(\mu )=\inf \pi _{q,k,m,\mathbf{e}}(\mu ,(\mu
_{n})_{n})  \label{011}
\end{equation}%
the infimum being over all the sequences of measures $\mu _{n},n\in {\mathbb{%
N}}$ which are absolutely continuous. It is easy to check that $\rho _{q,k,m,%
\mathbf{e}}$ is a norm on the space $\mathcal{S}_{q,k,m,\mathbf{e}}$ defined
by%
\begin{equation}
\mathcal{S}_{q,k,m,\mathbf{e}}=\{\mu \in \mathcal{M}:\rho _{q,k,m,\mathbf{e}%
}(\mu )<\infty \}.  \label{O13}
\end{equation}

\smallskip

The following result gives the key estimate in our paper. We prove it in
Appendix \ref{sect-Hermite}.

\begin{proposition}
\label{8bis}Let $q,k\in {\mathbb{N}},m\in {\mathbb{N}}_{\ast }$ and $\mathbf{%
e}\in \mathcal{E}.$ There exists a universal constant $C$ (depending on $%
q,k,m,d$ and $\mathbf{e}$) such that for every $f\in C^{2m+q}({\mathbb{R}}%
^{d})$ one has
\begin{equation}
\left\Vert f\right\Vert _{q,\mathbf{e}}\leq C\rho _{q,k,m,\mathbf{e}}(\mu )
\label{Oo10bis}
\end{equation}%
where $\mu (dx)=f(x)dx.$
\end{proposition}

We state now our main theorem:

\begin{theorem}
\label{2C} Let $q,k\in {\mathbb{N}},m\in {\mathbb{N}}_{\ast }$ and let $%
\mathbf{e}\in \mathcal{E}$.

$i)$ Take $q=0.$ Then%
\begin{equation*}
\mathcal{S}_{0,k,m,\mathbf{e}}\subset L^{\mathbf{e}}
\end{equation*}%
in the sense that if $\mu \in \mathcal{S}_{0,k,m,\mathbf{e}}$ then $\mu $ is
absolutely continuous and the density $p_{\mu }$ belongs to $L^{\mathbf{e}}.$
Moreover there exists a universal constant $C$ such that
\begin{equation*}
\left\Vert p_{\mu }\right\Vert _{L^{\mathbf{e}}}\leq C\rho _{0,k,m,\mathbf{e}%
}(\mu ).
\end{equation*}

$ii)$ Take $q\geq 1.$ Then
\begin{equation*}
\mathcal{S}_{q,k,m,\mathbf{e}}\subset W^{q,\mathbf{e}}\quad \mbox{and}\quad
\left\Vert p_{\mu }\right\Vert _{q,\mathbf{e}}\leq C\rho _{q,k,m,\mathbf{e}%
}(\mu ),\quad \mu \in \mathcal{S}_{q,k,m,\mathbf{e}}.
\end{equation*}
\end{theorem}

\textbf{Proof.} We consider a function $\phi \in C_{b}^{\infty }({\mathbb{R}}%
^{d})$ such that $0\leq \phi \leq 1_{B_{1}}$ and, for $\delta \in (0,1),$ we
define $\phi _{\delta }(x)=\delta ^{-d}\phi (\delta ^{-1}x).$ For a measure $%
\mu $ we define $\mu \ast \phi _{\delta }$ by $\int fd\mu \ast \phi _{\delta
}=\int f\ast \phi _{\delta }d\mu .$ Since $\left\Vert f\ast \phi _{\delta
}\right\Vert _{k,\infty }\leq \left\Vert f\right\Vert _{k,\infty }$ it
follows that $d_{k}(\mu \ast \phi _{\delta },\nu \ast \phi _{\delta })\leq
d_{k}(\mu ,\nu ).$ We will also prove that%
\begin{equation}
\left\Vert f\ast \phi _{\delta }\right\Vert _{2m+q,2m,{\mathbf{e}}}\leq
2^{2m}\left\Vert f\right\Vert _{2m+q,2m,{\mathbf{e}}}.  \label{1}
\end{equation}%
Suppose for a moment that (\ref{1}) holds. Then%
\begin{equation*}
\pi _{q,k,m,\mathbf{e}}(\mu \ast \phi _{\delta },(\mu _{n}\ast \phi _{\delta
})_{n})\leq 2^{2m}\pi _{q,k,m,\mathbf{e}}(\mu ,(\mu _{n})_{n})\leq
2^{2m}\rho _{q,k,m,{\mathbf{e}}}(\mu ).
\end{equation*}%
Let $p_{\delta }$ be the density of the measure $\mu \ast \phi _{\delta }.$
The above inequality and (\ref{Oo10bis}) prove that
\begin{equation*}
\sup_{0<\delta \leq 1}\left\Vert p_{\delta }\right\Vert _{q,e}\leq C\rho
_{q,k,m,{\mathbf{e}}}(\mu )<\infty .
\end{equation*}%
So the family $p_{\delta },\delta \in (0,1)$ is bounded in $W^{q,\mathbf{e}}$
which is a reflexive space. So it is weakly relatively compact. Consequently
we may find a sequence $\delta _{n}\rightarrow 0$ such that $p_{\delta
_{n}}\rightarrow p$ weakly for some $p\in W^{q,\mathbf{e}}$. Since $\mu \ast
\phi _{\delta }\rightarrow \mu $ weakly $\mu (dx)=p(x)dx.$ And $\left\Vert
p\right\Vert _{q,{\mathbf{e}}}\leq C\rho _{q,k,m,\mathbf{e}}(\mu ).$ So the
proof is completed.

Let us check (\ref{1}). For $\lambda >0$ we denote $g_{\lambda
}(x)=(1+\left\vert x\right\vert )^{\lambda }g(x).$ Notice that for $\delta
\leq 1$%
\begin{eqnarray*}
\left\vert (g\ast \phi _{\delta })_{\lambda }(x)\right\vert &\leq
&(1+\left\vert x\right\vert )^{\lambda }\int \left\vert g(x-y)\right\vert
\phi _{\delta }(y)dy\leq \int (1+\left\vert x-y\right\vert +\delta
)^{\lambda }\left\vert g(x-y)\right\vert \phi _{\delta }(y)dy \\
&\leq &2^{\lambda }\int (1+\left\vert x-y\right\vert )^{\lambda }\left\vert
g(x-y)\right\vert \phi _{\delta }(y)dy=2^{\lambda }\left\vert g_{\lambda
}\right\vert \ast \phi _{\delta }(x).
\end{eqnarray*}%
Then, by (\ref{Oo2}) $\left\Vert (g\ast \phi _{\delta })_{\lambda
}\right\Vert _{\mathbf{e}}\leq 2^{\lambda }\left\Vert \left\vert g_{\lambda
}\right\vert \ast \phi _{\delta }\right\Vert _{\mathbf{e}}\leq 2^{\lambda
}\left\Vert \phi _{\delta }\right\Vert _{1}\left\Vert \left\vert g_{\lambda
}\right\vert \right\Vert _{\mathbf{e}}=2^{\lambda }\left\Vert g_{\lambda
}\right\Vert _{\mathbf{e}}$. Using this inequality (with $\lambda =2m)$ for $%
g=\partial _{\alpha }f$ we obtain (\ref{1}). $\square $

\medskip

We consider now a special class of Orlicz norms which verify a supplementary
condition: given $\alpha ,\gamma \geq 0$ we define%
\begin{equation}  \label{restriction}
\mathcal{E}_{\alpha ,\gamma }=\Big\{{\mathbf{e}}:\limsup_{R\rightarrow
\infty } \frac{\beta _{{\mathbf{e}}}(R)}{R^{\alpha }(\ln R)^{\gamma }}%
<\infty \Big\}.
\end{equation}

In this case we have:

\begin{theorem}
\label{ThLlog} Let $q,k\in {\mathbb{N}},m\in {\mathbb{N}}_{\ast }$ and let $%
\mathbf{e}\in \mathcal{E}_{\alpha ,\gamma }$. If $2m>d$, $\gamma \geq 0$ and
$0\leq \alpha <\frac{2m+q+k}{d(2m-1)}$ then
\begin{equation*}
W^{q+1,2m,{\mathbf{e}}}\subset \mathcal{S}_{q,k,m,\mathbf{e}}\subset W^{q,%
\mathbf{e}}
\end{equation*}%
and there exists some constant $C$ such that%
\begin{equation}
\frac{1}{C}\left\Vert p_{\mu }\right\Vert _{q,{\mathbf{e}}}\leq \rho _{q,k,m,%
\mathbf{e}}(\mu )\leq C\left\Vert p_{\mu }\right\Vert _{q+1,2m,{\mathbf{e}}}.
\label{Norm}
\end{equation}%
In particular this is true for ${\mathbf{e}}_{\log }$ and for ${\mathbf{e}}%
_p $ with $\frac{p-1}{p}<\frac{2m+q+k}{d(2m-1)}.$
\end{theorem}

\textbf{Proof.} The first inequality in (\ref{Norm}) is proved in Theorem (%
\ref{2C}). As for the second, we use Lemma \ref{kernel copy(1)} in Appendix %
\ref{app-superkernels}. Let $f\in W^{q+1,2m,{\mathbf{e}}}$ and $\mu
_{f}(dx)=f(x)dx.$ We have to prove that $\rho _{q,k,m,\mathbf{e}}(\mu
_{f})<\infty .$ We consider a super kernel $\phi $ (see (\ref{kk1})) and we
define $f_{\delta }=f\ast \phi _{\delta }.$ We take $\delta _{n}=2^{-\theta
n}$ with $\theta $ to be chosen in a moment and we choose $n_{\ast }$ such
that for $n\geq n_{\ast }$ one has $\beta _{{\mathbf{e}}}(2^{nd})\leq
C2^{nd\alpha }n^{\gamma }.$ Using (\ref{kk2}) with $l=2m$, we obtain $%
d_{k}(\mu _{f},\mu _{f_{\delta _{n}}})\leq C\left\Vert f\right\Vert _{q+1,2m,%
{\mathbf{e}}}\delta _{n}^{q+k+1}$ and using (\ref{kk3}) we obtain $%
\left\Vert f_{\delta _{n}}\right\Vert _{2m+q,2m,\mathbf{e}}\leq C\left\Vert
f\right\Vert _{q+1,2m,{\mathbf{e}}}\delta _{n}^{2m-1}$. Then we can write
\begin{eqnarray*}
\pi _{q,k,m,\mathbf{e}}(\mu _{f},\mu _{f_{\delta _{n}}}) &=
&\sum_{n=0}^{\infty }2^{n(q+k)}\beta _{\mathbf{e}}(2^{nd})d_{k}(\mu _{f},\mu
_{f_{\delta _{n}}})+\sum_{n=0}^{\infty }\frac{1}{2^{2nm}}\left\Vert
f_{\delta _{n}}\right\Vert _{2m+q,2m,\mathbf{e}} \\
&\leq &C\left\Vert f\right\Vert _{q+1,2m,{\mathbf{e}}}\Big(1+\sum_{n\geq
n_{\ast }}^{\infty }2^{n(q+k+d\alpha -\theta (q+k+1))}n^{\gamma
}+\sum_{n=0}^{\infty }\frac{1}{2^{n(2m-\theta (2m-1))}}\Big).
\end{eqnarray*}%
In order to obtain the convergence of the above series we need to choose $%
\theta $ such that%
\begin{equation*}
\frac{q+k+d\alpha }{q+k+1}<\theta <\frac{2m}{2m-1}
\end{equation*}%
and this is possible under our restriction on $\alpha .$ $\square $

\medskip

We give now a criterion in order to check that $\mu \in \mathcal{S}_{q,k,m,%
\mathbf{e}}.$

\begin{theorem}
\label{Criterion} Let $q,k\in {\mathbb{N}},m\in {\mathbb{N}}_{\ast }$ and
let $\mathbf{e}\in \mathcal{E}_{\alpha ,\gamma }$. We consider a non
negative finite measure $\mu $ and we suppose that there exists a family of
measures $\mu _{\delta }(dx)=f_{\delta }(x)dx,\delta >0$ which verifies the
following assumptions. There exist $C,r>0$ and a function $\lambda
_{q,m}(\delta )$, $\delta \in (0,1)$, which is right-continuous and non
increasing such that
\begin{equation*}
\Vert f_{\delta }\Vert _{2m+q,2m,{\mathbf{e}}}\leq \lambda _{q,m}(\delta
)\leq C\delta ^{-r}.
\end{equation*}%
We consider some $\eta >0$ and $\kappa \geq 0$ and we assume that%
\begin{equation}
\lambda _{q,m}^{\eta }(\delta )d_{k}(\mu ,\mu _{\delta })\leq \frac{C}{(\ln
(1/\delta ))^{\kappa }}.  \label{Balance}
\end{equation}%
If (\ref{Balance}) holds with
\begin{equation}
\eta >\frac{q+k+\alpha d}{2m},\qquad \kappa =0  \label{i2}
\end{equation}%
then
\begin{equation*}
\mu \in \mathcal{S}_{q,k,m,\mathbf{e}}\subset W^{q,\mathbf{e}}.
\end{equation*}%
The same conclusion holds if
\begin{equation}
\eta =\frac{q+k+\alpha d}{2m}\qquad and\qquad \kappa >1+\gamma +\eta .
\label{i3}
\end{equation}
\end{theorem}

\textbf{Proof}. Let $\varepsilon _{0}>0.$ We define
\begin{equation*}
\delta _{n}=\inf \{\delta >0:\lambda _{q,m}(\delta )\leq \frac{2^{2mn}}{%
n^{1+\varepsilon _{0}}}\}.
\end{equation*}%
Let $0<\theta <2m/r$ where $r$ is the one in the growth condition on $%
\lambda _{q,m}.$ Since $\delta^r\lambda_{q,m}(\delta)\leq C$, we have
\begin{equation*}
\lambda _{q,m}(2^{-\theta n})\leq C2^{n\theta r}\leq \frac{2^{2mn}}{%
n^{1+\varepsilon _{0}}}
\end{equation*}%
which means that $\delta _{n}\leq 2^{-\theta n}.$ Since $\mathbf{e}\in
\mathcal{E}_{\alpha ,\gamma }$ we have
\begin{equation*}
\pi _{q,k,m,{\mathbf{e}}}(\mu ,(\mu _{\delta _{n}})_{n})\leq
C\sum_{n=1}^{\infty }2^{n(q+k+\alpha d)}n^{\gamma }d_{k}(\mu ,\mu _{\delta
_{n}})+C\sum_{n=1}^{\infty }2^{-2mn}\left\Vert f_{\delta _{n}}\right\Vert
_{2m+q,2m,\mathbf{e}}.
\end{equation*}%
Since $\lambda _{q,m}$ is right continuous, $\lambda
_{q,m}(\delta_{n})=2^{2mn}n^{-(1+\varepsilon _{0})}$ so $\sum_{n=1}^{\infty }%
\frac{1}{2^{2mn}}\lambda _{q,m}(\delta _{n})<\infty .$

By recalling that $\ln (1/\delta _{n})\geq C\theta n$ and by using (\ref%
{Balance}), we obtain
\begin{eqnarray}
2^{n(q+k+\alpha d)}n^{\gamma }d_{k}(\mu ,\mu _{\delta _{n}}) &\leq
&2^{n(q+k+\alpha d)}\frac{Cn^{\gamma }}{\lambda _{q,m}^{\eta }(\delta
_{n})(\ln (1/\delta _{n}))^{\kappa }}  \label{i1} \\
&\leq &C\times 2^{n(q+k+\alpha d-2m\eta )}n^{\gamma +\eta (1+\varepsilon
_{0})-\kappa }.  \notag
\end{eqnarray}%
If $q+k+\alpha d<2\eta m$ the series with the general term given in (\ref{i1}%
) is convergent. If $q+k+\alpha d=2\eta mn$ we need that $\kappa >1+\gamma
+\eta (1+\varepsilon _{0})$ in order to obtain the convergence of the
series. If $\kappa >1+\gamma +\eta $ then we may choose $\varepsilon _{0}$
sufficiently small in order to have $\gamma +\eta (1+\varepsilon
_{0})-\kappa >1$ and we are done. $\square $

\medskip

There are two important examples: ${\mathbf{e}}={\mathbf{e}}_p$ that we
discuss in a special subsection below and ${\mathbf{e}}={\mathbf{e}}_{\log }$
which we discuss now. We recall that ${\mathbf{e}}_{\log }\in \mathcal{E}%
_{\alpha ,\gamma }$ with $\alpha =0$ and $\gamma =1$ and $\left\Vert
f_{\delta }\right\Vert _{2m,2m,{\mathbf{e}}_{\log }}\leq C1\vee \left\Vert
f_{\delta }\right\Vert _{2m,2m,1+}$ where $\left\Vert f_{\delta }\right\Vert
_{2m,2m,1+}$\ is defined in (\ref{E22}). Then as a particular case of the
previous theorem we obtain:

\begin{theorem}
\label{ThLog} We consider a non negative finite measure $\mu $ and we
suppose that there exists a family of measures $\mu _{\delta }(dx)=f_{\delta
}(x)dx,\delta >0$ which verifies the following assumptions. There exist $%
m\in {\mathbb{N}}_{\ast }$, $C,r,\varepsilon >0$ and a function $\lambda
_{m}(\delta )$, $\delta \in (0,1)$, which is right-continuous and non
increasing such that
\begin{equation}
\left\Vert f_{\delta }\right\Vert _{2m,2m,1+}\leq \lambda _{m}(\delta )\leq
C\delta ^{-r}\qquad \mbox{and}\qquad \lambda _{m}^{\frac{1}{2m}}(\delta
)d_{1}(\mu ,\mu _{\delta })\leq \frac{C}{(\ln (1/\delta ))^{2+\frac{1}{2m}%
+\varepsilon }}.  \label{Llog}
\end{equation}%
Then $\mu (dx)=f(x)dx$ with $f\in L^{{\mathbf{e}}_{\log }}$.
\end{theorem}

\subsection{The $L^{p}$ criterion}

\label{sect:Lp}

In the case of the $L^{p}$ norms, that is ${\mathbf{e}}={\mathbf{e}}_p,$ our
result fits in the general theory of the interpolation spaces and we may
give a more precise characterization of the space $\mathcal{S}_{q,k,m,%
\mathbf{e}_p}=: \mathcal{S}_{q,k,m,p}$. We come back to the standard
notation and we denote $\left\Vert \cdot \right\Vert _{p}$ instead of $%
\left\Vert \cdot \right\Vert _{{\mathbf{e}}_p}, $ $W^{q,p}$ instead of $W^{q,%
{\mathbf{e}}_p}$ and so on. In Appendix \ref{sect-interp} we prove that in
this case the space $\mathcal{S}_{q,k,m,p} $ is related to the following
interpolation space. Let $X=W_{\ast }^{k,\infty } $ where $W_{\ast
}^{k,\infty }$ is the dual of $W^{k,\infty }$ (notice that one may look to $%
\mu \in \mathcal{M}$ as to an element of $W_\ast^{k,\infty }$ and then $%
d_{k}(\mu ,\nu )=\left\Vert \mu -\nu \right\Vert _{W_{\ast }^{k,\infty }}).$
We also take $Y=W^{q+2m,2m,p}$ and for $\gamma \in (0,1)$ we denote by $%
(X,Y)_{\gamma }$ the real interpolation space of order $\gamma $ between $X$
and $Y$ (see the Appendix \ref{sect-interp} for notations). Then we have%
\begin{equation*}
\mathcal{S}_{q,k,m,p}=(X,Y)_{\gamma }\qquad \mbox{with}\qquad \gamma =\frac{%
q+k+d/p_{\ast }}{2m}.
\end{equation*}
So Theorem \ref{ThLlog} reads
\begin{equation*}
W^{q+1,2m,p}\subset (W_{\ast }^{k,\infty },W^{q+2m,2m,p})_{\gamma }\subset
W^{q,p}.
\end{equation*}

We go now further and we notice that if (\ref{i2}) holds then the
convergence of the series in (\ref{i1}) is very fast. This allows us to
obtain some more regularity.

\begin{theorem}
\label{ThLp} Let $q,k\in {\mathbb{N}},m\in {\mathbb{N}}_{\ast }$, $p>1$ and
set
\begin{equation}  \label{etaBalance}
\eta>\frac{q+k+d/p_*}{2m}.
\end{equation}
We consider a non negative finite measure $\mu $ and a family of finite non
negative measures $\mu _{\delta }(dx)=f_{\delta }(x)dx,\delta >0$.

\medskip

\textbf{A.} We assume that there exist $C,r>0$ and a right-continuous and
non increasing function $\lambda _{q,m}(\delta )$, $\delta \in (0,1)$, such
that
\begin{equation*}
\Vert f_{\delta }\Vert _{2m+q,2m,p}\leq \lambda _{q,m}(\delta )\leq C\delta
^{-r}
\end{equation*}%
and moreover, with $\eta $ given in \eqref{etaBalance},
\begin{equation}
\lambda _{q,m}(\delta )^{\eta }d_{k}(\mu ,\mu _{\delta })\leq C.
\label{Balance'}
\end{equation}%
Then $\mu (dx)=f(x)dx$ with $f\in W^{q,p}$.

\medskip

\textbf{B.} We assume that (\ref{Balance'}) holds with $q+1$ instead of $%
q,\, $that is
\begin{equation*}
\lambda _{q+1,m}(\delta )^{\eta }d_{k}(\mu ,\mu _{\delta })\leq C.
\end{equation*}%
We denote%
\begin{equation}
s_{\eta }(q,k,m,p)=\frac{2m\eta -(q+k+d/p_{\ast })}{2m\eta }\wedge \frac{%
\eta }{1+\eta }.  \label{i11}
\end{equation}%
Then for every multi index $\alpha $ with $\left\vert \alpha \right\vert =q$
and every $s<s_{\eta }(q,k,m,p)$ we have $\partial _{\alpha }f\in \mathcal{B}%
^{s,p}$ where $\mathcal{B}^{s,p}$ is the Besov space of index $s.$
\end{theorem}

\textbf{Proof}. \textbf{A.} The fact that (\ref{Balance'}) implies $\mu
(dx)=f(x)dx$ with $f\in W^{q,p}$ is an immediate consequence of Theorem \ref%
{Criterion}.

\textbf{B.} We prove the regularity property: $g:=\partial _{\alpha }f\in
\mathcal{B}^{s,p}$ for $\left\vert \alpha \right\vert =q$ and $s<s_{\eta
}(q,k,m).$ In order to do it we will use Lemma \ref{Besov} so we have to
check (\ref{Int3}).

\smallskip

\textbf{Step 1.} We begin with the point $i)$ in (\ref{Int3}) so we have to
estimate $\left\Vert g\ast \partial _{i}\phi _{\varepsilon }\right\Vert
_{\infty }$. The reasoning is analogous with the one in the proof of Theorem %
\ref{Criterion} but we will use the first inequality in (\ref{Norm}) with $q$
replaced by $q+1$ and $k$ replaced by $k-1.$ So we define $\delta _{n}=\inf
\{\delta >0:\lambda _{q+1,m}(\delta )\leq n^{-2}2^{2mn}\}$ and we have $%
\delta _{n}\leq 2^{-\theta n}$ for $\theta <2m/r.$ We obtain
\begin{eqnarray*}
\left\Vert g\ast \partial _{i}\phi _{\varepsilon }\right\Vert _{p}
&=&\left\Vert \partial _{i}\partial _{\alpha }(f\ast \phi _{\varepsilon
})\right\Vert _{p}\leq \left\Vert f\ast \phi _{\varepsilon }\right\Vert
_{q+1,p}\leq \rho _{q+1,k-1,m,p}(\mu \ast \phi _{\varepsilon }) \\
&\leq &\sum_{n=1}^{\infty }2^{n(q+k+d/p_{\ast })}d_{k-1}(\mu \ast \phi
_{\varepsilon },\mu _{\delta _{n}}\ast \phi _{\varepsilon
})+\sum_{n=1}^{\infty }2^{-2mn}\left\Vert f_{\delta _{n}}\ast \phi
_{\varepsilon }\right\Vert _{2m+q+1,2m,p}.
\end{eqnarray*}%
By the choice of $\delta _{n}$
\begin{equation*}
\left\Vert f_{\delta _{n}}\ast \phi _{\varepsilon }\right\Vert
_{2m+q+1,2m,p}\leq \left\Vert f_{\delta _{n}}\right\Vert _{2m+q+1,2m,p}\leq
\lambda _{q+1,m}(\delta _{n})\leq \frac{1}{n^{2}}2^{2nm}
\end{equation*}%
so the second series is convergent. We estimate now the first sum. Since $%
\left\Vert f\ast \phi _{\varepsilon }\right\Vert _{k,\infty }\leq
\varepsilon ^{-1}\left\Vert f\right\Vert _{k-1,\infty }$ it follows that $%
d_{k-1}(\mu \ast \phi _{\varepsilon },\mu _{\delta _{n}}\ast \phi
_{\varepsilon })\leq \varepsilon ^{-1}d_{k}(\mu ,\mu _{\delta _{n}}).$ Then,
using (\ref{Balance'}) (with $q=1$ instead of $q)$ and the choice of $\delta
_{n}$ we obtain
\begin{eqnarray*}
2^{n(q+k+d/p_{\ast })}d_{k-1}(\mu \ast \phi _{\varepsilon },\mu _{\delta
_{n}}\ast \phi _{\varepsilon }) &\leq &\frac{C}{\varepsilon }%
2^{n(q+1+d/p_{\ast })}d_{k}(\mu ,\mu _{\delta _{n}})\leq \frac{C}{%
\varepsilon }2^{n(q+1+d/p_{\ast })}\lambda _{q+1,m}^{-\eta }(\delta _{n}) \\
&\leq &\frac{Cn^{2\eta }}{\varepsilon }2^{n(q+1+d/p_{\ast }-2m\eta )}.
\end{eqnarray*}%
We fix now $\varepsilon >0$ and we take some $n_{\varepsilon }\in {\mathbb{N}%
}$ (to be chosen in the sequel) and we write%
\begin{equation*}
\sum_{n=1}^{\infty }2^{n(q+k+d/p_{\ast })}d_{k-1}(\mu \ast \phi
_{\varepsilon },\mu _{\delta _{n}}\ast \phi _{\varepsilon })\leq
C\sum_{n=1}^{n_{\varepsilon }}2^{n(q+k+d/p_{\ast })}+\frac{C}{\varepsilon }%
\sum_{n=n_{\varepsilon }+1}^{\infty }n^{2\eta }2^{n(q+k+d/p_{\ast }-2\eta
m)}.
\end{equation*}%
We take $a>0$ and we upper bound the above series by%
\begin{equation*}
2^{n_{\varepsilon }(q+k+d/p_{\ast })}+\frac{C}{\varepsilon }%
2^{n_{\varepsilon }(q+k+d/p_{\ast }+a-2\eta m)}.
\end{equation*}%
In order to optimize we take $n_{\varepsilon }$ such that $%
2^{2mn_{\varepsilon }}=\frac{1}{\varepsilon }.$ With this choice we obtain%
\begin{equation*}
2^{n_{\varepsilon }(q+k+d/p_{\ast }+a)}\leq C\varepsilon ^{-\frac{%
q+k+d/p_{\ast }+a}{2m\eta }}.
\end{equation*}%
We conclude that
\begin{equation*}
\left\Vert g\ast \partial _{i}\phi _{\varepsilon }\right\Vert _{p}\leq
C\varepsilon ^{-\frac{q+k+d/p_{\ast }+a}{2m\eta }}
\end{equation*}%
which means (\ref{Int3}) $i)$ holds for $s<1-\frac{q+k+d/p_{\ast }}{2m\eta }%
. $

\textbf{Step 2}. We check now (\ref{Int3}) $ii)$ so we have to estimate $%
\left\Vert g\ast \phi _{\varepsilon }^{i}\right\Vert _{p}$ with $\phi
_{\varepsilon }^{i}(x)=x^{i}\phi _{\varepsilon }(x).$ We take $u\in (0,1)$
(to be chosen in a moment) and we define
\begin{equation*}
\delta _{n,\varepsilon }=\inf \{\delta >0:\lambda _{q+1,m}(\delta )\leq
n^{-2}2^{2mn}\times \varepsilon ^{-(1-u)}\}.
\end{equation*}%
Then we proceed as in the previous step:
\begin{eqnarray*}
\left\Vert \partial _{i}(g\ast \phi _{\varepsilon }^{i})\right\Vert _{p}
&\leq &\rho _{q+1,k-1,m,p}(\mu \ast \phi _{\varepsilon }^{i}) \\
&\leq &\sum_{n=1}^{\infty }2^{n(q+k+d/p_{\ast })}d_{k-1}(\mu \ast \phi
_{\varepsilon }^{i},\mu _{\delta _{n,\varepsilon }}\ast \phi _{\varepsilon
}^{i})+\sum_{n=1}^{\infty }2^{-2mn}\left\Vert f_{\delta _{n,\varepsilon
}}\ast \phi _{\varepsilon }^{i}\right\Vert _{2m+q+1,2m,p}.
\end{eqnarray*}%
It is easy to check that for every $h\in L^{p}$ one has $\left\Vert h\ast
\phi _{\varepsilon }^{i}\right\Vert _{p}\leq \varepsilon \left\Vert
h\right\Vert _{p}$ so that, by our choice of $\delta _{n,\varepsilon }$ we
obtain%
\begin{equation*}
\left\Vert f_{\delta _{n,\varepsilon }}\ast \phi _{\varepsilon
}^{i}\right\Vert _{2m+q+1,2m,p}\leq \varepsilon \left\Vert f_{\delta
_{n,\varepsilon }}\right\Vert _{2m+q+1,2m,p}\leq \varepsilon \times \frac{%
2^{2mn}}{n^{2}}\times \varepsilon ^{-(1-u)}.
\end{equation*}%
It follows that the second sum is upper bounded by $C\varepsilon ^{u}.$

Since $\left\Vert \partial _{j}h\ast \phi _{\varepsilon }^{i}\right\Vert
_{\infty }\leq C\left\Vert h\right\Vert _{\infty }$ it follows that
\begin{equation*}
d_{k-1}(\mu \ast \phi _{\varepsilon }^{i},\mu _{\delta _{n,\varepsilon
}}\ast \phi _{\varepsilon }^{i})\leq Cd_{k}(\mu ,\mu _{\delta
_{n,\varepsilon }})\leq \frac{C}{\lambda _{q+1,m}^{\eta }(\delta
_{n,\varepsilon })}=\frac{Cn^{2}}{2^{2mn\eta }}\varepsilon ^{\eta (1-u)}.
\end{equation*}%
Since $2m\eta >q+k+d/p_{\ast }$ the first sum is convergent also and is
upper bounded by $C\varepsilon ^{\eta (1-u)}.$ We conclude that
\begin{equation*}
\left\Vert \partial _{i}(g\ast \phi _{\varepsilon }^{i})\right\Vert _{p}\leq
C\varepsilon ^{\eta (1-u)}+C\varepsilon ^{u}.
\end{equation*}%
In order to optimize we take $u=\frac{\eta }{1+\eta }.\square $

\subsection{Convergence criteria in $W^{q,p}$ and $W^{q,{\mathbf{e}}_{\log}}$%
}

\label{sect:convLp}

For a function $f$, we denote $\mu _{f}(dx)=f(x)dx.$

\begin{theorem}
\label{Conv} Let $\eta :{\mathbb{R}}_{+}\rightarrow {\mathbb{R}}_{+}$ be a
non decreasing function and $a\geq 1$ be such that%
\begin{equation}
\lim_{n\rightarrow \infty }\eta (n)=+\infty \quad \mbox{and}\quad \eta
(n+1)\leq a\eta (n),\quad \mbox{for every }n\in {\mathbb{N}}.  \label{cbis1}
\end{equation}%
Let $m,k,q\in {\mathbb{N}}$ be fixed. Let $f_{n}$, $n\in {\mathbb{N}}$, be a
sequence of functions and $\mu \in \mathcal{M}$.

\medskip

$i)$ Let $p\geq 1$. If there exists $\alpha >\frac{q+k+d/p_{\ast }}{m}$ such
that
\begin{equation}
\left\Vert f_{n}\right\Vert _{q+2m,2m,p}\leq \eta ^{1/\alpha }(n)\quad %
\mbox{and}\quad d_{k}(\mu ,\mu _{f_{n}})\leq \frac{1}{\eta (n)},
\label{cbis1p}
\end{equation}%
then $\mu (dx)=f(x)dx$ for some $f\in W^{q,p}$. Moreover, there exists a
constant $C$ depending on $a,\alpha $ such that for every $n\in {\mathbb{N}}$%
\begin{equation}
\left\Vert f-f_{n}\right\Vert _{q,p}\leq C\eta^{-\theta }(n)\quad \mbox{with}%
\quad \theta =\frac{1}{\alpha }\wedge (1-\frac{q+k+d/p_{\ast }}{\alpha m}).
\label{cbis3p}
\end{equation}

$ii)$ If there exists $\alpha >\frac{q+k}{m}$ such that
\begin{equation}
\left\Vert f_{n}\right\Vert _{q+2m,2m,1+}\leq \eta ^{1/\alpha }(n)\quad %
\mbox{and}\quad d_{k}(\mu ,\mu _{f_{n}})\leq \frac{1}{\eta (n)},
\label{cbis1elog}
\end{equation}%
then $\mu (dx)=f(x)dx$ for some $f\in W^{q,{\mathbf{e}}_{\log }}$. Moreover,
there exists a constant $C$ depending on $a,\alpha $ such that for every $%
n\in {\mathbb{N}}$%
\begin{equation}
\left\Vert f-f_{n}\right\Vert _{q,{\mathbf{e}}_{\log }}\leq C(\eta
^{-1/\alpha }(n)+(\log _{2}\eta (n))\eta ^{-(1-\frac{q+k}{\alpha m}%
)}(n))=:\varepsilon _{n}(\alpha ).  \label{cbis3elog}
\end{equation}%
And if $\varepsilon _{n}(\alpha )\leq 1$ then
\begin{equation}
\sum_{0\leq \left\vert \alpha \right\vert \leq q}\int \left\vert (\partial
_{\alpha }f-\partial _{\alpha }f_{n})(x)\right\vert (1+\ln ^{+}\left\vert
(\partial _{\alpha }f-\partial _{\alpha }f_{n})(x)\right\vert dx\leq
2C_{\ast }\varepsilon _{n}(\alpha ).  \label{cbis4elog}
\end{equation}
\end{theorem}

\textbf{Proof.} $i)$ \textbf{Step 1}. For $r\in {\mathbb{N}}$, we define%
\begin{equation*}
n_{r}=\min \{n:\eta (n)\geq 2^{\alpha rm}\}\quad \mbox{and}\quad r_{n}=\min
\{r\in N:n_{r}\geq n\}.
\end{equation*}%
Then we have
\begin{equation}
\frac{1}{a}\eta (n)\leq 2^{\alpha r_{n}m}\leq C\eta (n).  \label{cbis6}
\end{equation}

\bigskip Since $\{r\in N:n_{r}\geq n\}$ is a discrete set, its minimum $%
r_{n} $ belongs to this set, so $n_{r_{n}}\geq n.$ Then $\eta (n)\leq \eta
(n_{r_{n}})\leq a\eta (n_{r_{n}}-1)\leq a2^{\alpha r_{n}m}.$ On the other
hand, since $r_{n}-1\notin \{r\in N:n_{r}\geq n\}$ one has $n>n_{r_{n}-1}$
and then $\eta (n)\geq \eta (n_{r_{n}-1})\geq 2^{\alpha
(r_{n}-1)m}=C^{-1}2^{\alpha r_{n}m}$ with $C=2^{\alpha m}.$ So, (\ref{cbis6}%
) holds.

\smallskip

\textbf{Step 2}. We fix $n\in N$ and for $r\in N$ we define%
\begin{equation*}
g_{r}=0\mbox{ if }r<r_{n}\mbox{ and }g_{r}=f_{n_{r}}-f_{n}\mbox{ if }r\geq
r_{n}
\end{equation*}%
and $\nu (dx)=\mu (dx)-f_{n}(x)dx,\nu _{r}(dx)=g_{r}(x)dx.$ Using (\ref%
{Oo10bis}) (recall that $\beta _{\mathbf{e}_{p}}=t^{1/p_{\ast }}$) we get%
\begin{equation*}
\rho _{q,k,m,p}(\nu )\leq \sum_{r=1}^{\infty }2^{r(q+k+d/p_{\ast
})}d_{k}(\nu ,\nu _{r})+\sum_{r=1}^{\infty }2^{-2mr}\left\Vert
g_{r}\right\Vert _{q+2m,2m,p}=:S_{1}+S_{2}.
\end{equation*}

We estimate $S_{1}.$ For $r<r_n$ we have $\nu _{r}=0$ so that $d_{k}(\nu
,\nu _{r})=d_{k}(\nu ,0)=d_{k}(\mu,\mu _{f_{n}})\leq \eta^{-1}(n).$ And for $%
r\geq r_n$ we have

\begin{equation*}
d_{k}(\nu ,\nu _{r})=d_{k}(\mu,\mu _{f_{n_{r}}})\leq \frac{1}{\eta (n_{r})}%
\leq \frac{1}{2^{rm\alpha }}.
\end{equation*}%
So, we obtain%
\begin{equation*}
S_{1}\leq 2^{r_n(q+k+d/p_{\ast })}\eta^{-1}(n)+\frac{C}{2^{r_nm\alpha(1 -%
\frac{(q+k+d/p_{\ast }}{\alpha m})}}
\end{equation*}%
and using (\ref{cbis6}),
\begin{equation*}
S_{1}\leq C\eta^{-(1-\frac{q+k+d/p_{\ast }}{\alpha m})}(n).
\end{equation*}

We estimate now $S_{2}.$ We have $g_{r}=0$ for $r<r_{n}$ and for $r\geq
r_{n} $
\begin{equation*}
\left\Vert g_{r}\right\Vert _{q+2m,2m,p}\leq \left\Vert f_{n_{r}}\right\Vert
_{q+2m,2m,p}+\left\Vert f_{n}\right\Vert _{q+2m,2m,p}\leq \eta
(n_{r})^{1/\alpha }+\eta (n)^{1/\alpha }.
\end{equation*}%
But $\eta (n_{r})\leq a\eta (n_{r}-1)\leq a2^{\alpha rm}$, so that
\begin{equation*}
\left\Vert g_{r}\right\Vert _{q+2m,2m,p}\leq a^{1/\alpha }2^{rm}+\eta
(n)^{1/\alpha }.
\end{equation*}%
It follows that
\begin{equation*}
S_{2}\leq a^{1/\alpha }\sum_{r\geq r_{n}}2^{-rm}+\eta (n)^{1/\alpha
}\sum_{r\geq r_{n}}2^{-2rm}\leq C\big(2^{-r_{n}m}+\eta (n)^{1/\alpha
}2^{-2r_{n}m}\big)
\end{equation*}%
and using (\ref{cbis6}) we get
\begin{equation*}
S_{2}\leq C\eta (n)^{-1/\alpha }.
\end{equation*}%
Then, we obtain
\begin{equation*}
\rho _{q,k,m,p}(\nu )\leq C(\eta ^{-1/\alpha }(n)+\eta ^{-(1-\frac{%
q+k+d/p_{\ast }}{\alpha m})}(n))
\end{equation*}%
and Theorem \ref{2C} allows one to conclude.

\smallskip

$ii)$ We take $n_{r}$ and $r_{n}$ as in Step 1 above, giving (\ref{cbis6}),
and we take $g_{r}$, $\nu $, $\nu _{r}$ as in Step 2 above. Then, by using (%
\ref{Oo10bis}) we get%
\begin{equation*}
\rho _{q,k,m,{\mathbf{e}}_{\log }}(\nu )\leq \sum_{r=1}^{\infty
}2^{r(q+k)}\beta _{{\mathbf{e}}_{\log }}(2^{rd})d_{k}(\nu ,\nu
_{r})+\sum_{r=1}^{\infty }2^{-2mr}\left\Vert g_{r}\right\Vert _{q+2m,2m,{%
\mathbf{e}}_{\log }}.
\end{equation*}%
By (\ref{E23}) and (\ref{E1}), we can write
\begin{equation*}
\rho _{q,k,m,{\mathbf{e}}_{\log }}(\nu )\leq C\sum_{r=1}^{\infty
}2^{r(q+k)}rd_{k}(\nu ,\nu _{r})+\sum_{r=1}^{\infty }2^{-2mr}1\vee
\left\Vert g_{r}\right\Vert _{q+2m,2m,1+}=:S_{1}+S_{2}.
\end{equation*}%
Concerning $S_{1}$, for $r<r_{n}$ we have $d_{k}(\nu ,\nu _{r})=d_{k}(\nu
,0)=d_{k}(\mu ,\mu _{f_{n}})\leq \eta ^{-1}(n)$ and for $r\geq r_{n}$ we
have $d_{k}(\nu ,\nu _{r})\leq \frac{1}{\eta (n_{r})}\leq \frac{1}{%
2^{rm\alpha }}.$ So, we obtain%
\begin{equation*}
S_{1}\leq C\Big(r_{n}2^{r_{n}(q+k)}\eta ^{-1}(n)+\frac{r_{n}}{%
2^{r_{n}m\alpha (1-\frac{q+k}{\alpha m})}}\Big).
\end{equation*}%
Using (\ref{cbis6}),
\begin{equation*}
S_{1}\leq Cr_{n}\eta ^{-(1-\frac{q+k}{\alpha m})}(n)\leq C(\log _{2}\eta
(n))\,\eta ^{-(1-\frac{q+k}{\alpha m})}(n).
\end{equation*}%
As for $S_{2}$, we proceed as in Step 2 above and we obtain $S_{2}\leq C\eta
(n)^{-1/\alpha }.$ Then,
\begin{equation*}
\rho _{q,k,m,{\mathbf{e}}_{\log }}(\nu )\leq C(\eta ^{-1/\alpha }(n)+\eta
^{-(1-\frac{q+k+d/p_{\ast }}{\alpha m})}(n))
\end{equation*}%
and the statement again follows from Theorem \ref{2C}. So (\ref{cbis3elog})
is proved. In order to check (\ref{cbis4elog}) we use (\ref{E24bis}) (notice
that, since $\left\Vert f-f_{n}\right\Vert _{q,{\mathbf{e}}_{\log }}\leq
\varepsilon _{n}(\alpha )\leq 1,$ we have $\ln ^{+}\left\Vert
f-f_{n}\right\Vert _{q,{\mathbf{e}}_{\log }}=0).$ $\square $

\subsection{Random variables and integration by parts}

\label{sect-RV}

In this section we work in the framework of random variables. For a random
variable $F$ we denote by $\mu _{F}$ the law of $F$ and if $\mu _{F}$ is
absolutely continuous we denote by $p_{F}$ its density. We will use Theorem %
\ref{ThLp} for $\mu _{F}$ so we will look for a family of random variables $%
F_{\delta },\delta >0$ such that $\mu _{F_{\delta }}$ satisfy the hypothesis
of this theorem. Sometimes it is easy to construct such a family with
explicit densities $p_{F_{\delta }}$ and then one may check (\ref{Balance'})
directly (this is the case in the examples in Section \ref{sect-pathdep} and %
\ref{sect-heat}). But sometimes one does not know $p_{F_{\delta }}$ and then
it is useful to use the integration by parts machinery in order to prove (%
\ref{Balance'}) - this is the case in the example given is Section \ref%
{jumps} or the application to a kind of generalization of the H\"ormander
condition to general Wiener functionals developed in \cite{bib:BC-Horm}.

We briefly recall the abstract definition of integration by parts formulae
and we give some useful properties (coming essentially from \cite{bib:[BCa]}%
). We consider two random variables $F=(F_{1},...,F_{d})$ and $G.$ Given a
multi index $\alpha =(\alpha _{1},...,\alpha _{k})\in \{1,...,d\}^{k}$ and
for $p\geq 1$ we say that $\mathrm{IP}_{\alpha ,p}(F,G)$ holds if we may
find a random variable $H_{\alpha }(F;G)\in L^{p}$ such that for every $f\in
C^{\infty }({\mathbb{R}}^{d})$ one has%
\begin{equation}
{\mathbb{E}}(\partial _{\alpha }f(F)G)={\mathbb{E}}(f(F)H_{\alpha }(F;G)).
\label{Lp5'}
\end{equation}%
The weight $H_{\alpha }(F;G)$ is not uniquely determined: the one with the
lowest variance is ${\mathbb{E}}(H_{\alpha }(F;G)\mid \sigma (F))$. This
quantity is uniquely determined. So we denote%
\begin{equation}
\theta _{\alpha }(F,G)={\mathbb{E}}(H_{\alpha }(F;G)\mid \sigma (F)).
\label{ip5'}
\end{equation}

For $m\in {\mathbb{N}}$ and $p\geq 1$ we denote by $\mathcal{R}_{m,p}$ the
class of random variables $F$ in ${\mathbb{R}}^{d}$ such that $\mathrm{IP}%
_{\alpha ,p}(F,1)$ holds for every \ multi index $\alpha $ with $\left\vert
\alpha \right\vert \leq m.$ We define
\begin{equation}
T_{m,p}(F)=\left\Vert F\right\Vert _{p}+\sum_{\left\vert \alpha \right\vert
\leq m}\left\Vert \theta _{\alpha }(F,1)\right\Vert _{p}.  \label{ip6}
\end{equation}%
Notice that by H\"{o}lder's inequality $\left\Vert {\mathbb{E}}(H_{\alpha
}(F;1)\mid \sigma (F))\right\Vert _{p}\leq \left\Vert H_{\alpha
}(F;1)\right\Vert _{p}.$ It follows that for every choice of the weights $%
H_{\alpha }(F;1)$ one has
\begin{equation}
T_{m,p}(F)\leq \left\Vert F\right\Vert _{p}+\sum_{\left\vert \alpha
\right\vert \leq m}\left\Vert H_{\alpha }(F;1)\right\Vert _{p}.  \label{ip7}
\end{equation}

\begin{theorem}
\label{Riez}Let $m,l\in {\mathbb{N}}$ and $p>d.$ If $F\in \mathcal{R}%
_{m+1,p} $ then the law of $F$ is absolutely continuous and the density $%
p_{F}$ belongs to $C^{m}({\mathbb{R}}^{d}).$ Moreover, suppose that $F\in
\mathcal{R}_{m+1,2(d+1)}.$ There exists a universal constant $C$ (depending
on $d,l$ and $m$ only) such that for every multi index $\alpha $ with $%
\left\vert \alpha \right\vert \leq m$%
\begin{equation}
\left\vert \partial _{\alpha }p_{F}(x)\right\vert \leq C
T^{d^{2}-1}_{1,2(d+1)}(F) T_{m+1,2(d+1)}(F)(1+\left\Vert F\right\Vert
_{l})(1+\left\vert x\right\vert )^{-l}.  \label{ip8}
\end{equation}%
In particular, for every $q\geq 1,k\in {\mathbb{N}}$ there exists a
universal constant $C$ (depending on $d,m,k,p$ and $q$) such that
\begin{equation}
\left\Vert p_{F}\right\Vert _{m,k,q}\leq
CT_{1,2(d+1)}^{d^{2}-1}(F)T_{m+1,2(d+1)}(F)(1+\left\Vert F\right\Vert
_{d+k+1}).  \label{ip9}
\end{equation}
\end{theorem}

\textbf{Proof}. The proof is an immediate consequence of the results in \cite%
{bib:[BCa]}. In order to see this we have to give the relation between the
notation used in that paper and the notation used here: we work with the
probability measure $\mu _{F}(dx)={\mathbb{P}}(F\in dx)$ and in \cite%
{bib:[BCa]} we use the notation $\partial _{\alpha }^{\mu _{F}}g(x)={\mathbb{%
E}}(H_{\alpha }(F;g(F))\mid F=x).$

The fact that $F\in \mathcal{R}_{m+1,p}$ implies that $F\sim p_{F}(x)dx$
with $p_{F}\in C^{m}({\mathbb{R}}^{d})$ is proved in \cite{bib:[BCa]}
Proposition 9. We consider now a function $\psi \in C_{b}^{\infty }(R^{d})$
such that $1_{B_{1}}\leq \psi \leq 1_{B_{2}}.$ In \cite{bib:[BCa]} Theorem 8
we have given the following representation formula:%
\begin{equation*}
\partial _{\alpha }p_{F}(x)=\sum_{i=1}^{d}{\mathbb{E}}(\partial
_{i}Q_{d}(F-x)\theta _{(\alpha ,i)}(F;\psi (F-x))1_{B_{2}}(F-x))
\end{equation*}%
where $B_{r}$ denotes the ball centered at $0$ with radius $r$, $Q_{d}$ is
the Poisson kernel on ${\mathbb{R}}^{d}$ and, if $\alpha =(\alpha
_{1},...,\alpha _{k})$, then $(\alpha ,i)=(\alpha _{1},...,\alpha _{k},i).$
Using H\"{o}lder's inequality we obtain (with $p_{\ast }$ the conjugate of $%
p)$
\begin{equation*}
\left\vert \partial _{\alpha }p_{F}(x)\right\vert \leq
\sum_{i=1}^{d}\left\Vert \partial _{i}Q_{d}(F-x)\right\Vert _{p}\left\Vert
\theta _{(\alpha ,i)}(F;\psi (F-x))1_{B_{2}}(F-x)\right\Vert _{p_{\ast }}.
\end{equation*}%
We take $p=d+1$ so that $p_{\ast }=(d+1)/d\leq 2.$ In \cite{bib:[BCa]}
Theorem 5 we proved that
\begin{equation*}
\left\Vert \partial _{i}Q_{d}(F-x)\right\Vert _{p}\leq C
T_{1,2(d+1)}^{d^{2}-1}(F).
\end{equation*}%
Moreover we have the following computational rule (Lemma 9 in \cite%
{bib:[BCa]})
\begin{equation*}
\theta _{i}(F,fg(F))=f(F)\theta _{i}(F,g(F))+(g\partial _{i}f)(F).
\end{equation*}%
Since $\psi \in C_{b}^{\infty }({\mathbb{R}}^{d})$ we may use the above
formula in order to get%
\begin{eqnarray*}
\left\Vert \theta _{(\alpha ,i)}(F;\psi (F-x))1_{B_{2}}(F-x)\right\Vert
_{p_{\ast }} &\leq &\left\Vert \theta _{(\alpha ,i)}(F;\psi
(F-x))\right\Vert _{2p_{\ast }}\sqrt{{\mathbb{P}}(\left\vert F-x\right\vert
\leq 2)} \\
&\leq &C_{\psi }T_{\left\vert \alpha \right\vert +1,2p_{\ast }}(F)\,\sqrt{{%
\mathbb{P}}(\left\vert F-x\right\vert \leq 2)}.
\end{eqnarray*}%
For $\left\vert x\right\vert \geq 4$
\begin{equation*}
{\mathbb{P}}(\left\vert F-x\right\vert \leq 2)\leq {\mathbb{P}}(\left\vert
F\right\vert \geq \frac{1}{2}\left\vert x\right\vert )\leq \frac{2^{k}}{%
\left\vert x\right\vert ^{k}}{\mathbb{E}}(\left\vert F\right\vert ^{k})
\end{equation*}%
so the proof of \eqref{ip8} is completed. $\square $

\medskip

We are now ready to rewrite Theorem \ref{ThLp}:

\begin{theorem}
\label{IP} Let $k,q\in {\mathbb{N}}$, $m\in{\mathbb{N}}_*$, $p>1$ and let
\begin{equation*}
\eta>\frac{q+k+d/p_*}{2m},
\end{equation*}
$p_*$ denoting the conjugate of $p$. Let $F$, $F_\delta$,$\delta>0$, be
random variables and let $\mu_F$, $\mu_{F_\delta}$, $\delta>0$, denote the
associated laws.

\medskip

\textbf{A.} Suppose that $F_{\delta }\in \mathcal{R}_{2m+q+1,2(d+1)}$, $%
\delta >0$ are uniformly bounded in $L^{2m+d+1}$ and that there exist $C>0$
and $\theta >0$ such that
\begin{align}
& T_{2m+q+1,2(d+1)}(F_{\delta })\leq C\delta ^{-\theta (2m+q+1)},
\label{i13} \\
& d_{k}(\mu _{F},\mu _{F_{\delta }})\leq C\delta ^{\theta \eta
d^{2}(2m+q+1)}.  \label{i14}
\end{align}%
Then $\mu _{F}(dx)=p_{F}(x)dx$ with $p_{F}\in W^{q,p}.$

\medskip

\textbf{B.} Suppose that $F_{\delta }\in \mathcal{R}_{2m+q+2,2(d+1)}$, $%
\delta >0$, and (\ref{i13}) holds with $q+1$ instead of $q.$Then for every
multi index $\alpha $ with $\left\vert \alpha \right\vert =q$ and every $%
s<s_{\eta }(q,k,m,p)$ we have $\partial _{\alpha }p_{F}\in \mathcal{B}^{s,p}$
where $\mathcal{B}^{s,p}$ is the Besov space of index $s$ and $s_{\eta
}(q,k,m,p)$ is given in (\ref{i11}).
\end{theorem}

\textbf{Proof}. \textbf{A.} Let $n,l\in {\mathbb{N}}$ and $p>1$ be fixed. By
using (\ref{i13}) and (\ref{ip9}) we obtain $\left\Vert p_{F_{\delta
}}\right\Vert _{2m+q,2m,p}\leq C\delta ^{-\theta d^{2}(2m+q+1)}.$ So, as a
consequence of (\ref{i14}) we obtain $\left\Vert p_{F_{\delta }}\right\Vert
_{2m+q,2m,p}^{\eta }d_{k}(\mu _{F},\mu _{F_{\delta }})\leq C.$ And we apply
Theorem \ref{ThLp} \textbf{A}. Similarly, \textbf{B} follows by applying
Theorem \ref{ThLp} \textbf{B.} $\square $

\section{Examples}

\label{sect-examples}

\subsection{Path dependent SDE's}

\label{sect-pathdep}

In this section we look to the SDE%
\begin{equation}
dX_{t}=\sum_{j=1}^{n}\sigma _{j}(t,X)dW_{t}^{j}+b(t,X)dt  \label{Diff1}
\end{equation}%
where $W=(W^{1},...,W^{n})$ is a standard Brownian motion and $\sigma
_{j},b:C({\mathbb{R}}_{+};{\mathbb{R}}^{d})\rightarrow C({\mathbb{R}}_{+};{%
\mathbb{R}}^{d})$, $j=1,...,n$. We use the notation $\sigma_j(t,\varphi)=%
\sigma_j(\varphi)(t)$ and $b(t,\varphi)=b(\varphi)(t)$, $\varphi\in C({%
\mathbb{R}}_+;{\mathbb{R}}^d)$. If $\sigma _{j}$ and $b$ satisfy some
Lipschitz continuity property with respect to the sup-norm on $C({\mathbb{R}}%
_{+};{\mathbb{R}}^{d})$ then this equation has a unique solution. But we do
not want to make such an hypothesis here so we just consider an adapted
process $X_{t},t\geq 0$ which verifies the above equation.

We set $\Delta _{s,t}(w):=\sup_{s\leq u\leq t}\left\vert
w_{u}-w_{s}\right\vert $

\begin{theorem}
\label{RegLog} Let $b$ and $\sigma _{j}$, $j=1,\ldots,n$, be bounded.
Suppose that there exists $\varepsilon ,C>0$ such that
\begin{equation}
\left\vert \sigma _{j}(t,w)-\sigma _{j}(s,w)\right\vert \leq C\Big(\ln \Big(%
\frac{1}{\Delta _{s,t}(w)}\Big)\Big)^{-(2+\varepsilon )},\qquad \forall
j=1,...,n  \label{Diff3}
\end{equation}%
and that there exists some $\lambda^*\geq \lambda_* >0$ such that
\begin{equation}
\lambda^*\geq \sigma \sigma ^{\ast }(t,w)\geq \lambda_* \qquad \forall t\geq
0,w\in C({\mathbb{R}}_{+};{\mathbb{R}}^{d}).  \label{Diff2}
\end{equation}
Then for every $T>0$ the law of $X_{T}$ is absolutely continuous with
respect to the Lebesgue measure and the density belongs to $L^{{\mathbf{e}}%
_{\log}}$.
\end{theorem}

\begin{remark}
We note that in the particular case of standard SDE's we have $\sigma
_{j}(t,w)=\sigma _{j}(w_{t})$ and a sufficient condition in order that (\ref%
{Diff3}) holds is $\left\vert \sigma _{j}(x)-\sigma _{j}(y)\right\vert \leq
C(\ln (\frac{1}{\left\vert x-y\right\vert }))^{-(2+\varepsilon )}.$ This is
weaker than H\"{o}lder continuity.
\end{remark}

\textbf{Proof}. For $\delta >0$ we construct
\begin{equation*}
X_{T}^{\delta }=X_{T-\delta }+\sum_{j=1}^{n}\sigma _{j}(T-\delta
,X)(W_{T}^{j}-W_{T-\delta }^{j}).
\end{equation*}%
We will use Theorem \ref{ThLog} so we check the hypotheses there.

\smallskip

\textbf{Step 1.} We write $X_{T}-X_{T}^{\delta }=\sum_{j=1}^{n}I_{\delta
}^{j}+J_{\delta }$ with
\begin{equation*}
I_{\delta }^{j}=\int_{T-\delta }^{T}(\sigma _{j}(t,X)-\sigma _{j}(T-\delta
,X))dW_{t}^{j}\quad\mbox{and}\quad J_{\delta }=\int_{T-\delta }^{T}b(t,W)dt.
\end{equation*}%
Since $b$ is bounded, we have
\begin{equation}
{\mathbb{E}}(\left\vert J_{\delta }\right\vert) \leq C\delta .  \label{Ito5}
\end{equation}%
Let $a_{\delta }=\sqrt{\delta}\, \ln \frac{1}{\delta }$ and $A_{\delta
}=\{\Delta _{T-\delta ,T}(X)\leq a_{\delta }\}.$ We write ${\mathbb{E}}%
(\vert I_{\delta }^{j}\vert ^{2})=K_{\delta }+L_{\delta }$ with
\begin{eqnarray*}
K_{\delta } &=&\int_{T-\delta }^{T}{\mathbb{E}}(1_{A_{\delta
}^{c}}\left\vert \sigma _{j}(t,X)-\sigma _{j}(T-\delta ,X)\right\vert ^{2})dt
\\
L_{\delta } &=&\int_{T-\delta }^{T}{\mathbb{E}}(1_{A_{\delta }}\left\vert
\sigma _{j}(t,X)-\sigma _{j}(T-\delta ,X)\right\vert ^{2})dt.
\end{eqnarray*}%
By using the Bernstein's inequality we obtain ${\mathbb{P}}(A_{\delta
}^{c})\leq C\exp (-\frac{a_{\delta }^{2}}{C^{\prime }\delta })$. And since $%
\sigma_j$ is bounded, for any small $\delta$ we get
\begin{equation*}
K_{\delta }\leq C\delta{\mathbb{P}}(A_{\delta }^{c})\leq C\delta \exp (-%
\frac{a_{\delta }^{2}}{2C ^{\prime }\delta })\leq C\delta ^{\frac{3}{2}}.
\end{equation*}%
Moreover using (\ref{Diff3}) and again for $\delta$ small enough,
\begin{equation*}
L_{\delta }\leq \frac{C\delta }{(\ln \frac{1}{a_{\delta }})^{2(2+\varepsilon
)}}\leq \frac{C^{\prime }\delta }{(\ln \frac{1}{\delta })^{2(2+\varepsilon )}%
}
\end{equation*}%
(notice that $\ln(\frac 1\delta)/\ln \frac 1{a_\delta}\to \frac 12>0$ for $%
\delta\to 0$). We conclude that
\begin{equation*}
{\mathbb{E}}(\vert I_{\delta }^{j}\vert ^{2})\leq \frac{C\delta }{(\ln \frac{%
1}{\delta })^{2(2+\varepsilon )}}
\end{equation*}%
so that, if $\mu $ is the law of $X_{T}$ and $\mu _{\delta }$ is the law of $%
X_{T}^{\delta }$ then for every $\delta$ small,
\begin{equation}
d_{1}(\mu ,\mu _{\delta })\leq {\mathbb{E}}(\vert X_{T}-X_{T}^{\delta
}\vert) \leq \frac{C\delta ^{1/2}}{(\ln \frac{1}{\delta })^{2+\varepsilon }}.
\label{Ito6}
\end{equation}

\textbf{Step 2}. Given a positive definite matrix $a$, we denote
\begin{equation*}
\gamma _{\delta ,a}(y)=\frac{1}{(2\pi \delta )^{d/2}(\det a)^{1/2}}\exp \Big(%
-\frac{1}{2\delta }\langle a^{-1}y,y\rangle \Big).
\end{equation*}
With $\mu_\delta$ denoting the law of $X^\delta_T$, we have $\mu _{\delta
}(dy)=p_{\delta }(y)dy$ where
\begin{equation*}
p_{\delta }(y)={\mathbb{E}}(\gamma _{\delta ,a_{T-\delta
}(X)}(y-X_{T-\delta}))\quad \mbox{with}\quad a_{t}(X)=\sigma \sigma ^{\ast
}(t,X).
\end{equation*}%
Let $\alpha$ denote a multi index $\left\vert \alpha \right\vert =q$, $k\in {%
\mathbb{N}}$ and $\delta \leq 1 $. By using (\ref{Diff2}) we have%
\begin{align}  \label{Ito7.1}
\left\vert \partial _{\alpha }p_{\delta }(y)\right\vert &\leq C\delta ^{-q/2}%
{\mathbb{E}}\Big(\Big(1+\frac{\left\vert y-X_{T-\delta}\right\vert }{\delta
^{1/2}}\Big)^{q}\gamma _{\delta ,a_{T-\delta }(X)}(y-X_{T-\delta})\Big)
\notag \\
&\leq C\delta^{-q/2}{\mathbb{E}}\Big(\Big(1+\frac{\left\vert
y-X_{T-\delta}\right\vert }{\delta ^{1/2}}\Big)^{q}\gamma _{\delta
,\lambda^*I}(y-X_{T-\delta})\Big).
\end{align}
We use the fact that $0<x\mapsto (1+x)^{q}e^{- x^{2}}$ is bounded. This
gives
\begin{equation*}
\left\vert \partial _{\alpha }p_{\delta }(y)\right\vert \leq C\delta
^{-(d+q)/2},
\end{equation*}
so that, for small values of $\delta$,
\begin{equation}  \label{Ito7}
\ln ^{+}\left\vert \partial _{\alpha }p_{\delta }(y)\right\vert \leq C\Big (%
1+\ln \frac{1}{\delta }\Big )\leq C\,\ln \frac{1}{\delta }.
\end{equation}
Let $m\in {\mathbb{N}}$. Using (\ref{Ito7.1}) and (\ref{Ito7}) we obtain%
\begin{align*}
\left\Vert \partial _{\alpha }p_{\delta }\right\Vert _{2m,1+} &=\int
(1+\left\vert y\right\vert )^{2m}\left\vert \partial _{\alpha }p_{\delta
}(y)\right\vert (1+\ln^+|y|+\ln ^{+}\left\vert \partial _{\alpha }p_{\delta
}(y)\right\vert )dy \\
&\leq C\delta ^{-q/2}\,\ln \frac{1}{\delta }{\mathbb{E}}\Big( \int
(1+|y|)^{2m+1} \Big(1+\frac{|y-X_{T-\delta}|}{\delta^{1/2}} \Big)^{q} \gamma
_{\delta ,\lambda^*I}(y-X_{T-\delta})dy\Big) \\
&= C\delta ^{-q/2}\,\ln \frac{1}{\delta }{\mathbb{E}}\Big( \int
(1+|X_{T-\delta}+\delta^{1/2}z| )^{2m+q+1} \gamma _{1,\lambda^*I}(z)dz\Big)
\\
&\leq C\delta ^{-q/2}\ln \frac{1}{\delta }.
\end{align*}
We conclude that
\begin{equation}
\left\Vert p_{\delta }\right\Vert _{2m,2m,1+}=\sum_{0\leq \left\vert \alpha
\right\vert \leq 2m}\left\Vert \partial _{\alpha }p_{\delta }\right\Vert
_{2m,1+}\leq C\delta ^{-m}\ln \frac{1}{\delta }.  \label{Ito8}
\end{equation}

\textbf{Step 3. }We are now ready to check (\ref{Llog}): the exists $\delta
_{0}\leq 1$ such that for $\delta <\delta _{0}$ one has
\begin{eqnarray*}
\left\Vert p_{\delta }\right\Vert _{2m,2m,1+}^{1/2m}d_{1}(\mu ,\mu _{\delta
}) &\leq &C\delta ^{-1/2}\Big(\ln \frac{1}{\delta }\Big)^{1/2m}\times \frac{%
\delta ^{1/2}}{(\ln \frac{1}{\delta })^{2+\varepsilon }} \\
&=&\frac{C}{(\ln \frac{1}{\delta })^{2+\varepsilon -\frac{1}{2m}}}\leq \frac{%
C}{(\ln \frac{1}{\delta })^{2+\frac{1}{2m}+\varepsilon /2}}
\end{eqnarray*}%
the last inequality holding true as soon as $\frac{1}{m}\leq \varepsilon /2.$
So (\ref{Llog}) holds and the conclusion follows from Theorem \ref{ThLog}. $%
\square $

\subsection{Stochastic heat equation}

\label{sect-heat}

In this section we investigate the regularity of the law of the solution to
the stochastic heat equation introduced by Walsh in \cite{bib:[W]}. Formally
this equation is%
\begin{equation}
\partial _{t}u(t,x)=\partial _{x}^{2}u(t,x)+\sigma (u(t,x))W(t,x)+b(u(t,x))
\label{S1}
\end{equation}%
where $W$ denotes a white noise on ${\mathbb{R}}_{+}\times [0,1].$ We
consider Neumann boundary conditions that is $\partial _{x}u(t,0)=\partial
_{x}u(t,1)=0 $ and the initial condition is $u(0,x)=u_{0}(x).$ The rigorous
formulation to this equation is given by the mild form constructed as
follows. Let $G_{t}(x,y)$ be the fundamental solution to the deterministic
heat equation $\partial _{t}v(t,x)=\partial _{x}^{2}v(t,x)$ with Neumann
boundary conditions. Then $u$ satisfies%
\begin{eqnarray}
u(t,x)
&=&\int_{0}^{1}G_{t}(x,y)u_{0}(y)dy+\int_{0}^{t}\int_{0}^{1}G_{t-s}(x,y)%
\sigma (u(s,y))dW(s,y)  \label{S2} \\
&&+\int_{0}^{t}\int_{0}^{1}G_{t-s}(x,y)b(u(s,y))dyds  \notag
\end{eqnarray}%
where $dW(s,y)$ is the It\^{o} integral introduced by Walsh. The function $%
G_{t}(x,y)$ is explicitly known (see \cite{bib:[W]} or \cite{bib:[BP]}) but
here we will use just few properties that we list below (see the appendix in
\cite{bib:[BP]} for the proof). More precisely, for $0<\varepsilon <t$ we
have%
\begin{equation}
\int_{t-\varepsilon }^{t}\int_{0}^{1}G_{t-s}^{2}(x,y)dyds\leq C\varepsilon
^{1/2}  \label{INEG1}
\end{equation}%
Moreover, for $0<x_{1}<...<x_{d}<1$ there exists a constant $C$ depending on
$\min_{i=1,d}(x_{i}-x_{i-1})$ such that%
\begin{equation}
C\varepsilon ^{1/2}\geq \inf_{\left\vert \xi \right\vert
=1}\int_{t-\varepsilon }^{t}\int_{0}^{1}\left( \sum_{i=1}^{d}\xi
_{i}G_{t-s}(x_{i},y)\right) ^{2}dyds\geq C^{-1}\varepsilon ^{1/2}.
\label{INEG2}
\end{equation}%
This is an easy consequence of the inequalities (A2) and (A3) from \cite%
{bib:[BP]}.

In \cite{bib:[PZ]} one gives sufficient conditions in order to obtain the
absolute continuity of the law of $u(t,x)$ for $(t,x)\in (0,\infty )\times
\lbrack 0,1]$ and in \cite{bib:[BP]}, under appropriate hypotheses, one
obtains a $C^{\infty }$ density for the law of the vector $%
(u(t,x_{1}),...,u(t,x_{d}))$ with $(t,x_{i})\in (0,\infty )\times \{\sigma
\neq 0\},i=1,...,d.$ The aim of this section is to obtain the same type of
results but under much weaker regularity hypothesis on the coefficients. One
may first discuss the absolute continuity of the law and further, under more
regularity hypothesis on the coefficients, one may discuss the regularity of
the density. Here, in order to avoid technicalities, we restrict ourselves
to the absolute continuity property. We assume global ellipticity that is
\begin{equation}
\sigma (x)\geq c_{\sigma }>0\qquad \text{for every }x\in \lbrack 0,1].
\label{S3.1}
\end{equation}%
A local ellipticity condition may also be used but again,\ this gives more
technical complications that we want to avoid. This is somehow a benchmark
for the efficiency of the method developed in the previous sections.

We assume the following regularity hypothesis: $\sigma ,b$ are measurable
and bounded functions and there exists $h>0$ such that

\begin{equation}
\left\vert \sigma (x)-\sigma (y)\right\vert \leq \left\vert \ln \left\vert
x-y\right\vert \right\vert ^{-(2+h)},\quad \text{for every }x,y\in [0,1].
\label{S3}
\end{equation}%
This hypothesis is not sufficient in order to ensure existence and
uniqueness for the solution to (\ref{S2}) (one needs $\sigma $ and $b$ to be
globally Lipschitz continuous in order to obtain it) - so in the following
we will just consider a random field $u(t,x),(t,x)\in (0,\infty )\times
\lbrack 0,1]$ which is adapted to the filtration generated by $W$ (see Walsh
\cite{bib:[W]} for precise definitions) and which solves (\ref{S2}).

\begin{proposition}
\label{SPDE1} Suppose that (\ref{S3.1}) and (\ref{S3}) hold. Then for every $%
0<x_{1}<...<x_{d}<1$ and $T>0$, the law of the random vector $%
U=(u(T,x_{1}),...u(T,x_{d}))$ is absolutely continuous with respect to the
Lebesgue measure.
\end{proposition}

\textbf{Proof}. Given $0<\varepsilon <T$ we decompose%
\begin{equation}
u(T,x)=u_{\varepsilon }(T,x)+I_{\varepsilon }(T,x)+J_{\varepsilon }(T,x)
\label{S4}
\end{equation}%
with%
\begin{align*}
u_{\varepsilon }(T,x)
=&\int_{0}^{1}G_{t}(x,y)u_{0}(y)dy+\int_{0}^{T}\int_{0}^{1}G_{T-s}(x,y)%
\sigma (u(s\wedge (T-\varepsilon ),y))dW(s,y) \\
&+\int_{0}^{T-\varepsilon }\int_{0}^{1}G_{T-s}(x,y)b(u(s,y))dyds, \\
I_{\varepsilon }(T,x) =& \int_{T-\varepsilon
}^{T}\int_{0}^{1}G_{T-s}(x,y)(\sigma (u(s,y))-\sigma (u(s\wedge
(T-\varepsilon ),y)))dW(s,y), \\
J_{\varepsilon }(T,x) =&\int_{T-\varepsilon
}^{T}\int_{0}^{1}G_{T-s}(x,y)b(u(s,y))dyds.
\end{align*}
\textbf{Step 1}. We prove that
\begin{equation}
{\mathbb{E}}\left\vert I_{\varepsilon }(T,x)\right\vert ^{2}+{\mathbb{E}}%
\left\vert J_{\varepsilon }(T,x)\right\vert ^{2}\leq C\left\vert \ln
\varepsilon \right\vert ^{-2(2+h)}\varepsilon ^{1/2}.  \label{S5}
\end{equation}%
Let $\mu $ and $\mu_\varepsilon$ be the law of $%
U=(u(T,x_{1}),...,u(T,x_{d})) $ and $U_{\varepsilon }=(u_{\varepsilon
}(T,x_{1}),...,u_{\varepsilon }(T,x_{d}))$ respectively. Using the above
estimate one easily obtains
\begin{equation}
d_{1}(\mu ,\mu _{\varepsilon })\leq C\left\vert \ln \varepsilon \right\vert
^{-(2+h)}\varepsilon ^{1/4}.  \label{S5a}
\end{equation}%
Using the isometry property%
\begin{equation*}
{\mathbb{E}}\left\vert I_{\varepsilon }(T,x)\right\vert
^{2}=\int_{T-\varepsilon }^{T}\int_{0}^{1}G_{T-s}^{2}(x,y){\mathbb{E}}%
(\sigma (u(s,y)-\sigma (u(s\wedge (T-\varepsilon ),y)))^{2})dyds.
\end{equation*}%
We consider the set $\Lambda _{\varepsilon ,\eta }(s,y)=\{\left\vert
u(s,y)-u(s\wedge (T-\varepsilon ),y)\right\vert \leq \eta \}$ and we split
the above term as ${\mathbb{E}}\left\vert I_{\varepsilon }(T,x)\right\vert
^{2}=A_{\varepsilon ,\eta }+B_{\varepsilon ,\eta }$ with

\begin{eqnarray*}
A_{\varepsilon } &=&\int_{T-\varepsilon }^{T}\int_{0}^{1}G_{T-s}^{2}(x,y){%
\mathbb{E}}(\sigma (u(s,y)-\sigma (u(s\wedge (T-\varepsilon
),y)))^{2}1_{\Lambda _{\varepsilon ,\eta }(s,y)})dyds \\
B_{\varepsilon } &=&\int_{T-\varepsilon }^{T}\int_{0}^{1}G_{T-s}^{2}(x,y){%
\mathbb{E}}(\sigma (u(s,y)-\sigma (u(s\wedge (T-\varepsilon
),y)))^{2}1_{\Lambda _{\varepsilon ,\eta }^{c}(s,y)})dyds.
\end{eqnarray*}%
Using (\ref{S3})%
\begin{equation*}
A_{\varepsilon }\leq C(\ln \eta )^{2(2+h)}\int_{T-\varepsilon
}^{T}\int_{0}^{1}G_{T-s}^{2}(x,y)dyds\leq C\left\vert \ln \eta \right\vert
^{-2(2+h)}\varepsilon ^{1/2}
\end{equation*}%
the last inequality being a consequence of (\ref{INEG1}). Moreover, coming
back to (\ref{S2}), we have
\begin{equation*}
{\mathbb{P}}(\Lambda _{\varepsilon ,\eta }^{c}(s,y))\leq \frac{1}{\eta ^{2}}{%
\mathbb{E}}\left\vert u(s,y)-u(s\wedge (T-\varepsilon ),y)\right\vert
^{2}\leq \frac{C}{\eta ^{2}}\int_{T-\varepsilon
}^{s}\int_{0}^{1}G_{s-r}^{2}(y,z)dzdr\leq \frac{C\varepsilon ^{1/2}}{\eta
^{2}}
\end{equation*}%
so that%
\begin{equation*}
B_{\varepsilon }\leq \frac{C\varepsilon ^{1/2}}{\eta ^{2}}%
\int_{T-\varepsilon }^{T}\int_{0}^{1}G_{T-s}^{2}(x,y)dyds\leq \frac{%
C\varepsilon }{\eta ^{2}}.
\end{equation*}%
Taking $\eta =\varepsilon ^{1/16}$ we obtain%
\begin{equation*}
{\mathbb{E}}\left\vert I_{\varepsilon }(T,x)\right\vert ^{2}\leq
C(\left\vert \ln \varepsilon \right\vert ^{-2(2+h)}+\varepsilon
^{1/4})\varepsilon ^{1/2}\leq C\left\vert \ln \varepsilon \right\vert
^{-2(2+h)}\varepsilon ^{1/2}.
\end{equation*}%
We estimate now%
\begin{equation*}
\left\vert J_{\varepsilon }(T,x)\right\vert \leq \left\Vert b\right\Vert
_{\infty }\int_{T-\varepsilon }^{T}\int_{0}^{1}G_{T-s}(x,y)dyds=\left\Vert
b\right\Vert _{\infty }\varepsilon
\end{equation*}%
so (\ref{S5}) is proved.

\textbf{Step 2}. Conditionally to $\mathcal{F}_{T-\varepsilon }$ the random
vector $U_{\varepsilon }=(u_{\varepsilon }(T,x_{1}),...,u_{\varepsilon
}(T,x_{d}))$ is Gaussian of covariance matrix%
\begin{equation*}
\Sigma ^{i,j}(U_{\varepsilon })=\int_{T-\varepsilon
}^{T}\int_{0}^{1}G_{T-s}(x_{i},y)G_{T-s}(x_{j},y)\sigma ^{2}(u(s\wedge
(T-\varepsilon ),y))dyds,\quad i,j=1,...,d.
\end{equation*}%
By (\ref{INEG2})
\begin{equation*}
C\sqrt{\varepsilon }\geq \Sigma (U_{\varepsilon })\geq \frac{1}{C}\sqrt{%
\varepsilon }
\end{equation*}%
where $C$ is a constant which depends on the upper bounds of $\sigma $ and
on $c_{\sigma }.$

We use now the criterion given in Theorem \ref{ThLog} $.$ Let $%
p_{U_{\varepsilon }}$\ be the density of the law of $U_{\varepsilon }.$\
Conditionally to $\mathcal{F}_{T-\varepsilon }$ this is a Gaussian density
and the same reasoning as in the proof of (\ref{Ito8}) gives
\begin{equation*}
\left\Vert p_{U_{\varepsilon }}\right\Vert _{2m,2m,1+}\leq C(\varepsilon
^{-1/4})^{2m}\ln \frac{1}{\varepsilon }.
\end{equation*}%
So (\ref{Llog}) reads%
\begin{equation*}
\left\Vert p_{U_{\varepsilon }}\right\Vert _{2m,2m,1+}^{1/2m}d_{1}(\mu ,\mu
_{\varepsilon })\leq C\varepsilon ^{-1/4}(\ln \frac{1}{\varepsilon }%
)^{1/2m}\times \left\vert \ln \varepsilon \right\vert ^{-(2+h)}\varepsilon
^{1/4}=C\frac{1}{(\ln \frac{1}{\varepsilon })^{2+h-1/2m}}\leq C\frac{1}{(\ln
\frac{1}{\varepsilon })^{2+1/2m}}
\end{equation*}%
the last inequality being true as soon as $h>\frac{1}{m}.$ $\square $

\subsection{Piecewise deterministic Markov Processes}

\label{jumps}

In this section we deal with a jump type stochastic differential equation
which has already been considered in \cite{bib:[BCl1]}: it is an example of
piecewise deterministic Markov processes. We consider a Poisson point
process $p$ with state space $(E,\mathcal{B}(E)),$ where $E=\mathbb{R}%
^{d}\times \mathbb{R}_{+}.$ We refer to \cite{[IW]} for the notations. We
denote by $N$ the counting measure associated to $p$, that is $N([0,t)\times
A)=\#\{0\leq s<t;p_{s}\in A\}$ for $t\geq 0$ and $A\in \mathcal{B}(E)$. We
assume that the associated intensity measure is given by $\widehat{N}%
(dt,dz,du)=dt\times dz\times 1_{[0,\infty )}(u)du$ where $(z,u)\in E=\mathbb{%
R}^{d}\times \mathbb{R}_{+}.$ We are interested in the solution to the $d$
dimensional stochastic equation
\begin{equation}
X_{t}=x+\int_{0}^{t}\int_{E}c(z,X_{s-})1_{\{u<\gamma
(z,X_{s-})\}}N(ds,dz,du)+\int_{0}^{t}g(X_{s})ds.  \label{eq1}
\end{equation}%
The coefficients $c,g,\gamma $ are smooth functions (see the hypothesis $%
(H_{i}),i=0,1,2$ below). We remark that the infinitesimal generator of the
Markov process $X_{t}$ is given by
\begin{equation*}
L\psi (x)=g(x)\nabla \psi (x)+\int_{\mathbb{R}^{d}}(\psi (x+c(z,x))-\psi
(x))\gamma (z,x)dz
\end{equation*}%
See \cite{Ff} for the proof of existence and uniqueness of the solution to (%
\ref{eq1}). We will deal with two problems related to this equation.

First we give sufficient conditions in order that ${\mathbb{P}}(X_{t}(x)\in
dy)=p_{t}(x,y)dy$ where $X_{t}(x)$ is the solution to (\ref{eq1}) which
starts from $x$, so $X_{0}(x)=x.$ And we prove that, if the coefficients of
the equation are smooth, then $(x,y)\mapsto p_{t}(x,y)$ is smooth. Notice
that the methodology from \cite{Ff}, \cite{DR}, \cite{DF} and \cite{F1}
seems difficult to implement in order to prove the regularity with respect
to the initial condition $x.$ So this is the main point here.

The second result concerns convergence. In \cite{bib:[BCl1]} it is
constructed an approximation scheme which allows one to compute ${\mathbb{E}}%
(f(X_{t}(x))$ using a Monte Carlo method. And it is proved that the
convergence takes place in total variation distance. We use here the method
developed in our paper in order to prove that the density functions and
their derivatives converge as well and to estimate the error.

In \cite{bib:[BCl1]} one gives a Malliavin type approach to the equation (%
\ref{eq1}) which we recall and which we will heavily use here. We describe
first the approximation procedure. We consider a non-negative and smooth
function $\varphi :\mathbb{R}^{d}\rightarrow \mathbb{R}_{+}$ such that $%
\varphi (z)=0$ for $\left\vert z\right\vert >1$ and $\int_{\mathbb{R}%
^{d}}\varphi (z)dz=1.$ And for $M\in \mathbb{N}$ we denote $\Phi
_{M}=\varphi \ast 1_{B_{M}}$ with $B_{M}=\{z\in \mathbb{R}^{d}:\left\vert
z\right\vert <M\}.$ Then $\Phi _{M}\in C_{b}^{\infty }$ and we have $%
1_{B_{M-1}}\leq \Phi _{M}\leq 1_{B_{M+1}}.$ We denote by $X_{t}^{M}$ the
solution of the equation
\begin{equation}
X_{t}^{M}=x+\int_{0}^{t}\int_{E}c(z,X_{s-}^{M})1_{\{u<\gamma
(z,X_{s-}^{M})\}}\Phi _{M}(z)N(ds,dz,du)+\int_{0}^{t}g(X_{s}^{M})ds.
\label{eq2}
\end{equation}%
In the following we will assume that $\left\vert \gamma (z,x)\right\vert
\leq \overline{\gamma }$ for some constant $\overline{\gamma }.$ Let $%
N_{M}(ds,dz,du):=1_{B_{M+1}}(z)\times 1_{[0,2\overline{\gamma }%
]}(u)N(ds,dz,du).$ Since $\{u<\gamma (z,X_{s-}^{M})\}\subset \{u<2\overline{%
\gamma }\}$ and $\Phi _{M}(z)=0$ for $\left\vert z\right\vert >M+1,$ we may
replace $N$ by $N_{M}$ in the above equation and consequently $X_{t}^{M}$ is
solution to the equation
\begin{align*}
X_{t}^{M}& =x+\int_{0}^{t}\int_{E}c_{M}(z,X_{s-}^{M})1_{\{u<\gamma
(z,X_{s-}^{M})\}}N_{M}(ds,dz,du)+\int_{0}^{t}g(X_{s}^{M})ds,\quad \mbox{with}
\\
c_{M}(z,x)& =\Phi _{M}(z)c(z,x).
\end{align*}%
Since the intensity measure $\widehat{N}_{M}$ is finite we may represent the
random measure $N_{M}$ by a compound Poisson process. Let $\lambda _{M}=2%
\overline{\gamma }\times \mu (B_{M+1})=t^{-1}{\mathbb{E}}(N_{M}(t,E))$ (with
$\mu $ the Lebesgue measure) and let $J_{t}^{M}$ a Poisson process of
parameter $\lambda _{M}.$ We denote by $T_{k}^{M},k\in \mathbb{N}$ the jump
times of $J_{t}^{M}$. We also consider two sequences of independent random
variables $(Z_{k})_{k\in \mathbb{N}}$ in $\mathbb{R}^{d}$ and $(U_{k})_{k\in
\mathbb{N}}$ in $\mathbb{R}_{+}$ which are independent of $J^{M}$ and such
that
\begin{equation*}
Z_{k}\sim \frac{1}{\mu (B_{M+1})}1_{B_{M+1}}(z)dz\quad \mbox{and}\quad
U_{k}\sim \frac{1}{2\overline{\gamma }}1_{[0,2\overline{\gamma }]}(u)du.
\end{equation*}%
To simplify the notation, we omit the dependence on $M$ for the variables $%
(T_{k}^{M})$. Then equation $(\ref{eq2})$ may be written as
\begin{equation}
X_{t}^{M}=x+\sum_{k=1}^{J_{t}^{M}}c_{M}(Z_{k},X_{T_{k}^-}^{M})1_{(U_{k},%
\infty )}(\gamma (Z_{k},X_{T_{k}^-}^{M}))+\int_{0}^{t}g(X_{s}^{M})ds.
\label{eq3}
\end{equation}

Now $X_{t}^{M}$ is an explicit functional of the $Z_{k},k\in {\mathbb{N}}$
but, because of the indicator function, this functional is not
differentiable. In order to overcome this difficulty, following \cite%
{bib:[BCl1]}, we consider an alternative representation of the law of $%
X_{t}^{M}$. Let $z_{M}^{\ast }\in \mathbb{R}^{d}$ such that $\left\vert
z_{M}^{\ast }\right\vert =M+3$. We define %\begin{equation}\label{qM}
\begin{equation}
\begin{array}{rcl}
q_{M}(x,z) & := & \displaystyle\varphi (z-z_{M}^{\ast })\theta _{M,\gamma
}(x)+\frac{1}{2\overline{\gamma }\mu (B_{M+1})}1_{B_{M+1}}(z)\gamma (z,x),%
\mbox{ with}\smallskip \\
\theta _{M,\gamma }(x) & := & \displaystyle\frac{1}{\mu (B_{M+1})}%
\int_{\{\left\vert z\right\vert \leq M+1\}}\Big(1-\frac{1}{2\overline{\gamma
}}\gamma (z,x)\Big)dz.%
\end{array}
\label{a9}
\end{equation}%
We recall that $\varphi $ is a non-negative and smooth function with $\int
\varphi =1$ and which is null outside the unit ball. Moreover since, $0\leq
\gamma (z,x)\leq \overline{\gamma }$ one has $1\geq \theta _{M,\gamma
}(x)\geq 1/2$. By construction the function $q_{M}$ satisfies $\int
q_{M}(x,z)dz=1.$ Hence we can easily check (see \cite{bib:[BCl1]} for the
proof) that
\begin{equation}
{\mathbb{E}}(f(X_{T_{k}}^{M})\mid X_{T_{k}^{-}}^{M}=x)=\int_{{\mathbb{R}}%
^{d}}f(x+c_{M}(z,x))q_{M}(x,z)dz.  \label{condXM}
\end{equation}%
From the relation (\ref{condXM}) we construct a process $(\overline{X}%
_{t}^{M}),$ equal in law to $(X_{t}^{M}),$ in the following way. We denote
by $\Psi _{t}(x)$ the solution of $\Psi _{t}(x)=x+\int_{0}^{t}g(\Psi
_{s}(x))ds.$ We assume that the times $T_{k},k\in \mathbb{N}$ are fixed and
we consider a sequence $(z_{k})_{k\in \mathbb{N}}$ with $z_{k}\in \mathbb{R}%
^{d}.$ Then we define $x_{t},t\geq 0$ by $x_{0}=x$ and, if $x_{T_{k}}$ is
given, then
\begin{eqnarray*}
x_{t} &=&\Psi _{t-T_{k}}(x_{T_{k}})\quad T_{k}\leq t<T_{k+1}, \\
x_{T_{k+1}} &=&x_{T_{k+1}^{-}}+c_{M}(z_{k+1},x_{T_{k+1}^{-}}).
\end{eqnarray*}%
We remark that for $T_{k}\leq t<T_{k+1},x_{t}$ is a function of $%
z_{1},...,z_{k}.$ Notice also that $x_{t}$ solves the equation
\begin{equation}
x_{t}=x+\sum_{k=1}^{J_{t}^{M}}c_{M}(z_{k},x_{T_{k}^{-}})+%
\int_{0}^{t}g(x_{s})ds.  \label{a8}
\end{equation}

We consider now a sequence of random variables $(\overline{Z}_{k}),k\in
\mathbb{N}^{\ast }$ and we denote $\mathcal{G}_{k}=\sigma (T_{p},p\in
\mathbb{N})\vee \sigma (\overline{Z}_{p},p\leq k)$ and $\overline{X}%
_{t}^{M}=x_{t}(\overline{Z}_{1},...,\overline{Z}_{J_{t}^{M}}).$ We assume
that the law of $\overline{Z}_{k+1}$ conditionally on $\mathcal{G}_{k}$ is
given by
\begin{equation*}
{\mathbb{P}}(\overline{Z}_{k+1}\in dz\mid \mathcal{G}%
_{k})=q_{M}(x_{T_{k+1}^{-}}(\overline{Z}_{1},...,\overline{Z}%
_{k}),z)dz=q_{M}(\overline{X}_{T_{k+1}^{-}}^{M},z)dz.
\end{equation*}%
Clearly $\overline{X}_{t}^{M}$ satisfies the equation
\begin{equation}
\overline{X}_{t}^{M}=x+\sum_{k=1}^{J_{t}^{M}}c_{M}(\overline{Z}_{k},%
\overline{X}_{T_{k}^{-}}^{M})+\int_{0}^{t}g(\overline{X}_{s}^{M})ds.
\label{eq4}
\end{equation}%
And by (\ref{condXM}) the law of $\overline{X}_{t}^{M}$ coincides with the
law of $X_{t}^{M}.$ So now on we work with $\overline{X}_{t}^{M}$ which is a
smooth functional of $\overline{Z}_{k},k\in {\mathbb{N}}.$ But one more
difficulty remains: if $T_{1}>t$ then $\overline{X}_{t}^{M}$ is
deterministic, so this functional is not non-degenerated. In order to
contouring this last difficulty we add a small noise. We define%
\begin{equation*}
F_{t}^{M}(x)=\overline{X}_{t}^{M}(x)+\sqrt{TU_{M}}\times \Delta,\quad 0\leq
t\leq T,
\end{equation*}%
where $\overline{X}_{t}^{M}(x)$ is the solution to (\ref{eq4}) which starts
from $x,$ $\Delta $ is a standard normal random variable which is
independent of $T_{k}$ and $\overline{Z}_{k},k\in {\mathbb{N}}$ and%
\begin{equation}
U_{M}=\underline{\gamma }\int_{B_{M-1}^{c}}\underline{c}^{2}(z)dz  \label{UM}
\end{equation}%
with $\underline{\gamma }$ and $\underline{c}$ from (\ref{a1}) and (\ref{a2}%
) below. The approximation scheme for $X_{t}(x)$ is given by $F_{t}^{M}(x).$

Let us give our hypotheses.

\begin{itemize}
\item[$(H_{0})$] We assume that $\gamma ,g$ and $c$ are infinitely
differentiable functions in both variables $z$ and $x$. Moreover we assume
that $g$ and its derivatives are bounded.

\item[$(H_{1})$] There exist $\overline{\gamma }\geq \underline{\gamma }$,
such that
\begin{equation}
\overline{\gamma }\geq \gamma (z,x)\geq \underline{\gamma }\geq 0,\quad
\forall x\in \mathbb{R}^{d}  \label{a1}
\end{equation}%
and, for every $l\in {\mathbb{N}}$ there exists $\overline{\gamma }_{l}$ and
$\overline{\gamma }_{\ln ,l}$ such that for $\left\vert \alpha \right\vert
+\left\vert \beta \right\vert \leq l$%
\begin{equation}
\left\vert \partial _{x}^{\alpha }\partial _{z}^{\beta }\gamma
(x,z)\right\vert \leq \overline{\gamma }_{l},\qquad \left\vert \partial
_{x}^{\alpha }\partial _{z}^{\beta }\ln \gamma (x,z)\right\vert \leq
\overline{\gamma }_{\ln ,l}.  \label{h1}
\end{equation}

\item[$(H_{2})$] Setting, for $0<a<b$ and $r>0$,
\begin{equation*}
\underline{c}(z)=\frac{a}{1+\left\vert z\right\vert ^{r}},\qquad \overline{c}%
(z)=\frac{b}{1+\left\vert z\right\vert ^{r}},
\end{equation*}%
we assume that%
\begin{equation}
\left\Vert \nabla _{x}c\times (I+\nabla _{x}c)^{-1}(z,x)\right\Vert
+\left\vert c(z,x)\right\vert +\left\vert \partial _{z}^{\beta }\partial
_{x}^{\alpha }c(z,x)\right\vert \leq \overline{c}(z)\quad \forall z,x\in
\mathbb{R}^{d}  \label{a2}
\end{equation}%
and%
\begin{equation}
\sum_{j=1}^{d}\left\langle \partial _{z_{j}}c(z,x),\xi \right\rangle
^{2}\geq \underline{c}^{2}(z)\left\vert \xi \right\vert ^{2},\quad \forall
\xi \in \mathbb{R}^{d}.  \label{h2}
\end{equation}
\end{itemize}

\begin{remark}
The above hypotheses represent a particular case of the hypotheses from \cite%
{bib:[BCl1]}, corresponding to Example 1,ii) page 634 in that paper. More
general hypotheses may be considered (see \cite{bib:[BCl1]}) but our aim is
just to give an example in order to illustrate our method, so we restrict
ourself to this case.
\end{remark}

The basic estimate in our approach is the following:

\begin{theorem}
\label{C} Suppose that Hypotheses $(H_{i}),i=0,1,2$ hold. Consider a
function $\psi \in C_{b}^{\infty }({\mathbb{R}}^{d})$ such that $%
1_{B_{1}}\leq \psi \leq 1_{B_{2}}.$ Then for every $t,R>0,q\in {\mathbb{N}}$
and every multi indexes $\alpha ,\beta $ with $\left\vert \alpha \right\vert
+\left\vert \beta \right\vert \leq q,$ one has
\begin{equation}
\sup_{\left\vert x\right\vert \leq R,\left\vert y\right\vert \leq
R}\left\vert \partial _{x}^{\alpha }{\mathbb{E}}((\partial ^{\beta }\phi
)(F_{t}^{M}(x))\psi (F_{t}^{M}(x)-y)\right\vert \leq C\|\phi\|_\infty M^{dq}.
\label{a13}
\end{equation}%
Here $C$ is a constant which depends on $t,R,q$ but not on $M.$ In
particular the density $p_{t}^{M}(x,y)$ of the law of $F_{t}^{M}(x)$ verifies%
\begin{equation}
\sup_{\left\vert x\right\vert \leq R,\left\vert y\right\vert \leq
R}\left\vert \partial _{x}^{\alpha }\partial _{y}^{\beta
}p_{t}^{M}(x,y)\right\vert \leq CM^{d(q+d)}.  \label{a15}
\end{equation}
\end{theorem}

The above theorem is an extension of estimate (42) in Proposition 4 page 640
in \cite{bib:[BCl1]} and the proof is similar, except for one point: here we
consider derivatives $\partial _{x}^{\alpha }$ also (while in \cite%
{bib:[BCl1]} $\partial _{y}^{\beta }$ only appears). So we just sketch the
proof and focus on this supplementary difficulty.

We use an integration by parts formula based on $\overline{Z}_{k},k\in {%
\mathbb{N}}_{\ast }$ and on $\overline{Z}_{0}=\Delta $ which is constructed
as follows (we follow \cite{bib:[BCl1]}). Here $J=J_{t}^{M}$ and $T_{k}$ are
fixed, so they appear as constants. A simple functional is a random variable
of the form $F=f(\overline{Z}_{0},\overline{Z}_{1},...,\overline{Z}_{J})$
where $f$ is a smooth function. We use the weights $\pi _{k}=\Phi _{M}(%
\overline{Z}_{k}),k\in {\mathbb{N}}_{\ast },\pi _{0}=1$ and the Malliavin
derivative is defined as
\begin{equation*}
D_{k,j}=\pi _{k}\partial _{\overline{Z}_{k}^{j}}.
\end{equation*}%
For a multi index $\alpha =(\alpha _{1},...,\alpha _{q})$ with $\alpha
_{i}=(k_{i},j_{i})$ one defines the iterated derivative%
\begin{equation*}
D_{\alpha }=D_{\alpha _{q}}...D_{\alpha _{1}}.
\end{equation*}%
Then one defines the Sobolev norms:%
\begin{equation*}
\left\vert F\right\vert _{q}^{2}=\left\vert F\right\vert ^{2}+\sum_{1\leq
\left\vert \alpha \right\vert \leq q}\left\vert D_{\alpha }F\right\vert
^{2},\qquad \left\Vert F\right\Vert _{q,p}=({\mathbb{E}}(\left\vert
F\right\vert _{q}^{p}))^{1/p}.
\end{equation*}%
For $F=(F^{1},...,F^{d})$ the Malliavin covariance matrix is given by
\begin{equation*}
\sigma _{F}^{i,j}=\left\langle DF^{i},DF^{j}\right\rangle
=\sum_{k=0}^{J}\sum_{l=1}^{d}D_{k,l}F^{i}\times D_{k,l}F^{j}.
\end{equation*}

We introduce now the operator $L.$ Notice that the law of $\overline{Z}=(%
\overline{Z}_{0},\overline{Z}_{1},...,\overline{Z}_{J})$ is absolutely
continuous and has the density%
\begin{equation}
p_{J,x}(z_{0},z_{1},...,z_{J})=N(z_{0})%
\prod_{k=1}^{J}q_{M}(x_{T_{k}}(x,z_{1},...,z_{k-1}),z_{k})  \label{a10}
\end{equation}%
where $N$ is the density of the standard normal law (so of $\Delta ),$ $%
q_{M} $ is defined in (\ref{a9}) and $x_{T_{k}}(x,z_{1},...,$ $z_{k-1})$ is
the solution of (\ref{a8}) which starts from $x.$ Then we define%
\begin{equation*}
LF=\sum_{k=0}^{J}\sum_{j=1}^{d}D_{k,j}D_{k,j}F+D_{k,j}F\times D_{k,j}\ln
p_{J,x}(\overline{Z}_{k}).
\end{equation*}%
The basic duality relation is the following: for two simple functionals $F,G$%
\begin{equation*}
{\mathbb{E}}(FLG)={\mathbb{E}}(GLF)={\mathbb{E}}(\left\langle
DF,DG\right\rangle ).
\end{equation*}%
Having these objects at hand one proves the following integration by parts
formula. Let $F=(F^{1},...,F^{d})$ and $G$ be simple functionals and let $%
\beta =(\beta _{1},...,\beta _{q})\in \{1,...,d\}^{q}$ be multi-index of
length $q.$ Then for every $\phi \in C^{\infty }({\mathbb{R}}^{d})$%
\begin{equation}
{\mathbb{E}}(\partial _{\beta }\phi (F)G)={\mathbb{E}}(\phi (F)H_{\beta
}(F,G))  \label{a6}
\end{equation}%
where $H_{\beta }(F,G)$ is a random variable which verifies%
\begin{equation}
\left\Vert H_{\beta }(F,G\right\Vert _{p}\leq C\left\Vert (\det \sigma
_{F})^{-1}\right\Vert _{4p}^{3q-1}(1+\left\Vert F\right\Vert
_{q+1,4p}^{(6d+1)q})(1+\left\Vert LF\right\Vert _{q-1,4p}^{q})\left\Vert
G\right\Vert _{q,4p}.  \label{a5}
\end{equation}%
This result is proved in Theorem 2 and Theorem 3 in \cite{bib:[BCl1]}.
Before going on we need the following estimates.

\begin{lemma}
\label{MP-lemma} For every multi-index $\beta =(\beta _{1},...,\beta
_{q})\in \{1,...,d\}^{q}$ and every $p,R,T\geq 1$%
\begin{equation}
\sup_{\left\vert x\right\vert \leq R}{\mathbb{E}}(\sup_{t\leq T}\left\vert
\partial _{x}^{\beta }F_{t}^{M}(x)\right\vert _{l}^{p})\leq C  \label{a11}
\end{equation}%
and
\begin{equation}
\sup_{\left\vert x\right\vert \leq R}\left\Vert \partial _{x}^{\beta }\ln
p_{J,x}(\overline{Z})\right\Vert _{l,q}\leq CM^{d}.  \label{a12}
\end{equation}
\end{lemma}

\textbf{Proof}. The proof of (\ref{a11}) is analogous to the proof of Lemma
7 and Lemma 9 in \cite{bib:[BCl1]} so we leave it out. Let us prove (\ref%
{a12}). Notice first that%
\begin{equation*}
\partial _{x}^{\beta }\ln
p_{J,x}(z_0,z_{1},...,z_{J})=\sum_{k=1}^{J}\partial _{x}^{\beta }\ln
q_{M}(x_{T_{k}}(x,z_{1},...,z_{k-1}),z_{k}).
\end{equation*}%
On the set $\{q_{M}>0\}$ we have
\begin{eqnarray*}
&&\partial _{x}^{\beta }\ln q_{M}(x_{T_{k}}(x,z_{1},...,z_{k-1}),z_{k}) \\
&=&1_{B_{M+1}}(z_{k})\partial _{x}^{\beta }\ln \gamma
(x_{T_{k}}(x,z_{1},...,z_{k-1}),z_{k})+1_{B_{M+1}^{c}}(z_{k})\partial
_{x}^{\beta }\ln \theta _{M,\gamma }(x_{T_{k}}(x,z_{1},...,z_{k-1}).
\end{eqnarray*}%
We will use the following easy inequality: for any function $f\in C_{b}^{l}$
and every simple functional $F$ in ${\mathbb{R}}^{d}$ one has $\left\vert
f(F)\right\vert _{l}\leq C\left\Vert f\right\Vert _{l,\infty }\left\vert
F\right\vert _{l}$ where $\left\Vert f\right\Vert _{l,\infty
}=\sup_{x}\max_{\left\vert \alpha \right\vert \leq l}\left\vert \partial
^{\alpha }f(x)\right\vert .$ Notice that for every multi-index $\alpha $ one
has%
\begin{equation*}
\partial _{x}^{\beta }\theta _{M,\gamma }(x)=-\frac{1}{2\overline{\gamma }%
\mu (B_{M+1})}\int_{B_{M+1}}\partial _{x}^{\beta }\gamma (x,z)dz
\end{equation*}%
and moreover $\theta _{M,\gamma }(x)\geq 1/2.$ It follows that $\left\Vert
\ln \theta _{M,\gamma }\right\Vert _{l,\infty }\leq C\overline{\gamma }_{l}/%
\overline{\gamma }.$ One also has $\left\Vert \partial _{x}^{\beta }\ln
\gamma \right\Vert _{l,\infty }\leq \overline{\gamma }_{l+\left\vert \beta
\right\vert }$ so finally $\| \ln q_{M}(\cdot ,z)\|_{l,\infty} \leq C$ with $%
C$ a constant which depends on $\overline{\gamma },\overline{\gamma }_{l},%
\overline{\gamma }_{\ln l}.$ Then, using the above remark we obtain%
\begin{equation*}
\left\vert \partial _{x}^{\beta }\ln q_{M}(x_{T_{k}}(x,\overline{Z}_{1},...,%
\overline{Z}_{k-1}),\overline{Z}_{k})\right\vert _{l}\leq C\left\vert
F_{T_{k}}^{M}(x)\right\vert _{l}.
\end{equation*}%
Consequently%
\begin{equation*}
\left\vert \partial _{x}^{\beta }\ln p_{J,x}(\overline{Z}_{1},...,\overline{Z%
}_{J^M_t})\right\vert _{l}\leq C\sum_{k=1}^{J_{t}^{M}}\left\vert
F_{T_{k}}^{M}(x)\right\vert _{l}\leq J_{t}^{M}\times \sup_{s\leq
t}\left\vert F_{s}^{M}(x)\right\vert _{l}
\end{equation*}%
Since $({\mathbb{E}}(| J_{t}^{M}| ^{2}))^{1/2}=CM^{d}$ this, together with (%
\ref{a11}), gives%
\begin{equation*}
\left\Vert \partial _{x}^{\beta }\ln p_{J,x}(\overline{Z}_{1},...,\overline{Z%
}_{J^M_t})\right\Vert _{l,p}\leq CM^{d}.
\end{equation*}%
$\square $

\medskip

We are now ready to proceed to the

\smallskip

\textbf{Proof of Theorem \ref{C}}. In order to avoid notational
complications we just look to a particular case (the general case is
obviously similar). We assume that we are in the one dimensional case $d=1$
and $\left\vert \alpha \right\vert =\left\vert \beta \right\vert =1.$ Then
we look to
\begin{equation*}
\partial _{x}^{\alpha }{\mathbb{E}}((\partial ^{\beta }\phi
)(F_{t}^{M}(x))\psi (F_{t}^{M}(x)-y))=\partial _{x}{\mathbb{E}}(\phi
^{\prime }(F_{t}^{M}(x))\psi (F_{t}^{M}(x)-y)).
\end{equation*}%
Let $\nu (du)$ be the standard normal law and $z=(z_{1},...,z_{J}).$ Then,
with $\delta =\sqrt{TU_{M}}$ and $J=J_{t}^{M},$ we have%
\begin{eqnarray*}
&&\partial _{x}{\mathbb{E}}(\phi ^{\prime }(F_{t}^{M}(x))\psi
(F_{t}^{M}(x)-y)) \\
&=&\partial _{x}{\mathbb{E}}\int \nu (du)\int \phi ^{\prime }(\delta
u+x_{t}(x,z))\psi (\delta u+x_{t}(x,z)-y))p_{J,x}(z)dz \\
&=&I_{1}+I_{2}+I_{3}
\end{eqnarray*}%
with
\begin{eqnarray*}
I_{1} &=&{\mathbb{E}}\int \nu (du)\int \phi ^{\prime \prime }(\delta
u+x_{t}(x,z))\partial _{x}x_{t}(x,z)\psi (\delta u+x_{t}(x,z)-y))p_{J,x}(z)dz
\\
I_{2} &=&{\mathbb{E}}\int \nu (du)\int \phi ^{\prime }(\delta
u+x_{t}(x,z))\psi ^{\prime }(\delta u+x_{t}(x,z)-y))\partial
_{x}x_{t}(x,z)p_{J,x}(z)dz \\
I_{3} &=&{\mathbb{E}}\int \nu (du)\int \phi ^{\prime }(\delta
u+x_{t}(x,z))\psi (\delta u+x_{t}(x,z)-y))\partial _{x}p_{J,x}(z)dz.
\end{eqnarray*}%
We stress that $x_{t}(x,z)$ is defined as the solution of the equation (\ref%
{a8}) and so it depends on $T_{k},k\leq J_{t}^{M}.$ This is why ${\mathbb{E}}
$ appears in the previous expressions. Let us treat $I_{1}.$ Using the
integration by parts formula (\ref{a6})
\begin{eqnarray*}
I_{1} &=&{\mathbb{E}}(\phi ^{\prime \prime }(F_{t}^{M}(x))\partial
_{x}F_{t}^{M}(x)\psi (F_{t}^{M}(x)-y)) \\
&=&{\mathbb{E}}((\phi (F_{t}^{M}(x))H_{2}(F_{t}^{M}(x),\partial
_{x}F_{t}^{M}(x)\psi (F_{t}^{M}(x)-y)).
\end{eqnarray*}%
We use now some results from \cite{bib:[BCl1]}: according to Lemma 13 from
we have
\begin{equation}
\left\Vert LF_{t}^{M}(x)\right\Vert _{l,p}\leq CM;  \label{e1}
\end{equation}%
according to Lemma 9 we have
\begin{equation}
\left\Vert F_{t}^{M}(x)\right\Vert _{l,p}\leq C;  \label{e2}
\end{equation}%
Lemma 16 gives
\begin{equation}
{\mathbb{E}}((\det \sigma _{F_{t}^{M}(x)})^{-p})\leq C  \label{e3}
\end{equation}%
(notice that in Lemma 16 one asks that $2dp/t<\theta $ with $\theta $
defined in Hypothesis 3.2, iii) pg 630 in \cite{bib:[BCl1]}; but as said in
Example 1, ii) from the above paper, under our hypothesis we have $\theta
=\infty $ so our inequality holds for every $t>0$). Moreover, taking a look
to the proofs of the above results, one can see that the estimates (\ref{e1}%
),(\ref{e2}),(\ref{e3}) are uniform with respect to $x\in B_{R}.$ Then,
using (\ref{a5})
\begin{equation*}
\left\vert I_{1}\right\vert \leq C\Vert \phi \Vert _{\infty }M^{2}
\end{equation*}

and the estimate is uniform with respect to $x,y\in B_{R}.$ A similar
reasoning gives the same inequality for $I_{2}.$

We come now to $I_{3}.$ We write $\partial _{x}p_{J,x}(z)=\partial _{x}\ln
p_{J,x}(z)\times p_{J,x}(z)$ so that%
\begin{eqnarray*}
I_{3} &=&{\mathbb{E}}(\phi ^{\prime }(F_{t}^{M}(x))\psi
(F_{t}^{M}(x)-y)\partial _{x}\ln p_{J,x}(\overline{Z}_{1},...,\overline{Z}%
_{J})) \\
&=&{\mathbb{E}}((\phi (F_{t}^{M}(x))H_{1}(F_{t}^{M}(x),\psi
(F_{t}^{M}(x)-y)\partial _{x}\ln p_{J,x}(\overline{Z}_{1},...,\overline{Z}%
_{J})).
\end{eqnarray*}%
Using (\ref{a5}) and (\ref{a12}) we obtain%
\begin{equation*}
\left\vert I_{3}\right\vert \leq C\|\phi\|_\infty M^{2}.
\end{equation*}%
$\square $

\medskip

We will use the following approximation result:

\begin{lemma}
\label{MP-lemma2} Let $(H_{2})$ holds with $r>d$. For every Lipschitz
continuous function $f$ with Lipschitz constant less or equal to one, one
has
\begin{equation}
\left\vert {\mathbb{E}}(f(F_{t}^{M}(x))-{\mathbb{E}}(f(X_{t}(x))\right\vert
\leq CM^{-(r-d)}.  \label{a14}
\end{equation}%
where $C$ is a constant which is independent of $M$.
\end{lemma}

\textbf{Proof}. We have%
\begin{equation*}
\left\vert {\mathbb{E}}(f(F_{t}^{M}(x))-{\mathbb{E}}(f(\overline{X}%
_{t}^{M}(x))\right\vert \leq \sqrt{TU_{M}}\,{\mathbb{E}}(\left\vert \Delta
\right\vert )\leq CM^{-(r-d/2)},
\end{equation*}%
in which we have used $(H_{2})$ in order to estimate $U_{M}$ in (\ref{UM}).

Since the law of $\overline{X}_{t}^{M}(x)$ and $X_{t}^{M}(x)$ coincide, we
use Lemma 4 from \cite{bib:[BCl1]} and $(H_2)$. So, we obtain%
\begin{eqnarray*}
\left\vert {\mathbb{E}}(f(F_{t}^{M}(x))-{\mathbb{E}}(f(X_{t}(x))\right\vert
&\leq &CM^{-(r-d/2)}+\left\vert {\mathbb{E}}(f(X_{t}^{M}(x))-{\mathbb{E}}%
(f(X_{t}(x))\right\vert \\
&\leq &CM^{-(r-d/2)}+C\overline{\gamma }\int_{\{\left\vert z\right\vert >M\}}%
\overline{c}(z)dz \\
&\leq &CM^{-(r-d)}.
\end{eqnarray*}%
$\square $

\medskip

We are now able to present our main result.

\begin{theorem}
\label{MP-main} Assume Hypotheses $(H_i)$, $i=0,1,2$, hold. Let $q\in {%
\mathbb{N}}$ and $p>1$ be such that $d+2d(q+1+d/p_{\ast })<r$, where $r$ is
the constant in $(H_2)$. Then, for every $x\in {\mathbb{R}}^{d}$ and $t>0$
the law of $X_{t}(x)$ is absolutely continuous with respect to the Lebesgue
measure. We denote by $p_{t}(x,y)$ the density. Moreover, for every $R>0,$ $%
(x,y)\mapsto p_{t}(x,y)$ belongs to $W^{q,p}(B_{R}\times B_{R})$ and there
exists a constant $C$ (depending on $R)$ such that, for every $M\in {\mathbb{%
N}}$ and $\varepsilon >0$%
\begin{equation*}
\left\Vert p_{t}-p_{t}^{M}\right\Vert _{W^{q,p}(B_{R}\times B_{R})}\leq
\frac{C}{M^{r-d-2d(q+1+d/p_{\ast })-\varepsilon }}.
\end{equation*}
\end{theorem}

\begin{remark}
If $r>3d+2d^{2}$ then Sobolev embedding theorem ensures that $(x,y)\mapsto
p_{t}(x,y)$ is a continuous function. As a consequence, for every $x_{0}\in {%
\mathbb{R}}^{d}$ one may find $y_{0}\in {\mathbb{R}}^{d},\delta >0$ such that%
\begin{equation*}
\inf_{\left\vert y-y_{0}\right\vert \leq \delta }\inf_{\left\vert
x-x_{0}\right\vert \leq \delta }p_{t}(x,y)>0.
\end{equation*}%
This property is crucial in order to use Nummelin's splitting method in
order to prove convergence to equilibrium, see e.g. \cite{EL} , \cite{V1}
and \cite{V2}.
\end{remark}

\textbf{Proof}. We will use Theorem \ref{Conv} for the following measures.
Given $R>0$ we denote by $\Psi _{R}(x)$ a smooth function which verifies $%
1_{B_{R}}\leq \Psi _{R}\leq 1_{B_{R+1}}$ and we define
\begin{equation*}
f_{R,M}(x,y)=\Psi_R(x)\Psi_R(y)p_t^M(x,y) \quad\mbox{and}\quad
f_{R}(x,y)=\Psi_R(x)\Psi_R(y)p_t(x,y).
\end{equation*}
We note that
\begin{equation*}
\|p_t-p^M_t\|_{W^{q,p}(B_R\times B_R)} \leq \|f_R-f_{R,M}\|_{W^{q,p}({%
\mathbb{R}}^d\times {\mathbb{R}}^d)}.
\end{equation*}
We will use Theorem \ref{Conv} to estimate the term in the above r.h.s. Let
\begin{equation*}
\mu_{R,M}(dx,dy)=f_{R,M}(x,y)dxdy \quad\mbox{and}\quad
\mu_{R}(dx,dy)=f_{R}(x,y)dxdy.
\end{equation*}
For a Lipschitz continuous function with Lipschitz constant $\leq 1$, one
has
\begin{align*}
\left\vert \int gd\mu _{R}-\int gd\mu _{R}^{M}\right\vert &= \left\vert \int
\Psi_R(x)\Big({\mathbb{E}}(g(x,X_t(x))\Psi_R(X_t(x))- {\mathbb{E}}%
(g(x,X^M_t(x))\Psi_R(X^M_t(x))\Big)dx\right\vert \\
&\leq C M^{-(r-d)},
\end{align*}
in which we have used (\ref{a14}). Then, $d_{1}(\mu _{R},\mu _{R}^{M})\leq
CM^{-(r-d)}.$ By (\ref{a15}) we also have
\begin{equation*}
\|f_{R,M}\|_{2m+q,2m,p}\leq CM^{d(2m+q+d)}.
\end{equation*}
Now, we fix $m$ and we apply Theorem \ref{Conv} $i)$ with
\begin{equation*}
\alpha =\alpha (m)=\frac{r-d}{d(q+2m+d)}
\end{equation*}%
and $\eta (M)=M^{r-d}$. So, we obtain that $\mu _{R}$ is absolutely
continuous and if $f_{R}$ denotes its density, we also get
\begin{equation*}
\left\Vert f_{R}-f_{R,M}\right\Vert _{W^{q,p}({\mathbb{R}}^{d}\times {%
\mathbb{R}}^{d})}\leq C\frac 1{M^{(r-d)\theta}} \quad \mbox{with}\quad
\theta=\frac 1\alpha\wedge \Big(1-\frac{q+1+d/p_{\ast }}{\alpha m}\Big).
\end{equation*}%
Since $\lim_{m}m\alpha (m)=\frac{r-d}{2d}$ we obtain%
\begin{equation*}
(r-d)(1-\frac{q+1+d/p_{\ast }}{\alpha m})\rightarrow r-d-2d(q+1+d/p_{\ast })
\end{equation*}%
and%
\begin{equation*}
\frac{r-d}{\alpha }=d(q+2m+d)\rightarrow \infty
\end{equation*}%
So, taking $m$ sufficiently large we obtain, for each $\varepsilon >0$%
\begin{equation*}
\left\Vert f_{R}-f_{R,M}\right\Vert _{W^{q,p}({\mathbb{R}}^{d}\times {%
\mathbb{R}}^{d})}\leq \frac{C}{M^{r-d-2d(q+1+d/p_{\ast })-\varepsilon }}.
\end{equation*}%
$\square $

\begin{corollary}
\label{cor-MP-main} Suppose that $r\geq 3d+2d^{2}$ and set $%
k=\lfloor(r-3d-2d^{2})/2d\rfloor$. Then for every $R>0$ and every $%
\varepsilon >0$ there exists a constant $C_{R,\varepsilon }\geq 1$ such that
for every multi-indexes $\alpha ,\beta $ with $\left\vert \alpha \right\vert
+\left\vert \beta \right\vert \leq k$
\begin{equation*}
\sup_{\left\vert x\right\vert \leq R,\left\vert y\right\vert \leq
R}\left\vert \partial _{x}^{\alpha }\partial _{y}^{\beta
}p_{t}(x,y)-\partial _{x}^{\alpha }\partial _{y}^{\beta
}p_{t}^{M}(x,y)\right\vert \leq \frac{C_{R,\varepsilon }}{%
M^{r-d-2d(q+1+d/p_{\ast })-\varepsilon }}.
\end{equation*}
\end{corollary}

\textbf{Proof}. We take $p>1$ very close to $1$ (so that $p_{\ast }$ is very
large) and%
\begin{equation*}
q=\frac{r-d}{2d}-1-\frac{d}{p_{\ast }},\qquad k=\Big\lfloor q-\frac{d}{p}%
\Big\rfloor=\Big\lfloor \frac{r-3d-d^{2}}{2d}\Big\rfloor.
\end{equation*}%
Then Sobolev embedding theorem says that for $\left\vert \alpha \right\vert
+\left\vert \beta \right\vert \leq k$
\begin{equation*}
\sup_{\left\vert x\right\vert \leq R,\left\vert y\right\vert \leq
R}\left\vert \partial _{x}^{\alpha }\partial _{y}^{\beta }f(x,y)\right\vert
\leq C_{R}\left\Vert f\right\Vert _{W^{q,p}(B_{R}\times B_{R})}
\end{equation*}%
and we are done. $\square $

\appendix

\section{Hermite expansions and density estimates}

\label{sect-Hermite} The aim of this section is to give the proof of
Proposition \ref{8bis}. We recall that for $\mu \in \mathcal{M}$ and $\mu
_{n}(x)=f_{n}(x)dx,n\in {\mathbb{N}}$,
\begin{equation*}
\pi _{q,k,m,\mathbf{e}}(\mu ,(\mu _{n})_{n})=\sum_{n=0}^{\infty
}2^{n(q+k)}\beta _{\mathbf{e}}(2^{nd})d_{k}(\mu ,\mu
_{n})+\sum_{n=0}^{\infty }\frac{1}{2^{2nm}}\left\Vert f_{n}\right\Vert
_{2m+q,2m,\mathbf{e}}.
\end{equation*}

Our proposal for this section is to prove the following

\begin{proposition}
\label{8}Let $q,k\in {\mathbb{N}},m\in {\mathbb{N}}_{\ast }$ and $\mathbf{e}%
\in \mathcal{E}.$ There exists a universal constant $C$ (depending on $%
q,k,m,d$ and $e$) such that for every $f,f_{n}\in C^{2m+q}({\mathbb{R}}%
^{d}),n\in {\mathbb{N}}$, one has
\begin{equation}
\begin{array}{ll}
& \left\Vert f\right\Vert _{q,\mathbf{e}}\leq C\pi _{q,k,m,\mathbf{e}}(\mu
,(\mu _{n})_{n}).%
\end{array}
\label{Oo10}
\end{equation}%
where $\mu (x)=f(x)dx$ and $\mu _{n}(x)=f_{n}(x)dx.$
\end{proposition}

The proof of Proposition \ref{8} will follow from the next results and
properties of Hermite polynomials, so we postpone it at the end of this
section.

\medskip

We begin with a review of some basic properties of Hermite polynomials and
functions. The Hermite polynomials on ${\mathbb{R}}$ are defined by%
\begin{equation*}
H_{n}(t)=(-1)^{n}e^{t^{2}}\frac{d^{n}}{dt}e^{-t^{2}},\quad n=0,1,...
\end{equation*}%
They are orthogonal with respect to $e^{-t^{2}}dt.$ We denote the $L^{2}$
normalized Hermite functions by
\begin{equation*}
h_{n}(t)=(2^{n}n!\sqrt{\pi })^{-1/2}H_{n}(t)e^{-t^{2}/2}
\end{equation*}%
and we have%
\begin{equation*}
\int_{{\mathbb{R}}}h_{n}(t)h_{m}(t)dt=(2^{n}n!\sqrt{\pi })^{-1}\int_{{%
\mathbb{R}}}H_{n}(t)H_{m}(t)e^{-t^{2}}dt=\delta _{n,m}.
\end{equation*}%
The Hermite functions form an orthonormal basis in $L^{2}({\mathbb{R}}).$
For a multi index $\alpha =(\alpha _{1},...,\alpha _{d})\in {\mathbb{N}}^{d}$
we define the $d$-dimensional Hermite function
\begin{equation*}
\mathcal{H}_{\alpha }(x):=\prod_{i=1}^{d}h_{\alpha _{i}}(x_{i}),\quad
x=(x_{1},...,x_{d}).
\end{equation*}%
The $d$-dimensional Hermite functions form an orthonormal basis in $L^{2}({%
\mathbb{R}}^{d}).$ This corresponds to the chaos decomposition in dimension $%
d$ (but the notation we gave above is slightly different from the one used
in probability; see \cite{bib:[N]}, \cite{bib:[Sa]} and \cite{bib:[M]},
where Hermite polynomials are used. One may come back by a renormalization).
The Hermite functions are the eigenvectors of the Hermite operator $%
D=-\Delta +\left\vert x\right\vert ^{2}$, $\Delta $ denoting the Laplace
operator, and one has%
\begin{equation}
D\mathcal{H}_{\alpha }=(2\left\vert \alpha \right\vert +d)\mathcal{H}%
_{\alpha }\quad \mbox{with}\quad \left\vert \alpha \right\vert =\alpha
_{1}+...+\alpha _{d}.  \label{Her1}
\end{equation}%
We denote $W_{n}=\mathrm{Span}\{\mathcal{H}_{\alpha }:\left\vert \alpha
\right\vert =n\}$ and we have $L^{2}({\mathbb{R}}^{d})=\oplus _{n=0}^{\infty
}W_{n}$.

For a function $\Phi :{\mathbb{R}}^{d}\times {\mathbb{R}}^{d}\rightarrow {%
\mathbb{R}}$ and a function $f:{\mathbb{R}}^{d}\rightarrow {\mathbb{R}}$ we
use the notation%
\begin{equation*}
\Phi \diamond f(x)=\int_{{\mathbb{R}}^{d}}\Phi (x,y)f(y)dy.
\end{equation*}%
We denote by $J_{n}$ the orthogonal projection on $W_{n}$ and we have
\begin{equation}
J_{n}v(x)=\bar{\mathcal{H}}_{n}\diamond v(x)\quad \mbox{with}\quad \bar{%
\mathcal{H}}_{n}(x,y):=\sum_{\left\vert \alpha \right\vert =n}\mathcal{H}%
_{\alpha }(x)\mathcal{H}_{\alpha }(y).  \label{Her2}
\end{equation}%
Moreover, we consider a function $a:{\mathbb{R}}_{+}\rightarrow {\mathbb{R}}$
whose support is included in $[\frac{1}{4},4]$\ and we define%
\begin{equation*}
\bar{\mathcal{H}}_{n}^{a}(x,y)=\sum_{j=0}^{\infty }a\Big(\frac{j}{4^{n}}\Big)%
\bar{\mathcal{H}}_{j}(x,y) =\sum_{j=4^{n-1}+1}^{4^{n+1}-1}a\Big(\frac{j}{%
4^{n}}\Big)\bar{\mathcal{H}}_{j}(x,y),\quad x,y\in {\mathbb{R}}^d,
\end{equation*}%
the last equality being a consequence of the support property of the
function $a.$

The following estimate is a crucial point in our approach. It has been
proved in \cite{bib:[E]}, \cite{bib:[Dz]} and then in \cite{bib:[PY]}. We
refer to Corollary 2.3, inequality (2.17), in \cite{bib:[PY]} (we thank to
G. Kerkyacharian who signaled us this paper).

\begin{theorem}
\label{5}Let $a:{\mathbb{R}}_{+}\rightarrow {\mathbb{R}}_{+}$ be a non
negative $C^{\infty }$ function with the support included in $[\frac{1}{4}%
,4].$ We denote $\left\Vert a\right\Vert _{l}=\sum_{i=0}^{l}\sup_{t\geq
0}\left\vert a^{(i)}(t)\right\vert .$ For every multi-index $\alpha $ and
every $k\in {\mathbb{N}}$ there exists a constant $C_{k}$ (depending on $%
k,\alpha ,d)$ such that for every $n\in {\mathbb{N}}$ and every $x,y\in {%
\mathbb{R}}^{d}$%
\begin{equation}
\left\vert \frac{\partial ^{\left\vert \alpha \right\vert }}{\partial
x^{\alpha }}\bar{\mathcal{H}}_{n}^{a}(x,y)\right\vert \leq C_{k}\left\Vert
a\right\Vert _{k}\frac{2^{n(\left\vert \alpha \right\vert +d)}}{%
(1+2^{n}\left\vert x-y\right\vert )^{k}}.  \label{Her3}
\end{equation}
\end{theorem}

Following the ideas in \cite{bib:[PY]} we consider a function $a:{\mathbb{R}}%
_{+}\rightarrow {\mathbb{R}}_{+} $ of class $C_{b}^{\infty }$ with the
support included in $[\frac{1}{4},4]$ and such that $a(t)+a(4t)=1$ for $t\in
\lbrack \frac{1}{4},1].$ We may construct $a$ in the following way: we take
a function $a:[0,1]\rightarrow {\mathbb{R}}_{+}$ with $a(t)=0$ for $t\leq
\frac{1}{4}$ and $a(1)=1.$ We may choose $a$ such that $a^{(l)}(\frac{1}{4}%
)=a^{(l)}(1-)=0 $ for every $l\in {\mathbb{N}}$. Then we define $a(t)=1-a(%
\frac{t}{4})$ for $t\in \lbrack 1,4]$ and $a(t)=0$ for $t\geq 4.$ This is
the function we will use in the following$.$ Notice that $\ a$ has the
property:%
\begin{equation}
\sum_{n=0}^{\infty }a\Big(\frac{t}{4^{n}}\Big)=1\quad \forall t\geq 1.
\label{Her0}
\end{equation}%
In order to check the above equality we fix $n_{t}$ such that $%
4^{n_{t}-1}\leq t<4^{n_{t}}$ and we notice that $a(\frac{t}{4^{n}})=0$ if $%
n\notin \{n_{t}-1,n_{t}\}.$ So $\sum_{n=0}^{\infty }a(\frac{t}{4^{n}}%
)=a(4s)+a(s)=1$ with $s=t/4^{n_{t}}\in \lbrack \frac{1}{4},1].$ In the
following we fix a function $a$ and the constants in our estimates will
depend on $\left\Vert a\right\Vert _{l}$ for some fixed $l.$ Using this
function we obtain the following representation formula:

\begin{proposition}
\label{6}For every $f\in L^{2}({\mathbb{R}}^{d})$
\begin{equation*}
f=\sum_{n=0}^{\infty }\bar{\mathcal{H}}_{n}^{a}\diamond f
\end{equation*}%
the series being convergent in $L^{2}({\mathbb{R}}^{d}).$
\end{proposition}

\textbf{Proof}. We fix $N$ and we denote
\begin{equation*}
S_{N}^{a}=\sum_{n=1}^{N}\bar{\mathcal{H}}_{n}^{a}\diamond f,\quad
S_{N}=\sum_{j=1}^{4^{N}}\bar{\mathcal{H}}_{j}\diamond f\quad and\quad
R_{N}^{a}=\sum_{j=4^{N}+1}^{4^{N+1}}(\bar{\mathcal{H}}_{j}\diamond f)a\Big(%
\frac{j}{4^{N+1}}\Big).
\end{equation*}%
Let $j\leq 4^{N+1}.$ For $n\geq N+2$ one has $a(\frac{j}{4^{n}})=0.$ So
using (\ref{Her0}) we obtain $\sum_{n=1}^{N}a(\frac{j}{4^{n}}%
)=\sum_{n=1}^{\infty }a(\frac{j}{4^{n}})-a(\frac{j}{4^{N+1}})=1-a(\frac{j}{%
4^{N+1}}).$ And for $j\leq 4^{N}$ one has $a(\frac{j}{4^{N+1}})=0.$\ It
follows that%
\begin{eqnarray*}
S_{N}^{a} &=&\sum_{n=1}^{N}\sum_{j=0}^{\infty }a\Big(\frac{j}{4^{n}}\Big)%
\bar{\mathcal{H}}_{j}\diamond f =\sum_{n=1}^{N}\sum_{j=0}^{4^{N+1}}a\Big(%
\frac{j}{4^{n}}\Big)\bar{\mathcal{H}}_{j}\diamond f =\sum_{j=0}^{4^{N+1}}(%
\bar{\mathcal{H}}_{j}\diamond f)\sum_{n=1}^{N}a\Big(\frac{j}{4^{n}}\Big) \\
&=&\sum_{j=0}^{4^{N+1}}\bar{\mathcal{H}}_{j}\diamond
f-\sum_{j=4^{N}+1}^{4^{N+1}}(\bar{\mathcal{H}}_{j}\diamond f)a\Big(\frac{j}{%
4^{N+1}}\Big)=S_{N+1}-R_{N}^{a}.
\end{eqnarray*}%
One has $S_{N}\rightarrow f$ in $L^{2}$ and $\left\Vert R_{N}^{a}\right\Vert
_{2}\leq \left\Vert a\right\Vert _{\infty
}\sum_{j=4^{N}+1}^{4^{N+1}}\left\Vert \bar{\mathcal{H}}_{j}\diamond
f\right\Vert _{2}\rightarrow 0$ so the proof is completed. $\square $

\smallskip

We will need the following lemma concerning properties of the Luxembourg
norms.

\begin{lemma}
Let $\rho \geq 0$ be a measurable function. Then for every measurable
function $f$ one has
\begin{equation}
\left\Vert \rho \ast f\right\Vert _{\mathbf{e}}\leq \left\Vert \rho
\right\Vert _{1}\left\Vert f\right\Vert _{\mathbf{e}}.  \label{Oo2}
\end{equation}
\end{lemma}

\textbf{Proof}. Let $c=m\left\Vert f\right\Vert _{\mathbf{e}}$ with $%
m=\left\Vert \rho \right\Vert _{1}=\int \rho (x-y)dy.$ Since $\mathbf{e}$ is
convex we obtain%
\begin{eqnarray*}
\int \mathbf{e}\Big(\frac{1}{c}(\rho \ast f)(x)\Big)dx &=&\int \mathbf{e}%
\Big(\int \frac{\rho (x-y)}{m}\times \frac{m}{c}f(y)dy\Big)dx \\
&\leq &\int dx\int \frac{\rho (x-y)}{m}\times \mathbf{e}\Big(\frac{m}{c}f(y)%
\Big)dy \\
&=&\int \mathbf{e}\Big(\frac{m}{c}f(y)\Big)\int \frac{\rho (x-y)}{m}%
dxdy=\int \mathbf{e}\Big(\frac{m}{c}f(y)\Big)dy \\
&=&\int \mathbf{e}\Big(\frac{1}{\left\Vert f\right\Vert _{\mathbf{e}}}f(y)%
\Big)dy\leq 1
\end{eqnarray*}%
and this means that $\left\Vert \rho \ast f\right\Vert _{\mathbf{e}}\leq
c=\left\Vert \rho \right\Vert _{1}\left\Vert f\right\Vert _{\mathbf{e}}.$ $%
\square $

\begin{lemma}
Let $\mathbf{e}\in \mathcal{E}$ and $\rho _{n,p}(z)=(1+2^{n}\left\vert
z\right\vert )^{-p}$, with $p>d$. There exists a constant $C_{p}$ depending
on $p$ and $d$ such that
\begin{equation}
\left\Vert \rho _{n,p}\right\Vert _{\mathbf{e}}\leq \frac{1}{\mathbf{e}^{-1}(%
\frac{1}{C_{p}}2^{nd})}.  \label{Oo1}
\end{equation}%
In particular, for $p=d+1$ there exists a constant $C$ depending on $d$ and
on the doubling constant of $\mathbf{e}$ such that (with $\phi _{\mathbf{e}}$
defined in (\ref{O3}))%
\begin{equation}
\left\Vert \rho _{n,d+1}\right\Vert _{\mathbf{e}}\leq \frac{C}{\mathbf{e}%
^{-1}(2^{nd})}=C2^{-nd}\beta _{\mathbf{e}}(2^{nd})=C\phi _{\mathbf{e}}(\frac{%
1}{2^{nd}}).  \label{Oo1a}
\end{equation}
\end{lemma}

\textbf{Proof}. Let $c>0$. By passing in polar coordinates and by using the
change of variable $s=2^{n}r$, we obtain%
\begin{eqnarray*}
\int_{{\mathbb{R}}^{d}}\mathbf{e}\Big(\frac{1}{c}\rho _{n,p}(z)\Big)dz
&=&A_{d}\int_{0}^{\infty }r^{d-1}\mathbf{e}\Big(\frac{1}{c}\times \frac{1}{%
(1+2^{n}r)^{p}}\Big)dr \\
&=&2^{-nd}A_{d}\int_{0}^{\infty }s^{d-1}\mathbf{e}\Big(\frac{1}{c}\times
\frac{1}{(1+s)^{p}}\Big)ds
\end{eqnarray*}%
where $A_{d}$ is the surface of the unit sphere in ${\mathbb{R}}^{d}.$ Using
the property (\ref{Y}) ii) we upper bound the above term by%
\begin{equation*}
2^{-nd}\mathbf{e}\Big(\frac{1}{c}\Big)A_{d}\int_{0}^{\infty }s^{d-1}\times
\frac{1}{(1+s)^{p}}ds=C_{p}2^{-nd}\mathbf{e}\Big(\frac{1}{c}\Big).
\end{equation*}%
In order to prove that $\left\Vert \rho _{n,p}\right\Vert _{\mathbf{e}}\leq
c $ we have to check that $\int_{{\mathbb{R}}^{d}}\mathbf{e}(\frac{1}{c}\rho
_{n,p}(z))dz\leq 1.$ In view of the above inequalities it suffices that $%
\mathbf{e}(\frac{1}{c})\leq 2^{nd}/C_{p}$ that is $c\geq 1/\mathbf{e}%
^{-1}(2^{nd}/C_{p}).$ $\square $

\begin{proposition}
\label{7a} Let $\mathbf{e}\in\mathcal{E}$ and $\mathbf{e}_*$ be the
conjugate of $\mathbf{e}$. Set $\alpha $ as a multi index.

\smallskip

$i)$ There exists a universal constant $C$ (depending on $\alpha ,d$ and $%
\mathbf{e}$) such that%
\begin{equation}  \label{Oo3}
\begin{array}{ll}
a)\quad & \left\Vert \partial _{\alpha }\bar{\mathcal{H}}_{n}^{a}\diamond
f\right\Vert _{\mathbf{e}} \leq C\left\Vert a\right\Vert _{d+1}\times
2^{n\left\vert \alpha \right\vert }\left\Vert f\right\Vert _{\mathbf{e}},
\smallskip \\
b)\quad & \left\Vert \partial _{\alpha }\bar{\mathcal{H}}_{n}^{a}\diamond
f\right\Vert _{\infty } \leq C\left\Vert a\right\Vert _{d+1}\times
2^{n\left\vert \alpha \right\vert }\beta _{\mathbf{e}}(2^{nd})\left\Vert
f\right\Vert _{\mathbf{e}_{\ast }}%
\end{array}%
\end{equation}

$ii)$ Let $m\in {\mathbb{N}}_{\ast }.$ There exists a universal constant $C$
(depending on $\alpha ,m,d$ and $\mathbf{e}$) such that
\begin{equation}
\left\Vert \bar{\mathcal{H}}_{n}^{a}\diamond \partial _{\alpha }f\right\Vert
_{ \mathbf{e}}\leq \frac{C\left\Vert a\right\Vert _{d+1}^{2}}{4^{nm}}%
\left\Vert f\right\Vert _{2m+\left\vert \alpha \right\vert ,2m,\mathbf{e}}
\label{Oo4}
\end{equation}

$iii)$ Let $k\in {\mathbb{N}}.$ There exists a universal constant $C$
(depending on $\alpha ,k,d$ and $e)$ such that
\begin{equation}
\left\Vert \bar{\mathcal{H}}_{n}^{a}\diamond \partial _{\alpha
}(f-g)\right\Vert _{\mathbf{e}}\leq C\left\Vert a\right\Vert _{d+1}\times
2^{n(\left\vert \alpha \right\vert +k)}\beta (2^{nd})d_{k}(\mu _{f},\mu _{g})
\label{Oo5}
\end{equation}
\end{proposition}

\textbf{Proof}. $i)$ By using (\ref{Her3}) with $k=d+1$ we get
\begin{equation}
\left\vert \partial _{\alpha }\bar{\mathcal{H}}_{n}^{a}\diamond
f(x)\right\vert \leq C2^{n(\left\vert \alpha \right\vert +d)}\left\Vert
a\right\Vert _{d+1}\int \rho _{n,d+1}(x-y)\left\vert f(y)\right\vert dy.
\label{Oo5'}
\end{equation}%
Since $\mathbf{e}$ is symmetric, i.e. $\mathbf{e}(|x|)=\mathbf{e}(x)$, one
has $\Vert f\Vert _{\mathbf{e}}=\Vert |f|\Vert _{\mathbf{e}}.$ Moreover, if $%
0\leq f(x)\leq g(x)$ then $\Vert f\Vert _{\mathbf{e}}\leq \Vert g\Vert _{%
\mathbf{e}}.$ Using these properties in addition to (\ref{Oo5'}) and (\ref%
{Oo2}), we obtain
\begin{eqnarray*}
\left\Vert \partial _{\alpha }\bar{\mathcal{H}}_{n}^{a}\diamond f\right\Vert
_{\mathbf{e}} &=&\left\Vert \left\vert \partial _{\alpha }\bar{\mathcal{H}}%
_{n}^{a}\diamond f\right\vert \right\Vert _{\mathbf{e}}\leq C2^{n(\left\vert
\alpha \right\vert +d)}\left\Vert a\right\Vert _{d+1}\left\Vert \rho
_{n,d+1}\ast \left\vert f\right\vert \right\Vert _{\mathbf{e}} \\
&\leq &C2^{n(\left\vert \alpha \right\vert +d)}\left\Vert a\right\Vert
_{d+1}\left\Vert \rho _{n,d+1}\right\Vert _{1}\left\Vert \left\vert
f\right\vert \right\Vert _{\mathbf{e}}.
\end{eqnarray*}%
Using (\ref{Oo1a}) with $\mathbf{e}(x)=\left\vert x\right\vert $ we obtain $%
\left\Vert \rho _{n,d+1}\right\Vert _{1}\leq C/2^{nd}.$ So we conclude that
\begin{equation*}
\left\Vert \partial _{\alpha }\bar{\mathcal{H}}_{n}^{a}\diamond f\right\Vert
_{\mathbf{e}}\leq C\left\Vert a\right\Vert _{d+1}2^{n\left\vert \alpha
\right\vert }\left\Vert \left\vert f\right\vert \right\Vert _{\mathbf{e}}
\end{equation*}%
so $a)$ is proved. Again by (\ref{Oo5'})
\begin{align*}
\left\vert \partial _{\alpha }\bar{\mathcal{H}}_{n}^{a}\diamond
f(x)\right\vert & \leq C\left\Vert a\right\Vert _{d+1}2^{n(\left\vert \alpha
\right\vert +d)}\int \rho _{n,d+1}(x-y)\left\vert f(y)\right\vert dy \\
& \leq C\left\Vert a\right\Vert _{d+1}2^{n(\left\vert \alpha \right\vert
+d)}\left\Vert \rho _{n,d+1}\right\Vert _{\mathbf{e}}\left\Vert f\right\Vert
_{\mathbf{e}_{\ast }},
\end{align*}%
the second inequality being a consequence of the H\"{o}lder inequality %
\eqref{O2}. Using (\ref{Oo1a}), $b)$ is proved as well.

\smallskip

$ii)$ We define the functions $a_{m}(t)=a(t)t^{-m}.$ Since $a(t)=0$ for $%
t\leq \frac{1}{4}$ and for $t\geq 4$ we have $\left\Vert a_{m}\right\Vert
_{d+1}\leq C_{m,d}\left\Vert a\right\Vert _{d+1}.$ Moreover $D\bar{\mathcal{H%
}}_{j}\diamond v=(2j+d)\bar{\mathcal{H}}_{j}\diamond v$ so we obtain%
\begin{equation*}
\bar{\mathcal{H}}_{j}\diamond v=\frac{1}{2j}(D-d)\bar{\mathcal{H}}%
_{j}\diamond v.
\end{equation*}%
We denote $L_{m,\alpha }=(D-d)^{m}\partial _{\alpha }$ and we notice that $%
L_{m,\alpha }=\sum_{\left\vert \beta \right\vert \leq 2m}\sum_{\left\vert
\gamma \right\vert \leq 2m+\left\vert \alpha \right\vert }c_{\beta ,\gamma
}x^{\beta }\partial _{\gamma }$ where $c_{\beta ,\gamma }$ are universal
constants. It follows that there exists some universal constant $C$ such
that
\begin{equation}
\left\Vert L_{m,\alpha }f\right\Vert _{\mathbf{e}}\leq C\left\Vert
f\right\Vert _{2m+\left\vert \alpha \right\vert ,2m,\mathbf{e}}.  \label{Oo9}
\end{equation}%
We take now $v\in L^{\mathbf{e}_{\ast }}$ and we write
\begin{eqnarray*}
\left\langle v,\bar{\mathcal{H}}_{n}^{a}\diamond (\partial _{\alpha
}f)\right\rangle &=&\left\langle \bar{\mathcal{H}}_{n}^{a}\diamond
v,\partial _{\alpha }f\right\rangle =\sum_{j=0}^{\infty }a\Big(\frac{j}{4^{n}%
}\Big)\left\langle \bar{\mathcal{H}}_{j}\diamond v,\partial _{\alpha
}f\right\rangle \\
&=&\sum_{j=1}^{\infty }a\Big(\frac{j}{4^{n}}\Big)\frac{1}{(2j)^{m}}%
\left\langle (D-d)^{m}\bar{\mathcal{H}}_{j}\diamond v,\partial _{\alpha
}f\right\rangle \\
&=&\frac{1}{2^{m}}\times \frac{1}{4^{nm}}\sum_{j=1}^{\infty }a_{m}(\frac{j}{%
4^{n}})\left\langle \bar{\mathcal{H}}_{j}\diamond v,L_{m,\alpha
}f\right\rangle \\
&=&\frac{1}{2^{m}}\times \frac{1}{4^{nm}}\left\langle \bar{\mathcal{H}}%
_{n}^{a_{m}}\diamond v,L_{m,\alpha }f\right\rangle .
\end{eqnarray*}%
By using the decomposition in Proposition \ref{6}, we write $L_{m,\alpha
}f=\sum_{j=0}^{\infty }\bar{\mathcal{H}}_{j}^{a}\diamond L_{m,\alpha }f.$
For $\left\vert j-n\right\vert \geq 2,$ by the support property of $a,$ one
has $a(\frac{k}{4^{n}})a(\frac{k}{4^{j}})=0$ for every $k\in {\mathbb{N}}$.
One also has $\langle \mathcal{H}_{\alpha }\diamond v,\mathcal{H}_{\beta
}\diamond L_{m,\alpha }f\rangle$ $=0$ if $\left\vert \alpha \right\vert \neq
\left\vert \beta \right\vert .$ Then a straightforward decomposition gives $%
\langle \bar{\mathcal{H}}_{n}^{a_{m}}\diamond v,\bar{\mathcal{H}}%
_{j}^{a}\diamond L_{m,\alpha }f\rangle =0.$ So using H\"{o}lder's inequality%
\begin{eqnarray*}
\left\vert \left\langle v,\bar{\mathcal{H}}_{n}^{a}\diamond (\partial
_{\alpha }f)\right\rangle \right\vert &\leq &\frac{1}{2^{m}}\times \frac{1}{%
4^{nm}}\sum_{j=n-1}^{n+1}\left\vert \left\langle \bar{\mathcal{H}}%
_{n}^{a_{m}}\diamond v,\bar{\mathcal{H}}_{j}^{a}\diamond L_{m,\alpha
}f\right\rangle \right\vert \\
&\leq &\frac{1}{2^{m}}\times \frac{1}{4^{nm}}\sum_{j=n-1}^{n+1}\left\Vert
\bar{\mathcal{H}}_{n}^{a_{m}}\diamond v\right\Vert _{\mathbf{e}_{\ast
}}\left\Vert \bar{\mathcal{H}}_{j}^{a}\diamond L_{m,\alpha }f\right\Vert _{%
\mathbf{e}}.
\end{eqnarray*}%
Using point $i)$ $a)$ with $\alpha $ equal to the void index, we obtain $%
\left\Vert \bar{\mathcal{H}}_{n}^{a_{m}}\diamond v\right\Vert _{\mathbf{e}%
_{\ast }}\leq C\left\Vert a_{m}\right\Vert _{d+1}\left\Vert v\right\Vert _{%
\mathbf{e}_{\ast }}\leq C\times C_{m,d}\left\Vert a\right\Vert
_{d+1}\left\Vert v\right\Vert _{\mathbf{e}_{\ast }}$. Moreover, we have $%
\left\Vert \bar{\mathcal{H}}_{j}^{a}\diamond L_{m,\alpha }f\right\Vert _{%
\mathbf{e}}\leq C\left\Vert a\right\Vert _{d+1}\left\Vert L_{m,\alpha
}f\right\Vert _{\mathbf{e}}\leq C\left\Vert a\right\Vert _{d+1}\times $ $%
\times \left\Vert f\right\Vert _{2m+\left\vert \alpha \right\vert ,2m,%
\mathbf{e}},$ the last inequality being a consequence of (\ref{Oo9}). We
obtain%
\begin{equation*}
\left\vert \left\langle v,\bar{\mathcal{H}}_{n}^{a}\diamond (\partial
_{\alpha }f)\right\rangle \right\vert \leq \frac{C\left\Vert a\right\Vert
_{d+1}^{2}}{4^{nm}}\left\Vert v\right\Vert _{\mathbf{e}_{\ast }}\left\Vert
f\right\Vert _{2m+\left\vert \alpha \right\vert ,2m,\mathbf{e}}
\end{equation*}%
and, since $L^{\mathbf{e}}$ is reflexive, (\ref{Oo4}) is proved.

\smallskip

$iii)$ We write
\begin{eqnarray*}
\left\vert \left\langle v,\bar{\mathcal{H}}_{n}^{a}\diamond (\partial
_{\alpha }(f-g))\right\rangle \right\vert &=&\left\vert \left\langle \bar{%
\mathcal{H}}_{n}^{a}\diamond v,\partial _{\alpha }(f-g)\right\rangle
\right\vert =\left\vert \left\langle \partial _{\alpha }\bar{\mathcal{H}}%
_{n}^{a}\diamond v,f-g)\right\rangle \right\vert \\
&=&\left\vert \int \partial _{\alpha }\bar{\mathcal{H}}_{n}^{a}\diamond
vd\mu _{f}-\int \partial _{\alpha }\bar{\mathcal{H}}_{n}^{a}\diamond vd\mu
_{g}\right\vert .
\end{eqnarray*}%
We use the definition of $d_{k}$ and (\ref{Oo3}) $b)$ and we obtain%
\begin{eqnarray*}
&&\left\vert \int \partial _{\alpha }\bar{\mathcal{H}}_{n}^{a}\diamond vd\mu
_{f}-\int \partial _{\alpha }\bar{\mathcal{H}}_{n}^{a}\diamond vd\mu
_{g}\right\vert \leq \left\Vert \partial _{\alpha }\bar{\mathcal{H}}%
_{n}^{a}\diamond v\right\Vert _{k,\infty }d_{k}(\mu _{f},\mu _{g}) \\
&\leq &\left\Vert \bar{\mathcal{H}}_{n}^{a}\diamond v\right\Vert
_{k+\left\vert \alpha \right\vert ,\infty }d_{k}(\mu _{f},\mu _{g})\leq
C\left\Vert a\right\Vert _{d+1}2^{n(k+\left\vert \alpha \right\vert )}\beta
_{\mathbf{e}}(2^{nd})\left\Vert v\right\Vert _{\mathbf{e}_{\ast }}d_{k}(\mu
_{f},\mu _{g})
\end{eqnarray*}%
which implies (\ref{Oo5}). $\square $

\medskip

We are now ready for the

\smallskip

\textbf{Proof of Proposition \ref{8}.} Let $\alpha $ with $\left\vert \alpha
\right\vert \leq q.$ Using Proposition \ref{6}

\begin{equation*}
\partial _{\alpha }f=\sum_{n=1}^{\infty }\bar{\mathcal{H}}_{n}^{a}\diamond
\partial _{\alpha }f=\sum_{n=1}^{\infty }\bar{\mathcal{H}}_{n}^{a}\diamond
\partial _{\alpha }(f-f_{n})+\sum_{n=1}^{\infty }\bar{\mathcal{H}}%
_{n}^{a}\diamond \partial _{\alpha }f_{n}
\end{equation*}%
and using (\ref{Oo5}) and (\ref{Oo4})%
\begin{eqnarray*}
\left\Vert \partial _{\alpha }f\right\Vert _{\mathbf{e}} &\leq
&\sum_{n=1}^{\infty }\left\Vert \bar{\mathcal{H}}_{n}^{a}\diamond \partial
_{\alpha }(f-f_{n})\right\Vert _{\mathbf{e}}+\sum_{n=1}^{\infty }\left\Vert
\bar{\mathcal{H}}_{n}^{a}\diamond \partial _{\alpha }f_{n}\right\Vert _{%
\mathbf{e}} \\
&\leq &C\sum_{n=1}^{\infty }2^{n(\left\vert \alpha \right\vert +k)}\beta _{%
\mathbf{e}}(2^{nd})d_{k}(\mu _{f},\mu _{f_{n}})+C\sum_{n=1}^{\infty }\frac{1%
}{2^{2nm}}\left\Vert f_{n}\right\Vert _{2m+|\alpha |,2m,\mathbf{e}}
\end{eqnarray*}%
so (\ref{Oo10}) is proved. $\square $

%%%%%%%%%%%%%%%%%%%%%%%%%%%%%%%%%%%%%%%%%%%%%%%%%%%%%%%%%%%%%%%%%%%

\section{Interpolation spaces}

\label{sect-interp}

In this section we prove that, in the case of the $L^{p}$ norms, (that is ${%
\mathbf{e}}={\mathbf{e}}_p)$ the space $\mathcal{S}_{q,k,m,{\mathbf{e}}_p}$
is an interpolation space between $W_{\ast }^{k,\infty }$ (the dual of $%
W^{k,\infty })$\ and $W^{q,2m,p}.$ A similar interpretation holds for ${%
\mathbf{e}}_{\log }$ but this case is more exotic and we do not enter into
details here.

To begin we recall the framework of interpolation spaces. We are given two
Banach spaces $(X,\left\Vert \cdot \right\Vert _{X})$ and $(Y,\left\Vert
\cdot \right\Vert _{Y})$ with $X\subset Y$ (with continuous embedding). We
denote $\mathcal{L}(X,X)$ the space of the linear bounded operators from $X$
into itself and we denote by $\left\Vert L\right\Vert _{X,X}$ the operator
norm. A Banach space $(W,\left\Vert \cdot \right\Vert _{W})$ such that $%
X\subset W\subset Y$ is called an interpolation space for $X$ and $Y$ if $%
\mathcal{L}(X,X)\cap \mathcal{L}(Y,Y)\subset \mathcal{L}(W,W).$ Let $\gamma
\in (0,1).$ If there exists a constant $C$ such that $\left\Vert
L\right\Vert _{W,W}\leq C\left\Vert L\right\Vert _{X,X}^{\gamma }\left\Vert
L\right\Vert _{Y,Y}^{1-\gamma }$ for every $L\in \mathcal{L}(X,X)\cap
\mathcal{L}(Y,Y)$ then $W$ is an interpolation space of order $\gamma .$ And
if one may take $C=1$ then $W$ is an exact interpolation space of order $%
\gamma .$ There are several methods for constructing interpolation spaces.
We focus here on the so called $K$-method. For $y\in Y$ and $t>0$ one
defines $K(y,t)=\inf_{x\in X}(\left\Vert y-x\right\Vert _{Y}+t\left\Vert
x\right\Vert _{X})$ and
\begin{equation*}
\left\Vert y\right\Vert _{\gamma }=\int_{0}^{\infty }t^{-\gamma }K(y,t)\frac{%
dt}{t},\qquad (X,Y)_{\gamma }=\{y\in Y:\left\Vert y\right\Vert _{\gamma
}<\infty \}.
\end{equation*}%
Then one proves that $(X,Y)_{\gamma }$ is an exact interpolation space of
order $\gamma .$ One may also use the following discrete variant of the
above norm. Let $\gamma \geq 0.$ For $y\in Y$ and for a sequence $x_{n}\in
X,n\in {\mathbb{N}}$ we define
\begin{equation}
\pi _{\gamma }(y,(x_{n})_{n})=\sum_{n=1}^{\infty }2^{n\gamma }\left\Vert
y-x_{n}\right\Vert _{Y}+\frac{1}{2^{n}}\left\Vert x_{n}\right\Vert _{X}
\label{Int1}
\end{equation}%
and
\begin{equation*}
\rho _{\gamma }^{X,Y}(y)=\inf \pi _{\gamma }(y,(x_{n})_{n})
\end{equation*}%
with the infimum taken over all the sequences $x_{n}\in X,n\in {\mathbb{N}}.$
Then a standard result in interpolation theory (the proof is elementary)
says that there exists a constant $C>0$ such that
\begin{equation}
\frac{1}{C}\left\Vert y\right\Vert _{\gamma }\leq \rho _{\gamma
}^{X,Y}(y)\leq C\left\Vert y\right\Vert _{\gamma }  \label{Int1'}
\end{equation}%
so that%
\begin{equation*}
\mathcal{S}_{\gamma }(X,Y)=:\{y:\rho _{\gamma }^{X,Y}(y)<\infty
\}=(X,Y)_{\gamma }
\end{equation*}%
Take now $q,k\in N,m\in {\mathbb{N}}_{\ast }$ and $p>1$ and set $Y=W_{\ast
}^{k,\infty }$ and $X=W^{q,2m,p}.$ Then with the notation from (\ref{011})
and (\ref{O13})
\begin{equation}
\rho _{q,k,m,{\mathbf{e}}_{p}}(\mu )=\rho _{\gamma }^{X,Y}(\mu )\quad %
\mbox{and}\quad \mathcal{S}_{q,k,m,{\mathbf{e}}_{p}}=\mathcal{S}_{\gamma
}(X,Y),\quad \mbox{with}\quad \gamma =\frac{q+k+d/p_{\ast }}{2m}
\label{Int2}
\end{equation}%
Notice that in the definition of $\mathcal{S}_{q,k,m,{\mathbf{e}}_{p}}$ one
does not use precisely $\pi _{\gamma }(y,(x_{n})_{n})$ but $\pi _{\gamma
}^{(m)}(y,(x_{n})_{n})$ defined by%
\begin{eqnarray*}
\pi _{\gamma }^{(m)}(y,(x_{n})_{n}) &=&\sum_{n=1}^{\infty
}2^{n(q+k+d/p_{\ast })}\left\Vert y-x_{n}\right\Vert _{Y}+\frac{1}{2^{2mn}}%
\left\Vert x_{n}\right\Vert _{X} \\
&=&\sum_{n=1}^{\infty }2^{2mn\gamma }\left\Vert y-x_{n}\right\Vert _{Y}+%
\frac{1}{2^{2mn}}\left\Vert x_{n}\right\Vert _{X}
\end{eqnarray*}%
with $\gamma =\frac{q+k+d/p_{\ast }}{2m}.$ The fact that one uses $2^{2mn}$
instead of $2^{n}$ has no impact except that it changes the constants in (%
\ref{Int1'}). So the spaces are the same.

\bigskip

We turn now to a different point. For $p>1$ and $0<s<1$ we denote by $%
\mathcal{B}^{s,p}$ the Besov space and by $\left\Vert f\right\Vert _{%
\mathcal{B}^{s,p}}$ the Besov norm (see Triebel \cite{Tbis} for definitions
and notations). Our aim is to give a criterion which guarantees that a
function $f$ belongs to $\mathcal{B}^{s,p}.$ We will use the classical
equality $(W^{1,p},L^{p})_{s}=\mathcal{B}^{s,p}$.

\begin{lemma}
\label{Besov} Let $p>1$ and $0<s^{\prime }<s<1.$ Consider a function $\phi
\in C^{\infty }$ such that $\int_{{\mathbb{R}}^d}\phi(x)dx=1$ and let $\phi
_{\delta }(x)=\frac{1}{\delta ^{d}}\phi (\frac{x}{\delta })$ and $\phi
_{\delta }^{i}(x)=x^{i}\phi _{\delta }(x).$ We assume that $f\in L^{p}$
verifies the following hypothesis: for every $i=1,...,d$
\begin{equation}  \label{Int3}
\begin{aligned} i)&\qquad \limsup_{\delta \rightarrow 0}\delta
^{1-s}\left\Vert \partial _{i}(f\ast \phi _{\delta })\right\Vert _{p}
<\infty \\ ii)&\qquad \limsup_{\delta \rightarrow 0}\delta ^{-s}\left\Vert
\partial _{i}(f\ast \phi _{\delta }^{i})\right\Vert _{p} <\infty .
\end{aligned}
\end{equation}%
Then $f\in \mathcal{B}^{s^{\prime },p}$ for every $s^{\prime }<s$.
\end{lemma}

\textbf{Proof}. Let $f\in C^{1}.$ We use a Taylor expansion of order one and
we obtain%
\begin{eqnarray*}
f(x)-f\ast \phi _{\varepsilon }(x) &=&\int (f(x)-f(x-y))\phi _{\varepsilon
}(y)dy=\int_{0}^{1}d\lambda \int \left\langle \nabla f(x-\lambda
y),y\right\rangle \phi _{\varepsilon }(y)dy \\
&=&\int_{0}^{1}d\lambda \int \left\langle \nabla f(x-z),z\right\rangle \frac{%
1}{\lambda }\phi _{\varepsilon }\Big(\frac{z}{\lambda }\Big)\frac{dz}{%
\lambda ^{d}}=\int_{0}^{1}d\lambda \int \left\langle \nabla
f(x-z),z\right\rangle \phi _{\varepsilon \lambda }(z)\frac{dz}{\lambda } \\
&=&\sum_{i=1}^{d}\int_{0}^{1}\partial _{i}(f\ast \phi _{\varepsilon \lambda
}^{i})(x)\frac{d\lambda }{\lambda }.
\end{eqnarray*}

It follows that
\begin{equation*}
\left\Vert f-f\ast \phi _{\varepsilon }\right\Vert _{p}\leq
\sum_{i=1}^{d}\int_{0}^{1}\left\Vert \partial _{i}(f\ast \phi _{\varepsilon
\lambda }^{i})\right\Vert _{p}\frac{d\lambda }{\lambda }\leq d\,\varepsilon
^{s}\int_{0}^{1}\frac{d\lambda }{\lambda ^{1-s}}=C\varepsilon ^{s}.
\end{equation*}%
We also have $\left\Vert f\ast \phi _{\varepsilon }\right\Vert
_{W^{1,p}}\leq C(1+\left\Vert f\right\Vert _{\infty })\varepsilon ^{-(1-s)}$
so that
\begin{equation*}
K(f,\varepsilon )\leq \left\Vert f-f\ast \phi _{\varepsilon }\right\Vert
_{p}+\varepsilon \left\Vert f\ast \phi _{\varepsilon }\right\Vert
_{W^{1,p}}\leq C\varepsilon ^{s}.
\end{equation*}%
We conclude that for $s^{\prime }<s$ we have%
\begin{equation*}
\int_{0}^{1}\frac{1}{\varepsilon ^{s^{\prime }}}K(f,\varepsilon )\frac{%
d\varepsilon }{\varepsilon }\leq C\int_{0}^{1}\frac{\varepsilon ^{s}}{%
\varepsilon ^{s^{\prime }}}\frac{d\varepsilon }{\varepsilon }<\infty
\end{equation*}%
so $f\in (W^{1,p},L^{p})_{s^{\prime }}=\mathcal{B}^{s^{\prime },p}.\ \square
$

\section{Super kernels}

\label{app-superkernels}

A super kernel $\phi :{\mathbb{R}}^{d}\rightarrow {\mathbb{R}}$ is a
function which belongs to the Schwartz space $\mathcal{S}$\ (infinitely
differentiable functions which decrease in a polynomial way to infinity), $%
\int \phi (x)dx=1,$ and such that for every non null multi index $\alpha
=(\alpha _{1},...,\alpha _{d})\in {\mathbb{N}}^{d}$ one has
\begin{equation}
\int y^{\alpha }\phi (y)dy=0\qquad y^{\alpha }=\prod_{i=1}^{d}y_{i}^{\alpha
_{i}}.  \label{kk1}
\end{equation}%
See \cite{[KK]} Section 3, Remark 1 for the construction of a superkernel.\
The corresponding $\phi _{\delta }$, $\delta \in (0,1)$, is defined by
\begin{equation*}
\phi _{\delta }(y)=\frac{1}{\delta ^{d}}\phi \Big(\frac{y}{\delta }\Big).
\end{equation*}%
For a function $f$ we denote $f_{\delta }=f\ast \phi _{\delta }.$ We will
work with the norms $\left\Vert f\right\Vert _{k,\infty }$ and $\left\Vert
f\right\Vert _{q,l,{\mathbf{e}}}$ defined in (\ref{O4}) and in (\ref{O4a}).
And we have

\begin{lemma}
\label{kernel copy(1)} i) Let $k,q\in {\mathbb{N}},l>d$ and ${\mathbf{e}}\in%
\mathcal{E}$. There exists a universal constant $C$ such that for every $%
f\in W^{q,l,{\mathbf{e}}}$ one has%
\begin{equation}
\left\Vert f-f_{\delta }\right\Vert _{W_{\ast }^{k,\infty }}\leq C\left\Vert
f\right\Vert _{q,l,{\mathbf{e}}}\delta ^{q+k}.  \label{kk2}
\end{equation}%
ii) Let $l>d,n,q\in {\mathbb{N}}$, with $n\geq q$, and ${\mathbf{e}}\in%
\mathcal{E}$. There exists a universal constant $C$ such that%
\begin{equation}
\left\Vert f_{\delta }\right\Vert _{n,l,p}\leq C\left\Vert f\right\Vert
_{q,l,{\mathbf{e}}}\delta ^{-(n-q)}.  \label{kk3}
\end{equation}
\end{lemma}

\textbf{Proof}. i) We may suppose without loss of generality that $f\in
C_{b}^{\infty }.$ Using Taylor expansion of order $q+k$
\begin{eqnarray*}
f(\emph{x)}-f_{\delta }(x) &=&\int (f(\emph{x)}-f(y))\phi _{\delta }(x-y)dy
\\
&=&\int I(x,y)\phi _{\delta }(x-y)dy+\int R(x,y)\phi _{\delta }(x-y)dy
\end{eqnarray*}%
with
\begin{eqnarray*}
I(x,y) &=&\sum_{i=1}^{q+k-1}\frac{1}{i!}\sum_{\left\vert \alpha \right\vert
=i}\partial ^{\alpha }f(x)(x-y)^{\alpha }, \\
R(x,y) &=&\frac{1}{(q+k)!}\sum_{\left\vert \alpha \right\vert
=q+k}\int_{0}^{1}\partial ^{\alpha }f(x+\lambda (y-x))(x-y)^{\alpha
}d\lambda .
\end{eqnarray*}%
Using (\ref{kk1}) we obtain $\int I(x,y)\phi _{\delta }(x-y)dy=0$ and by a
change of variable we get
\begin{equation*}
\int R(x,y)\phi _{\delta }(x-y)dy=\frac{1}{(q+k)!}\sum_{\left\vert \alpha
\right\vert =q+k}\int_{0}^{1}\int dz\phi _{\delta }(z)\partial ^{\alpha
}f(x+\lambda z)z^{\alpha }d\lambda .
\end{equation*}%
We consider now $g\in W^{k,\infty }$ and we write%
\begin{equation*}
\int (f(\emph{x)}-f_{\delta }(x))g(x)dx=\frac{1}{(q+k)!}\sum_{\left\vert
\alpha \right\vert =q+k}\int_{0}^{1}d\lambda \int dz\phi _{\delta
}(z)z^{\alpha }\int \partial ^{\alpha }f(x+\lambda z)g(x)dx.
\end{equation*}%
Let us denote $f_{a}(x)=f(x+a).$ We have $(\partial ^{\alpha
}f)(x+a)=(\partial ^{\alpha }f_{a})(x).$ Let $\alpha $ with $\left\vert
\alpha \right\vert =\sum_{i=1}^{d}\alpha _{i}=q+k.$ We split $\alpha $ into
two multi indexes $\beta $ and $\gamma $ such that $\left\vert \beta
\right\vert =k,\left\vert \gamma \right\vert =q$ and $\partial ^{\beta
}\partial ^{\gamma }=\partial ^{\alpha }$ (this may be done in several ways
but any one of them is good for us). Then using integration by parts%
\begin{eqnarray*}
\left\vert \int \partial ^{\alpha }f(x+\lambda z)g(x)dx\right\vert
&=&\left\vert \int \partial ^{\beta }\partial ^{\gamma }f_{\lambda
z}(x)g(x)dx\right\vert \\
&\leq &\int \left\vert \partial ^{\gamma }f_{\lambda z}(x)\right\vert
\left\vert \partial ^{\beta }g(x)\right\vert dx\leq \left\Vert g\right\Vert
_{k,\infty }\int \left\vert \partial ^{\gamma }f_{\lambda z}(x)\right\vert dx
\\
&=&\left\Vert g\right\Vert _{k,\infty }\int \left\vert \partial ^{\gamma
}f(x)\right\vert dx.
\end{eqnarray*}%
We write $\partial ^{\gamma }f(x)=u_{l}(x)v_{\gamma }(x)$ with $%
u_{l}(x)=(1+\left\vert x\right\vert ^{2})^{-l/2}$ and $v_{\gamma
}(x)=(1+\left\vert x\right\vert ^{2})^{l/2}\partial ^{\gamma }f(x).$ Using H%
\"{o}lder inequality%
\begin{equation*}
\int \left\vert \partial ^{\gamma }f(x)\right\vert dx\leq C\left\Vert
u_{l}\right\Vert _{{\mathbf{e}}_{\ast }}\left\Vert v_{\gamma }\right\Vert _{{%
\mathbf{e}}}\leq C\left\Vert u_{l}\right\Vert _{{\mathbf{e}}_{\ast
}}\left\Vert f\right\Vert _{q,l,{\mathbf{e}}}.
\end{equation*}%
By Remark \ref{U} $\left\Vert u_{l}\right\Vert _{{\mathbf{e}}_{\ast
}}<\infty .$ So we obtain
\begin{eqnarray*}
\left\vert \int_{0}^{1}\int dz\phi _{\delta }(z)z^{\alpha }\int \partial
^{\alpha }f(x+\lambda z)g(x)dxd\lambda \right\vert &\leq &C\left\Vert
f\right\Vert _{q,l,{\mathbf{e}}}\left\Vert g\right\Vert _{k,\infty }\int
\phi _{\delta }(z)\left\vert z\right\vert ^{k+q}dz \\
&\leq &C\left\Vert f\right\Vert _{q,l,{\mathbf{e}}}\left\Vert g\right\Vert
_{k,\infty }\delta ^{k+q}.
\end{eqnarray*}

ii) Let $\alpha $ be a multi index with $\left\vert \alpha \right\vert =n$
and let $\beta ,\gamma $ be a splitting of $\alpha $ with $\left\vert \beta
\right\vert =q$ and $\left\vert \gamma \right\vert =n-q.$ Using the triangle
inequality, for every $y$ we have $1+\left\vert x\right\vert \leq
(1+\left\vert y\right\vert )(1+\left\vert x-y\right\vert ).$ Then
\begin{align*}
u(x)& :=(1+\left\vert x\right\vert )^{l}\left\vert \partial ^{\alpha
}f_{\delta }(x)\right\vert =(1+\left\vert x\right\vert )^{l}\left\vert
\partial ^{\beta }f\ast \partial ^{\gamma }\phi _{\delta }(x)\right\vert \\
& \leq \int (1+\left\vert x\right\vert )^{l}\left\vert \partial ^{\beta
}f(y)\right\vert \left\vert \partial ^{\gamma }\phi _{\delta
}(x-y)\right\vert dy\leq \alpha \ast \beta (x)
\end{align*}%
with%
\begin{equation*}
\alpha (y)=(1+\left\vert y\right\vert )^{l}\left\vert \partial ^{\beta
}f(y)\right\vert ,\qquad \beta (z)=(1+\left\vert z\right\vert
)^{l}\left\vert \partial ^{\gamma }\phi _{\delta }(z)\right\vert .
\end{equation*}%
Using (\ref{Oo2})\ we obtain
\begin{equation*}
\left\Vert u\right\Vert _{{\mathbf{e}}}\leq \left\Vert \alpha \ast \beta
\right\Vert _{{\mathbf{e}}}\leq \left\Vert \beta \right\Vert _{1}\left\Vert
\alpha \right\Vert _{{\mathbf{e}}}\leq \frac{C}{\delta ^{n-q}}\left\Vert
\alpha \right\Vert _{{\mathbf{e}}}=\frac{C}{\delta ^{n-q}}\left\Vert
f_{\beta ,l}\right\Vert _{{\mathbf{e}}}.
\end{equation*}%
$\square $


\addcontentsline{toc}{section}{References}

\begin{thebibliography}{99}
\bibitem{bib:[BCa]} \textsc{V. Bally, L. Caramellino} (2011). Riesz
transform and integration by parts formulas for random variables. \emph{%
Stochastic Process. Appl.}, \textbf{121}, 1332-1355.

\bibitem{BCb} \textsc{V. Bally, L. Caramellino} (2014). On the distance
between probability density functions. Preprint \texttt{arXiv:1311.7555}.

\bibitem{bib:BCpreprint} \textsc{V. Bally, L. Caramellino} (2012).
Regularity of probability laws by using an interpolation method. Preprint
\texttt{arXiv:1211.0052}.

\bibitem{bib:BC-Horm} \textsc{V. Bally, L. Caramellino} (2013). Regularity
of Wiener functionals under an H\"{o}rmander type condition of order one.
Preprint \texttt{arXiv:1307.3942}.

\bibitem{bib:[BCl1]} \textsc{V. Bally, E. Cl\'{e}ment} (2011). Integration
by parts formulas and applications to equations with jumps. \emph{Probab.
Theory Related Fields} \textbf{151}, 613-657.

\bibitem{bib:[BCl2]} \textsc{V. Bally, E. Cl\'{e}ment} (2011). Integration
by parts formulas with respect to jump times and stochastic differential
equations. \emph{Stochastic Analysis 2010}, ed. Dan Crisan, Springer.

\bibitem{bib:[BF]} \textsc{V. Bally, N. Fournier} (2011). Regularization
properties od the 2D homogeneous Boltzmann equation without cutoff. \emph{%
Probab. Theory Related Fields} \textbf{151}, 659-704.

\bibitem{bib:[BP]} \textsc{V. Bally, E. Pardoux} (1998). Malliavin Calculus
for White Noise Driven Parabolic SPDE's. \emph{Potential Anal.} \textbf{9},
27-64.

\bibitem{bib:[BS]} \textsc{C. Bennett, R. Sharpley} (1988). Interpolation of
operators. Academic Press INC.

\bibitem{DF} \textsc{A. Debussche, N. Fournier} (2013). Existence of
densities for stable-like driven SDE's with Holder continuous coefficients.
\emph{J. Funct. Anal.} \textbf{264}, 1757-1778.

\bibitem{DR} \textsc{A. Debussche, M. Romito} (2014). Existence of densities
for the 3D Navier--Stokes equations driven by Gaussian noise. \emph{Probab.
Theory Related Fields} \textbf{158}, 575-596.

\bibitem{bib:[De]} \textsc{S. De Marco} (2011). Smoothness and Asymptotic
Estimates of densities for SDEs with locally smooth coefficients and
Applications to square-root diffusions. \emph{Ann. Appl. Probab.} \textbf{21}%
, 1282-1321.

\bibitem{bib:[Dz]} \textsc{J. Dziubanski} (1997). Triebel-Lizorkin spaces
associated with Laguerre and Hermite expansions. \emph{Proc. Amer. Math. Soc.%
} \textbf{125}, 3547-3554.

\bibitem{bib:[E]} \textsc{J. Epperson} (1985). Hermite and Laguerre wave
packet expansions. \emph{Studia Math. J.} \textbf{34}, 777-799.

\bibitem{Ff} \textsc{N. Fournier} (2002). Jumping SDE's: absolute continuity
using monotonocity. \emph{Stochastic Process. Appl.}, \textbf{98}, 317-330.

\bibitem{bib:[F]} \textsc{N. Fournier} (2008). Smoothness of the law of some
one-dimensional jumping SDE's with non constant rate of jump. \emph{%
Electron. J. Probab.} \textbf{13}, 135-156.

\bibitem{F1} \textsc{N. Fournier} (2012). Finiteness of entropy for the
homogeneous Boltzmann equation with measure initial condition. To appear on
\emph{Ann. Appl. Probab.}; \texttt{arXiv:1203.0130}.

\bibitem{bib:[FP]} \textsc{N. Fournier, J. Printems} (2010). Absolute
continuity of some one-dimensional processes. \emph{Bernoulli} \textbf{16},
343-360.

\bibitem{[KK]} \textsc{A. Kebaier, A. Kohatsu-Higa} (2008). An optimal
control variance reduction method for density estimation. \emph{Stochastic
Process. Appl.} \textbf{118}, 2143--2180.

\bibitem{bib:[K.M-W]} \textsc{P. Kosmol, D. M\^{u}ller-Wichards} (2011).
\emph{Optimization in Functional Spaces with stability considerations in
Orlicz spaces.} De Gruyter Series in Non Linear Analysis and Applications
13. Walter de Gruyter.

\bibitem{[IW]} \textsc{N. Ikeda, S. Watanabe} (1989). \emph{Stochastic
Differential Equations and Diffusion processes}. North-Holland Mathematical
Library 24.

\bibitem{EL} \textsc{E. L\"ocherbach, D. Loukianova, O. Loukianov} (2011).
Polynomial bounds in the Ergodic Theorem for positive reccurent
one-dimensional diffusions and integrability of hitting times. \emph{Ann.
Inst. Henri Poincar\'e Probab. Stat.}, \textbf{47}, 425-449.

\bibitem{bib:[M]} \textsc{P. Malliavin} (1997). \emph{Stochastic Analysis}.
Springer.

\bibitem{bib:[VN]} \textsc{S. Ninomiya, N. Victoir} (2008). Weak
approximation of stochastic differential equations and applications to
derivative pricing. \emph{Appl. Math. Finance}, \textbf{15}, 107-121.

\bibitem{bib:NP} \textsc{I. Nourdin, G. Peccati} (2012). \emph{Normal
approximations with Malliavin calculus. From Stein's method to universality}%
. Cambridge Tracts in Mathematics, 192. Cambridge University Press.

\bibitem{bib:[N]} \textsc{D. Nualart} (2006) \emph{The Malliavin calculus
and related topics. Second Edition}. Springer-Verlag.

\bibitem{bib:[PY]} \textsc{P. Petrushev, Yuan Xu} (2008). Decomposition of
spaces of distributions induced by Hermite expansions. \emph{J. Fourier
Anal. and Appl.} \textbf{14}, 372-414.

\bibitem{bib:[PZ]} \textsc{E. Pardoux, T. Zhang} (1993). Absolute continuity
for the law of the solution of a parabolic SPDE. \emph{J. Funct. Anal.}
\textbf{112}, 447-458.

\bibitem{bib:[Sa]} \textsc{M. Sanz Sol\`e} (2005). \emph{Malliavin Calculus,
with Applications to Stochastic Partial Differential Equations}. EPFL Press.
Fundamental Sciences, Mathematics.

\bibitem{bib:[TT]} \textsc{D. Talay, L. Tubaro} (1990). Expansion of the
global error for numerical schemes solving stochastic differential
equations. \emph{Stochastic Anal. Appl.}, \textbf{8}, 94-120.

\bibitem{T} \textsc{H. Triebel} (1999). \emph{Interpolation Theory -
Function Spaces - Differential Operators}. John Wiley \& Sons, Incorporated.

\bibitem{Tbis} \textsc{H. Triebel} (2006). \emph{Theory of function spaces
III}. Birkh\"auser Verlag.

\bibitem{bib:[W]} \textsc{J.B. Walsh} (1986). \emph{An introduction to
stochastic partial differential equations}. In Ecole d'Et\'{e} de Probabilit%
\'{e}s de Saint Flour XV, Lecture Notes in Math. 1180, Springer; pp 226-437.

\bibitem{V1} \textsc{A. Yu. Veretenikov} (1997). On polynomial mixing bounds
for stochastic differential equations, \emph{Stochastic Process. Appl.}
\textbf{70}, 115-127.

\bibitem{V2} \textsc{A. Yu. Veretenikov, S.A. Klokov} (2004). Subexponential
mixing rate for a class of Markov processes (multidimensional case). \emph{%
Theory Probab. Appl.} \textbf{49}, 1-13 (translation from Russian of \emph{%
Teoriya Veroyatn. Primen.} \textbf{48}, 21-35, 2003).
\end{thebibliography}
\end{document}